\documentclass{article}
\usepackage[dvipsnames,svgnames,x11names]{xcolor}
\usepackage{arxiv}
\usepackage{amsmath}
\numberwithin{equation}{section}
\usepackage[utf8]{inputenc} 
\usepackage{hyperref}       
\usepackage{url}            
\usepackage{booktabs}       
\usepackage{amsfonts}       
\usepackage{nicefrac}       
\usepackage{microtype}      
\usepackage{cleveref}       
\usepackage{lipsum}         
\usepackage{graphicx}
\usepackage{natbib}
\usepackage{doi}
\usepackage{tikz}
\usepackage{amsthm}
\usepackage{esint}
\tikzset{>=stealth}
\usepackage{bm}
\usepackage{yhmath}
\usetikzlibrary{patterns}
\usepackage{amssymb}

\title{Reasonable mechanical model on shallow tunnel excavation to eliminate displacement singularity caused by unbalanced resultant}

\date{\today}
\author{
  \href{https://orcid.org/0000-0002-5143-2714}
  {\includegraphics[scale=0.06]{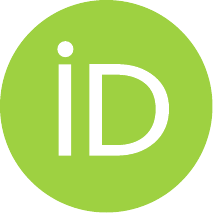}
    \bf {Luobin Lin}} \\
  Fujian Provincial Key Laboratory of Advanced Technology and Informatization in Civil Engineering\\
  College of Civil Engineering\\
  Fujian University of Technology\\
  No. 69 Xueyuan Road, Shangjie University Town, Fuzhou, 350118, Fujian, China \\
  \texttt{luobin\_lin@fjut.edu.cn} \\
  \href{https://orcid.org/0000-0002-5583-3734}
  {\includegraphics[scale=0.06]{orcid.pdf}
    \bf {Fuquan Chen}} \\
  College of Civil Engineering\\
  Fuzhou University\\
  No. 2 Xueyuan Road, Shangjie University Town, Fuzhou, 350108, Fujian, China \\
  \texttt{phdchen@fzu.edu.cn}\\
  \href{https://orcid.org/0000-0003-0679-9127}
  {\includegraphics[scale=0.06]{orcid.pdf}
    \bf {Xianhai Huang}}\\
  College of Civil Engineering\\
  Fujian University of Technology\\
  No. 69 Xueyuan Road, Shangjie University Town, Fuzhou, 350118, Fujian, China \\
  \texttt{hwhsh@163.com} \\
}

\renewcommand{\shorttitle}{Reasonable mechanical model on shallow tunnel excavation}

\hypersetup{
  pdftitle={Reasonable mechanical model on shallow tunnel excavation to eliminate displacement singularity caused by unbalanced resultant},
  pdfauthor={Luobin Lin, Fuquan Chen, Xianhai Huang},
  pdfkeywords={Shallow tunnel excavation, Unbalanced resultant, Displacement singularity, Analytic continuation, Riemann-Hilbert problem},
}

\begin{document}
\maketitle

\begin{abstract}
  When considering initial stress field in geomaterial, nonzero resultant of shallow tunnel excavation exists, which produces logarithmic items in complex potentials, and would further lead to a unique displacement singularity at infinity to violate geo-engineering fact in real world. The mechanical and mathematical reasons of such a unique displacement singularity in the existing mechanical models are elaborated, and a new mechanical model is subsequently proposed to eliminate this singularity by constraining far-field ground surface displacement, and the original unbalanced resultant problem is converted into an equilibrium one with mixed boundary conditions. To solve stress and displacement in the new model, the analytic continuation is applied to transform the mixed boundary conditions into a homogenerous Riemann-Hilbert problem with extra constraints, which is then solved using an approximate and iterative method with good numerical stability. The Lanczos filtering is applied to the stress and displacement solution to reduce the Gibbs phenomena caused by abrupt change of the boundary conditions along ground surface. Several numerical cases are conducted to verify the proposed mechanical model and the results strongly validate that the proposed mechanical model successfully eliminates the displacement singularity caused by unbalanced resultant with good convergence and accuracy to obtain stress and displacement for shallow tunnel excavation. A parametric investigation is subsequently conducted to study the influence of tunnel depth, lateral coefficient, and free surface range on stress and displacement distribution in geomaterial.
\end{abstract}

\keywords{Shallow tunnel excavation \and Unbalanced resultant \and Displacement singularity \and Analytic continuation \and Riemann-Hilbert problem}

\section{Introduction}
\label{sec:intro}

Shallow tunnel excavation is common in geo-engineering. Owing to initial stress field, the stress difference between top and bottom of a shallow tunnel always exists prior to excavation, and a subsequent excavation would alter stress and displacement distribution of geomaterial in a more mechanically complicated way than a deep tunnel. Moreover, the gravity of the excavated geomaterial can not be neglected comparing to the initial stress field, and the gravity gradient in geomaterial would also cause unbalanced resultant along tunnel periphery. 

An unbalanced resultant is mathematically challenging in the complex variable method, since it brings a logarithmic item in the complex potentials. Muskhelishvili~\cite{Muskhelishvili1966} has pointed out that the logarithmic item in complex potential may result in displacement singularity. Strack~\cite{Strack2002phdthesis} introduces an extra singularity in the upper half plane and the other corresponding logarithmic item into Muskhelishvili's complex potentials to meet the stress equilibrium in the lower half plane, but the displacement singularity at infinity caused by the unbalanced resultant along tunnel periphery still exists. With the modified complex potentials by Strack~\cite[]{Strack2002phdthesis} and Verruijt's conformal mapping~\cite[]{Verruijt1997traction,Verruijt1997displacement}, Verruijt and Strack~\cite[]{Strack_Verruijt2002buoyancy} and Lu et al.~\cite[]{Lu2016} propose solutions of different boundary conditions along tunnel periphery with unbalanced resultant of shallow tunnel excavation. Lu et al.~\cite[]{Lu2019new_solution} introduces the analytic continuation for free traction boundary along the whole gound surface to simplify the solution procedure and studies the stress and displacement in geomaterial for a shallow tunnel excavation subjected to underground water. Fang et al.~\cite[]{fang15:_compl} studies the stress and displacement of a underwater shallow tunnel excavation. Zeng et al.~\cite[]{Zengguisen2019} proposes a new conformal mapping for a noncircular cavity in a lower half plane, which is an extension of Verruijt's conformal mapping, and studies the stress and displacement in geomaterial due a noncircular shallow tunnel excavation. However, the displacement singularity at infinity is still not eliminated. The static equilibrium problems of shallow tunnel are also studied \cite[]{zhang2018complex,Ludechun2019,kong2021analytical}. The displacement singularity is not only observed in complex variable method, but also in the stress function method of real domain~\cite{Timoshenko_1951_Elasticity}, where a similar logarithmic item also exists in the stress potenial.

Unbalanced resultant may also be induced by surcharge loads along ground surface. The classic Flamant solution~\cite[]{flamant92:_sur} has shown that an unbalanced resultant on ground surface would cause displacement singularity at infinity as well. Wang et al.~\cite[]{wanghuaning2017flamant-english} propose a reasonable displacement model based on the Flamant's solution via symmetrical modification, and good results are obtained. Wang et al.~\cite[]{Wanghuaning2021shallow_rigid_lining,Wang2018surcharge_twin,Wang2018shallow_surcharge} then propose further studies on surcharge loads acting on ground surface of shallow tunnel, but no excavation process is considered, based on the symmetrical modification~\cite[]{wanghuaning2017flamant-english}. Zeng et al.~\cite[]{zeng2022analytical} also studies the visco-elastic mechanical behaviour of noncircular shallow tunnel on the symmetrical modification~\cite[]{wanghuaning2017flamant-english}. Zhang et al.~\cite[]{zhang2021analytical} study the visco-elastic deformation in visco-elastic geomaterial when surcharge loads is subjected on the ground surface. Inspired by Ref~\cite[]{wanghuaning2017flamant-english}, Lin et al. \cite[]{Self2020JEM} propose an extended displacement model of symmetrical modification for both surcharge load on ground surface and unbalanced resultant along shallow tunnel periphery for excavation, which indeed eliminates the displacement singularity, but displacement solution is highly dependent on the modification depth. In other words, the convergence of displacement is not as good as expected. Lu et al.~\cite[]{lu2021reasonable} modify the the coefficients of the Taylor expansion of the logarithmic item to obtain reasonable displacement distribution when a surcharge load of Gaussian distribution is applied on the ground surface. However, the solution~\cite[]{lu2021reasonable} is dependent on the constant in conformal mapping, which is corresponding to the complex coordinate of the mapping point in the physical plane of the unit disk origin. The value of this constant can be arbitrarily chosen, since a lower half plane without any cavity is simply connected. The solution by Lu et al.~\cite[]{lu2021reasonable} is unsurprisingly not suitable for unbalanced resultant problem of shallow tunnel excavation, since that constant should be unique for the remaining geomaterial after shallow tunnel excavation~\cite[]{Verruijt1997traction,Verruijt1997displacement}, which is a doubly connected region.

In summary, the displacement singularity caused by the logarithmic item in complex potentials is still not eliminated in the existing studies of shallow tunnel excavation, or the elimination is not as good as expected. For such reasons, we propose a new model in this study to confront the displacement singularity caused by unbalanced resultant of shallow tunnel excavation in a very straightforward and mechanical manner by introducing the analytic continuation and Riemann-Hilbert problem of mixed boundary value conditions into theoretical analyses in tunnel engineering.

\section{Problem formation}
\label{sec:problem}

\subsection{Existing mechanical model of shallow tunnel excavation}
\label{sec:problem-1}

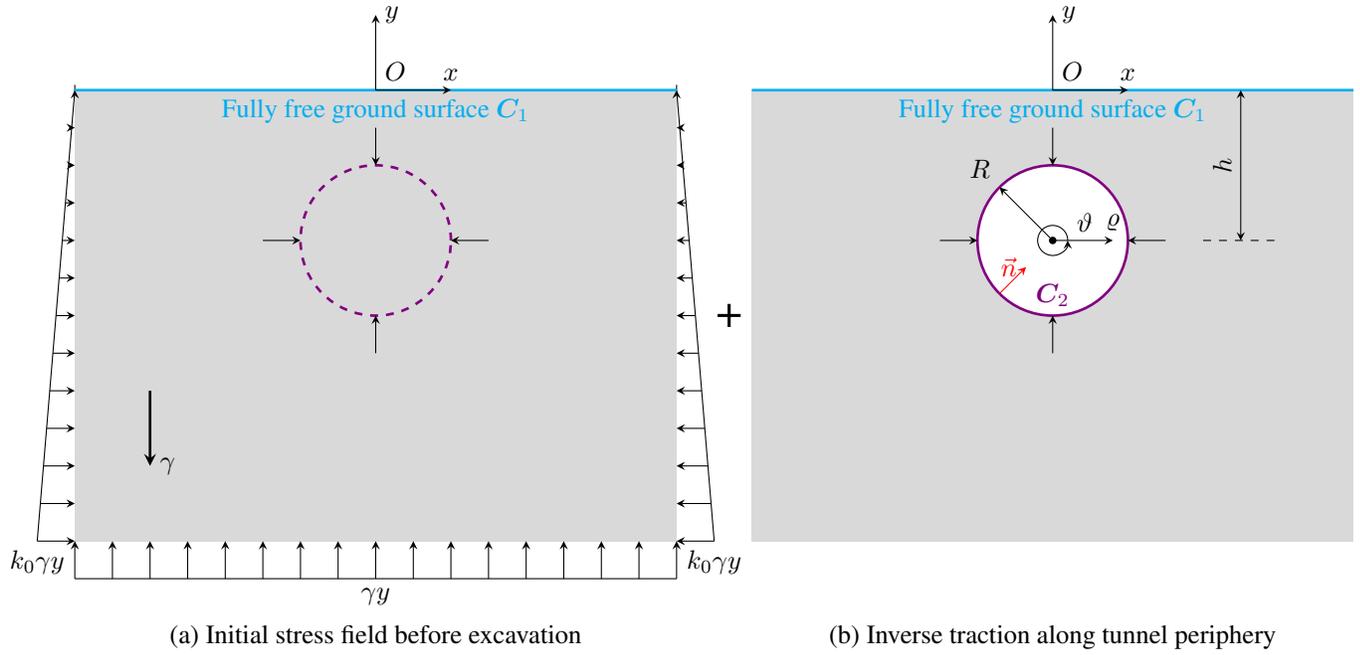
\begin{figure}[htb]
  \centering
  \begin{tikzpicture}
    \fill [gray!30] (0,0) rectangle (8,-6);
    \foreach \x in {0,1,2,...,16} \draw [->] ({\x*0.5},-6.5) -- ({\x*0.5},-6);
    \draw (0,-6.5) -- (8,-6.5);
    \node at (4,-6.5) [below] {$ \gamma y $};
    \foreach \x in {0,1,2,...,12} \draw [->] ({-\x*0.5/12},{-0.5*\x}) -- (0,{-0.5*\x});
    \draw (0,0) -- (-0.5,-6) node [below] {$ k_{0} \gamma y $};
    \foreach \x in {0,1,2,...,12} \draw [->] ({8+\x*0.5/12},{-0.5*\x}) -- (8,{-0.5*\x});
    \draw (8,0) -- (8.5,-6) node [below] {$ k_{0} \gamma y $};
    \draw [cyan, line width = 1pt] (0,0) -- (8,0);
    \draw [->, line width = 1pt] (1,-4) -- (1,-5) node [right] {$ \gamma $};
    \node at (4,0) [cyan, below] {Fully free ground surface $ {\bm C}_{1} $};
    \draw [->] (4,0) -- (5,0) node [above] {$ x $};
    \draw [->] (4,0) -- (4,1) node [right] {$ y $};
    \node at (4,0) [above right] {$ O $};
    \draw [violet, dashed, line width = 1pt] (4,-2) circle [radius = 1];
    \draw [->] (4,-0.5) -- (4,-1);
    \draw [->] (4,-3.5) -- (4,-3);
    \draw [->] (2.5,-2) -- (3,-2);
    \draw [->] (5.5,-2) -- (5,-2);
    \node at (4,-7) [below] {(a) Initial stress field before excavation};

    \fill [gray!30] (9,0) rectangle (17,-6);
    \draw [cyan, line width = 1pt] (9,0) -- (17,0);
    \node at (13,0) [cyan, below] {Fully free ground surface $ {\bm C}_{1} $};
    \draw [->] (13,0) -- (14,0) node [above] {$ x $};
    \draw [->] (13,0) -- (13,1) node [right] {$ y $};
    \node at (13,0) [above right] {$ O $};
    \fill [white] (13,-2) circle [radius = 1];
    \draw [violet, line width = 1pt] (13,-2) circle [radius = 1];
    \node at (13,-3) [violet, above] {$ {\bm C}_{2} $};
    \draw [->] (13,-0.5) -- (13,-1);
    \draw [->] (13,-3.5) -- (13,-3);
    \draw [->] (11.5,-2) -- (12,-2);
    \draw [->] (14.5,-2) -- (14,-2);
    \draw [->, red] ({13+cos(225)},{-2+sin(225)}) -- ({13+0.5*cos(225)},{-2+0.5*sin(225)}) node [left] {$ \vec{n} $};
    \draw [->] (13,-2) -- (13.8,-2) node [above] {$ \varrho $};
    \draw [->] (13.2,-2) arc [start angle = 0, end angle = 360, radius = 0.2];
    \fill (13,-2) circle [radius = 0.05];
    \node at (13.2,-2) [above right] {$ \vartheta $};
    \draw [dashed] (15,-2) -- (16,-2);
    \draw [<->] (15.5,0) -- (15.5,-2);
    \node at (15.5,-1) [rotate = 90, above] {$ h $};
    \draw [->] (13,-2) -- ({13+cos(135)},{-2+sin(135)}) node [above left] {$ R $};
    \node at (13,-7) [below] {(b) Inverse traction along tunnel periphery};

    \node at (8.7,-3) {\LARGE{+}};
  \end{tikzpicture}
  \caption{Initial stress field and inverse traction for shallow tunnel excavation}
  \label{fig:1}
\end{figure}

We start with a common model for shallow tunnelling in a heavy geomaterial in a lower half plane, similar to the cases in Refs~\cite{Strack2002phdthesis, Strack_Verruijt2002buoyancy, Lu2016, Zengguisen2019, Self2020JEM}. As shown in Fig.~\ref{fig:1}a, a geomaterial $ {\bm \varOmega} $ located in a complex lower half plane $ z(z = x + {\rm i}y) $ is linearly elastic, homogenerous, isotropic, of small deformation, and is subjected to a uniform intial stress field, which can be expressed as
\begin{equation}
  \label{eq:2.1}
  \left\{
    \begin{aligned}
      \sigma_{x}^{0} = & \; k_{0} \gamma y \\
      \sigma_{y}^{0} = & \; \gamma y \\
      \tau_{xy}^{0} = & \; 0
    \end{aligned}
  \right.
  , \quad y \leq 0
\end{equation}
where $ \sigma_{x}^{0} $, $ \sigma_{y}^{0} $, and $ \tau_{xy}^{0} $ denote horizontal, vertical, and shear stress components of the initial stress field in the rectangular coordinate system $ xOy $, $ k_{0} $ denotes lateral stress coefficient, $ \gamma $ denotes volumetric weight of geomaterial. The ground surface denoted by $ {\bm C}_{1} $ is fully free from any traction.

Then a shallow tunnel with radius of $R$ and buried depth of $h$ is excavated in the geomaterial of the lower half plane, and the tunnel periphery is denoted by $ {\bm C}_{2} $. The following external tractions are applied to tunnel periphery to cancel the tractions caused by the initial stress field as shown in Fig.~\ref{fig:1}b:
\begin{equation}
  \label{eq:2.2}
  \left\{
    \begin{aligned}
      X_{i}(S) = & \; - \sigma_{x}^{0}(S) \cdot \cos \langle \vec{n}, x \rangle - \tau_{xy}^{0}(S) \cdot \cos \langle \vec{n}, y \rangle \\
      Y_{i}(S) = & \; - \sigma_{y}^{0}(S) \cdot \cos \langle \vec{n}, y \rangle - \tau_{xy}^{0}(S) \cdot \cos \langle \vec{n}, x \rangle \\
    \end{aligned}
  \right.
  , \quad
  S: \; x^{2} + (y+h)^{2} = R^{2}
\end{equation}
where $ X_{i}(S) $ and $ Y_{i}(S) $ denote horizontal and vertical traction along tunnel periphery $ {\bm C}_{2} $ (denoted by $S$), respectively, $\vec{n}$ denotes the outward direction along the tunnel periphery, as shown in Fig.~\ref{fig:1}b, $ \langle \vec{n},x \rangle $ denotes angle between $ \vec{n} $ and $x$ axis, and $ \langle \vec{n},y \rangle $ denotes angle between $\vec{n}$ and $y$ axis. The initial stress field in Eq.~(\ref{eq:2.1}) and the tractions in Eq.~(\ref{eq:2.2}) together keep the tunnel periphery being traction-free. A local polar coordinate system $\varrho o \vartheta$ is introduced and located at the center of the tunnel, and Eq.~(\ref{eq:2.2}) turns to
\begin{equation}
  \label{eq:2.3}
  \left\{
    \begin{aligned}
      & X_{i}(S) = X_{i}(\vartheta) = k_{0} \gamma (-h+R\sin\vartheta) \cos\vartheta \\
      & Y_{i}(S) = Y_{i}(\vartheta) = \gamma (-h+R\sin\vartheta) \sin\vartheta \\
    \end{aligned}
  \right.
\end{equation}
Integrating the tractions along tunnel periphery in the clockwise direction (keeping the geomaterial on the left side) gives the following unbalanced resultants along tunnel periphery $ {\bm C}_{2} $:
\begin{equation}
  \label{eq:2.4}
  \left\{
    \begin{aligned}
      F_{x} = & \; \varointclockwise_{{\bm C}_{2}} X_{i}(S) |{\rm d}S| = - \varointclockwise_{{\bm C}_{2}} X_{i}(\vartheta) R {\rm d}\vartheta = 0 \\
      F_{y} = & \; \varointclockwise_{{\bm C}_{2}} Y_{i}(S) |{\rm d}S| = - \varointclockwise_{{\bm C}_{2}} Y_{i}(\vartheta) R {\rm d}\vartheta = \gamma \pi R^{2} \\
    \end{aligned}
  \right.
\end{equation}
where $ F_{x} $ and $ F_{y} $ denote horizontal and vertical components of the resultant acting along boundary ${\bm C}_{2}$, respectively; $ |{\rm d}S| = R|{\rm d}\vartheta| = - R{\rm d}\vartheta $ for clockwise length increment in physical plane. To be more specific, the value of the total gravity of the excavated geomaterial should be $-\gamma \pi R^{2}$, and the postive value of $F_{y}$ indicates the removal action of the excavated geomaterial and the corresponding external force applied to the tunnel periphery $ {\bm C}_{2} $.

\subsection{Unbalanced resultant and displacement singularity of existing mechanical model}
\label{sec:problem-2}

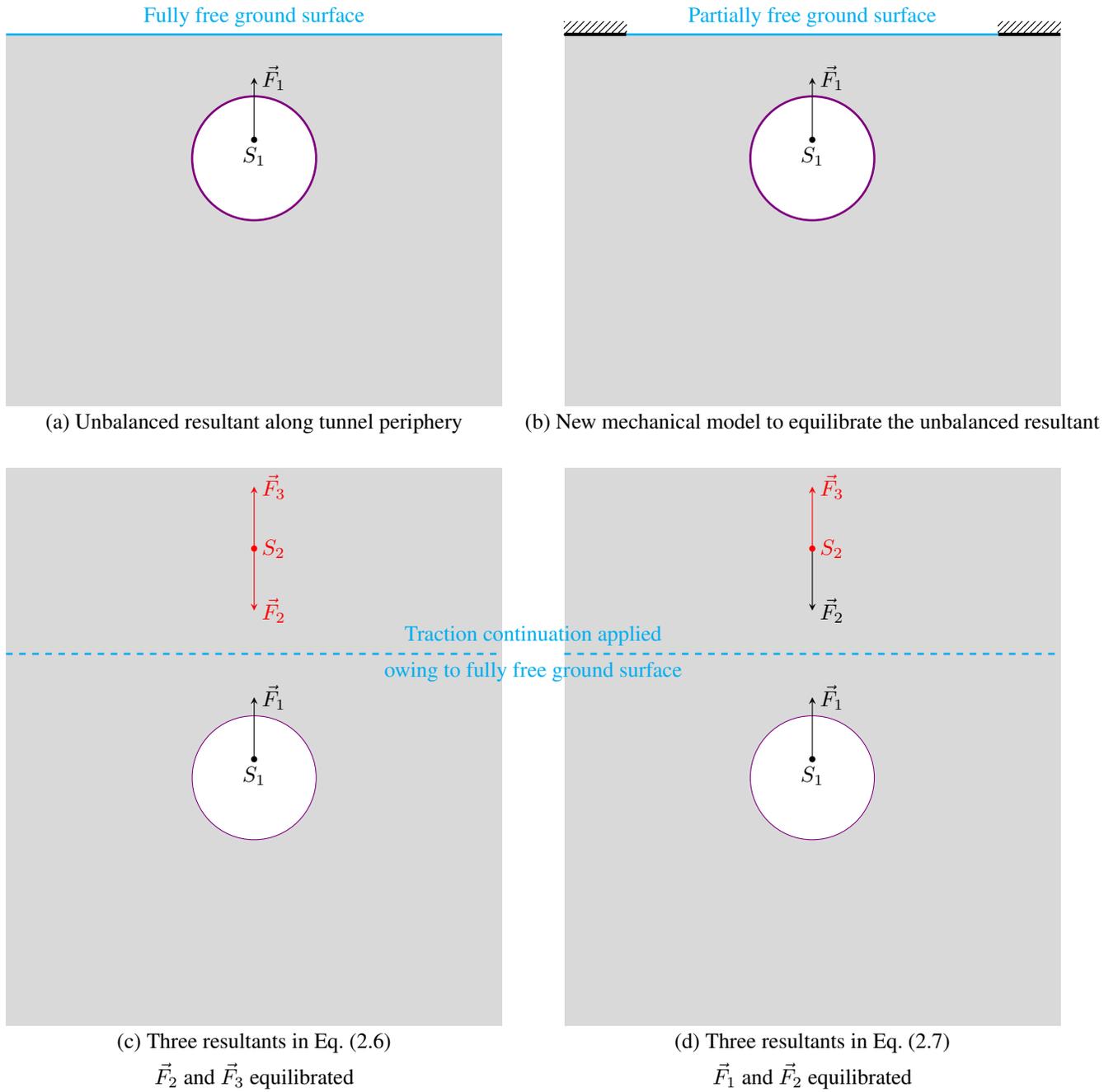
\begin{figure}[htb]
  \centering
  \begin{tikzpicture}
    \fill [gray!30] (0,0) rectangle (8,-6);
    \draw [cyan, line width = 1pt] (0,0) -- (8,0);
    \node at (4,0) [cyan, above] {Fully free ground surface};
    \fill [white] (4,-2) circle [radius = 1];
    \draw [violet, line width = 1pt] (4,-2) circle [radius = 1];
    \fill (4,-1.7) circle [radius = 0.05];
    \node at (4,-1.7) [below] {$ S_{1} $};
    \draw [->] (4,-1.7) -- (4,-0.7) node [right] {$ \vec{F}_{1} $};
    \node at (4,-6) [below] {(a) Unbalanced resultant along tunnel periphery};

    \fill [gray!30] (9,0) rectangle (17,-6);
    \draw [cyan, line width = 1pt] (10,0) -- (16,0);
    \draw [line width = 1.5pt] (9,0) -- (10,0);
    \draw [line width = 1.5pt] (16,0) -- (17,0);
    \fill [pattern = north east lines] (9,0) rectangle (10,0.2);
    \fill [pattern = north east lines] (16,0) rectangle (17,0.2);
    \node at (13,0) [cyan, above] {Partially free ground surface};
    \fill [white] (13,-2) circle [radius = 1];
    \draw [violet, line width = 1pt] (13,-2) circle [radius = 1];
    \fill (13,-1.7) circle [radius = 0.05];
    \node at (13,-1.7) [below] {$ S_{1} $};
    \draw [->] (13,-1.7) -- (13,-0.7) node [right] {$ \vec{F}_{1} $};
    \node at (13,-6) [below] {(b) New mechanical model to equilibrate the unbalanced resultant};

    \fill [gray!30] (0,-7) rectangle (8, -16);
    \draw [cyan, dashed, line width = 1pt] (0,-10) -- (8,-10);
    \fill [white] (4,-12) circle [radius = 1];
    \draw [violet] (4,-12) circle [radius = 1];
    \fill (4,-11.7) circle [radius = 0.05];
    \node at (4,-11.7) [below] {$ S_{1} $};
    \draw [->] (4,-11.7) -- (4,-10.7) node [right] {$ \vec{F}_{1} $};
    \draw [red, ->] (4,-8.3) -- (4,-9.3) node [right] {$ \vec{F}_{2} $};
    \draw [red, ->] (4,-8.3) -- (4,-7.3) node [right] {$ \vec{F}_{3} $};
    \fill [red] (4,-8.3) circle [radius = 0.05];
    \node at (4,-8.3) [right, red] {$ S_{2} $};
    \node at (4,-16) [below, align=center] {
      (c) Three resultants in Eq.~(\ref{eq:2.6})\\
      $ \vec{F}_{2} $ and $ \vec{F}_{3} $ equilibrated
    };

    \fill [gray!30] (9,-7) rectangle (17, -16);
    \draw [cyan, dashed, line width = 1pt] (9,-10) -- (17,-10);
    \fill [white] (13,-12) circle [radius = 1];
    \draw [violet] (13,-12) circle [radius = 1];
    \fill (13,-11.7) circle [radius = 0.05];
    \node at (13,-11.7) [below] {$ S_{1} $};
    \draw [->] (13,-11.7) -- (13,-10.7) node [right] {$ \vec{F}_{1} $};
    \draw [->] (13,-8.3) -- (13,-9.3) node [right] {$ \vec{F}_{2} $};
    \draw [red, ->] (13,-8.3) -- (13,-7.3) node [right] {$ \vec{F}_{3} $};
    \fill [red] (13,-8.3) circle [radius = 0.05];
    \node at (13,-8.3) [right, red] {$ S_{2} $};
    \node at (13,-16) [below, align = center] {
      (d) Three resultants in Eq.~(\ref{eq:2.7})\\
      $ \vec{F}_{1} $ and $ \vec{F}_{2} $ equilibrated
    };

    \node at (8.5,-10) [cyan, align = center] {
      Traction continuation applied\\
      owing to fully free ground surface
    };
  \end{tikzpicture}
  \caption{Unbalanced resultants and singularities in geomaterial and sketch of new mechanical model}
  \label{fig:2}
\end{figure}

With the unbalanced resultant in Eq.~(\ref {eq:2.4}), the excavation process is illustrated in Fig.~\ref{fig:2}a, where the unbalanced resultant denoted by $ \vec{F}_{1} (|\vec{F}_{1}| = F_{y}) $ is equivalently applied at singularity $ S_{1} $ within tunnel. Though such equivalence is not accurate enough, it helps to simplify our description below in a rough sense. Then the complex potentials $ \varphi (z) $ and $ \psi (z) $ related to the unbalanced resultant can be constructed according to Eqs. (4.1) and (4.2) in Ref~\cite{Strack2002phdthesis}, Eqs. (7) and (8) in Ref~\cite{Strack_Verruijt2002buoyancy}, Eqs. (17) and (18) in Ref~\cite{Lu2016}, or Eq.~(24) in Ref~\cite{Self2020JEM} as
\begin{subequations}
  \label {eq:2.5}
  \begin{equation}
    \label {eq:2.5a}
    \varphi(z) = -\frac{F_{x}+{\rm i}F_{y}}{2\pi(1+\kappa)} \left[ \kappa \ln (z-\overline{z}_{c}) + \ln(z-z_{c}) \right] + \varphi_{0}(z), \quad z \in {\bm \varOmega}
  \end{equation}
  \begin{equation}
    \label {eq:2.5b}
    \psi(z) = \frac{F_{x}-{\rm i}F_{y}}{2\pi(1+\kappa)} \left[ \ln(z-\overline{z}_{c}) + \kappa \ln(z-z_{c}) \right] + \psi_{0}(z), \quad z \in {\bm \varOmega}
  \end{equation}
\end{subequations}
where $ z_{c} $ denotes the complex coordinate of point $ S_{1} $ within the tunnel, $ \overline{z}_{c} $ denotes the complex coordinate of the corresponding conjugate point $ S_{2} $ in the upper half plane, $ \varphi_{0}(z) $ and $ \psi_{0}(z) $ denote the single-valued components of the complex potentials, which are always finite in the geomaterial.

With Eq.~(\ref{eq:2.4}), Eq.~(\ref{eq:2.5}) can be simplified and rewritten as
\begin{subequations}
  \label {eq:2.6}
  \begin{equation}
    \label {eq:2.6a}
    \varphi(z) = \left[ -\frac{{\rm i}F_{y}}{2\pi} \ln(z-\overline{z}_{c}) - \frac{-{\rm i} F_{y}}{2\pi(1+\kappa)} \ln (z-\overline{z}_{c}) \right] - \frac{{\rm i} F_{y}}{2\pi(1+\kappa)} \ln (z-z_{c}) + \varphi_{0} (z), \quad z \in {\bm \varOmega}
  \end{equation}
  \begin{equation}
    \label {eq:2.6b}
    \psi(z) = \left[ -\frac{{\rm i}F_{y}}{2\pi}\ln(z-\overline{z}_{c}) + \frac{{\rm i} \kappa F_{y}}{2\pi(1+\kappa)} \ln (z-\overline{z}_{c}) \right] + \frac{-{\rm i} \kappa F_{y}}{2\pi(1+\kappa)} \ln (z-z_{c}) + \psi_{0} (z), \quad z \in {\bm \varOmega}
  \end{equation}
\end{subequations}
where the first two logrithmic items take point $ S_{2} $ in the upper-half plane as singularity, and the last logarithmic item takes point $ S_{1} $ in the lower-half plane as singularity. Since the geomaterial has been assumed to be of small deformation, we can respectively examine the logarithmic items in Eq.~(\ref{eq:2.6}). Apparently, the second and third logarithmic items in Eq.~(\ref{eq:2.6}) indicate that an uppward and a downward concentrated resultant acting at singularities $ S_{2} $ and $ S_{1} $, respectively, according to $ \S $ 56a in Ref \cite[]{Muskhelishvili1966}, and the downward concentrated resultant is denoted by $ \vec{F}_{2} (|\vec{F}_{2}|=F_{y}) $ in Fig.~\ref{fig:2}c. Meanwhile, the first logarithmic item in Eq.~(\ref{eq:2.6}) indicates an upperward concentrated resultant acting at {\emph{some point on the boundary}} of the singularity $ S_{2} $ {\emph{after the existence of the singularity in the upper-half plane}}, according to $ \S $ 90 in Ref \cite[]{Muskhelishvili1966}, and such an upperward resultant is denoted by $ \vec{F}_{3}(|\vec{F}_{3}|=F_{y}) $ in Fig.~\ref{fig:2}c.

Since the singularity $ S_{2} $ is just a point, its boundary is itself. To guarantee the existence of its boundary, the only possible solution is that the upper-half plane should be fully defined except for the singularity $ S_{2} $, otherwise, the expression of the first logarithmic item in Eq.~(\ref{eq:2.6}) would be violated. In other words, the expression form of Eq.~(\ref{eq:2.6}) implicitly indicates that the analytic continuation principle has been applied {\emph{across the full length of the ground surface}}. Since the ground surface is always traction free, the analytic continution principle is reduced to be the traction continuation \cite[]{Muskhelishvili1966}, as shown in Fig.~\ref{fig:2}c. In Fig.~\ref{fig:2}c, $ \vec{F}_{2} $ and $ \vec{F}_{3} $ are a pair of equilibrium resultants. The application of the analytic continuation principle does not conflict with the deduction process of obtaining Eq.~(\ref{eq:2.5}) in Ref \cite[]{Strack2002phdthesis}, and the geomaterial remains a doubly connected region before and after application of the analytic continuation principle in Fig.~\ref{fig:2}c.

According associative law of addition, Eq.~(\ref{eq:2.6}) can be further modified as
\begin{subequations}
  \label {eq:2.7}
  \begin{equation}
    \label {eq:2.7a}
    \varphi(z) = -\frac{{\rm i}F_{y}}{2\pi} \ln(z-\overline{z}_{c}) + \frac{-{\rm i} F_{y}}{2\pi(1+\kappa)} \ln \frac{z-z_{c}}{z-\overline{z}_{c}} + \varphi_{0}(z), \quad z \in {\bm \varOmega}
  \end{equation}
  \begin{equation}
    \label {eq:2.7b}
    \psi(z) = -\frac{{\rm i}F_{y}}{2\pi}\ln(z-\overline{z}_{c}) + \frac{-{\rm i} \kappa F_{y}}{2\pi(1+\kappa)} \ln \frac{z-z_{c}}{z-\overline{z}_{c}} + \psi_{0}(z), \quad z \in {\bm \varOmega}
  \end{equation}
\end{subequations}
Eq.~(\ref{eq:2.7}) indicates that $ \vec{F}_{1} $ and $ \vec{F}_{2} $ are set to a new pair of equilibrium resultants, and $ \vec{F}_{3} $ is left alone, as shown in~\ref{fig:2}d. No matter for the mechanical models in Fig.~\ref{fig:2}c or~\ref{fig:2}d, no static equilibrium can be established, which is hazardous to obtain fully reasonable stress and displacement. To verify such a hazard, we could examine the stress and displacement for the complex potentials in Eq.~(\ref{eq:2.7}).

The stress and displacement components within geomaterial in the physical plane $ z(xOy) $ can be given as
\begin{subequations}
  \label{eq:2.8}
  \begin{equation}
    \label{eq:2.8a}
    \left\{
      \begin{aligned}
        \sigma_{y}(z) & + \sigma_{x}(z) = 2 \left[\displaystyle\frac{{\rm d} \varphi(z)}{{\rm d}z} + \overline{\displaystyle\frac{{\rm d} \varphi(z)}{{\rm d}z}}\right] \\
        \sigma_{y}(z) & - \sigma_{x}(z) + 2{\rm i}\tau_{xy}(z) = 2\left[\overline{z} \displaystyle\frac{{\rm d}^{2} \varphi(z)}{{\rm d}z^{2}} + \displaystyle\frac{{\rm d} \psi(z)}{{\rm d}z}\right] \\
      \end{aligned}
    \right.
    , \quad
    z \in {\bm \varOmega}
  \end{equation}
  \begin{equation}
    \label{eq:2.8b}
    g(z) = 2G \left[ u(z) + {\rm i} v(z) \right] = \kappa \varphi(z) - z \overline{\frac{{\rm d} \varphi(z)}{{\rm d}{z}}} - \overline{\psi(z)}, \quad z \in {\bm \varOmega}
  \end{equation}
\end{subequations}
where $ \sigma_{x}(z) $, $ \sigma_{y} (z) $, and $ \tau_{xy} (z) $ denote horizontal, vertical, and shear stress components, respectively; $ u (z) $ and $ v (z) $ denote horizontal and vertical displacement components, respectively; $ G = \frac{E}{2(1+\nu)} $ denotes shear modulus of geomaterial, $ E $ and $ \nu $ denote elastic modulus and poisson's ratio of geomaterial, respectively, $ \kappa $ denotes the Kolosov coefficient with $ \kappa = 3 - 4\nu $ for plane strain and $ \kappa = \frac{3-\nu}{1+\nu} $ for plane stress, and $ \nu $ denotes the Poisson's ratio of geomaterial.

When substituting Eq.~(\ref {eq:2.7}) into Eq.~(\ref {eq:2.8}), we can find that the stress components are always finite, whereas the displacement components would be infinite when $ z $ approach infinity, indicating a unique displacement singularity at infinity. Such a property has been reported in Ref~\cite{Self2020JEM}, where a symmetrical modification is provided to attempt to fix the unique displacement singularity. However, the modification strategy in Ref~\cite{Self2020JEM} still depends on the modification depth (an exogenous parameter), and consequently the convergence is not as good as expected, which is shown in the numerical cases in Section~\ref{sec:comp-with-exist}. 

Therefore, we should seek a more reasonable strategy to equilibrate the unbalanced resultant $ \vec{F}_{1} $ and to constrain the unique displacement singularity at infinity. To simultaneously achieve both goals, we propose a new mechanical model in Fig.~\ref{fig:2}b, where displacement along the far-field ground surface is symmetrically constrained. Apparently, the displacement constraint along far-field ground surface would provide a constraining force to equilibrate the unbalanced resultant $ \vec{F}_{1} $ acting at singularity $ S_{1} $, thus, the originally unbalanced problem in Refs~\cite{Strack2002phdthesis, Strack_Verruijt2002buoyancy, Lu2016, Zengguisen2019, Self2020JEM} would turn to a balanced one. Furthermore, the displacement constraint along far-field ground surface can be expected to literally and mechanically constrain the displacement singularity at infinity. The ground surface turns from a fully free one to a partially free one, as shown in Fig.~\ref{fig:2}b.

\subsection{Mixed boundary value problem and conformal mapping}
\label{sec:problem-3}

\begin{figure}[hbt]
  \centering
  \begin{tikzpicture}
    \fill [gray!30] (-5,3) rectangle (3,-3);
    \draw [cyan, line width = 1pt] (-4.5,3) -- (2.5,3);
    \draw [line width = 1.5pt] (-5,3) -- (-4.5,3);
    \draw [line width = 1.5pt] (2.5,3) -- (3,3);
    \fill [pattern = north east lines] (-5,3) rectangle (-4.5,3.2); 
    \fill [pattern = north east lines] (2.5,3) rectangle (3,3.2);
    \fill [white] (-1,0) circle [radius = 1.5];
    \draw [violet, line width = 1pt] (-1,0) circle [radius = 1.5];
    \fill [red] (-4.5,3) circle [radius = 0.05];
    \fill [red] (2.5,3) circle [radius = 0.05];
    \node at (-4.5,3) [below, red] {$ T_{1} $};
    \node at (2.5,3) [below, red] {$ T_{2} $};
    \node at (-2,3) [above, cyan] {$ {\bm C}_{12} $};
    \draw [->] (-1,3) -- (0,3) node [above] {$ x $};
    \draw [->] (-1,3) -- (-1,4) node [right] {$ y $};
    \node at (-1,3) [below] {$ O $};
    \node at (-4.5,3.2) [above] {$ {\bm C}_{11} $};
    \node at (2.5,3.2) [above] {$ {\bm C}_{11} $};
    \node at (-1,-1.5) [above, violet] {$ {\bm C}_{2} $};
    \draw [->] (-1,0) -- ({-1+1.5*cos(-135)},{1.5*cos(-135)}) node [above] {$ R $};
    \draw [<->] (1.5,0) -- (1.5,3);
    \draw [dashed] (0.5,0) -- (1.7,0);
    \node at (1.5,1.5) [rotate = 90, above] {$ h $};
    \fill (-1,0) circle [radius = 0.05];
    \node at (-1,-4) [below] {(a) Lower half plane containing a shallow tunnel $ {\bm \varOmega} $};

    \fill [gray!30] (9-0.5,0) circle [radius = 3];
    \draw [cyan, line width = 1pt] (9-0.5,0) circle [radius = 3];
    \fill [pattern = north east lines] (12-0.5,0) -- (12.2-0.5,0) arc [start angle = 0, end angle = 15, radius = 3.2] -- ({9+3*cos(15)-0.5},{3*sin(15)}) arc [start angle = 15, end angle = 0, radius = 3];
    \fill [pattern = north east lines] (12-0.5,0) -- (12.2-0.5,0) arc [start angle = 0, end angle = -15, radius = 3.2] -- ({9+3*cos(-15)-0.5},{3*sin(-15)}) arc [start angle = -15, end angle = 0, radius = 3];
    \draw [line width = 1.5pt] (12-0.5,0) arc [start angle = 0, end angle = 15, radius = 3];
    \draw [line width = 1.5pt] (12-0.5,0) arc [start angle = 0, end angle = -15, radius = 3];
    \fill [white] (9-0.5,0) circle [radius = 1.5];
    \draw [violet, line width = 1pt] (9-0.5,0) circle [radius = 1.5];
    \draw [->] (9-0.5,0) -- (10-0.5,0) node [above] {$ \rho $};
    \draw [->] (9.3-0.5,0) arc [start angle = 0, end angle = 360, radius = 0.3];
    \node at (9.3-0.5,0) [above right] {$\theta$};
    \fill (9-0.5,0) circle [radius = 0.05];
    \node at (9-0.5,0) [below] {$o$};
    \fill [red] ({9+3*cos(15)-0.5},{3*sin(15)}) circle [radius = 0.05];
    \fill [red] ({9+3*cos(-15)-0.5},{3*sin(-15)}) circle [radius = 0.05];
    \node at ({9+3*cos(15)-0.5},{3*sin(15)}) [left, red] {$ t_{2} $};
    \node at ({9+3*cos(-15)-0.5},{3*sin(-15)}) [left, red] {$ t_{1} $};
    \node at (6-0.5,0) [right, cyan] {$ {\bm c}_{12} $};
    \node at (12-0.5,0) [left] {$ {\bm c}_{11} $};
    \node at (9-0.5,-1.5) [above, violet] {$ {\bm c}_{2} $};
    \draw [->] (9-0.5,0) -- ({9-0.5+1.5*cos(-135)},{1.5*sin(-135)}) node [above] {$ r $};
    \node at (9-0.5,-4) [below] {(b) Unit annulus $ {\bm \omega} $};

    \draw [line width = 2pt, ->] (3.5,1) -- (5.5,1);
    \draw [line width = 2pt, ->] (5.5,-1) -- (3.5,-1);
    \node at (4.5,1) [above] {$ \zeta = \frac {z + {\rm i} a} {z - {\rm i} a} $};
    \node at (4.5,1) [below] {Forward};
    \node at (4.5,-1) [above] {$ z = a \frac {1 + \zeta} {1 - \zeta} $};
    \node at ({9+3*cos(-135)-0.5},{3*sin(-135)}) [below left] {$ {\bm \omega}^{-} (\rho > 1) $};
    \node at (9-0.5,-2) [below] {$ {\bm \omega}^{+} (\rho < 1) $};
    \node at (4.5,-1) [below] {Backward};
  \end{tikzpicture}
  \caption{Schematic diagram of mixed boundary conditions and bidiretional conformal mappings of lower half plane containing a shallow tunnel $ {\bm \varOmega} $ and unit annulus $ {\bm \omega} $}
  \label{fig:3}
\end{figure}
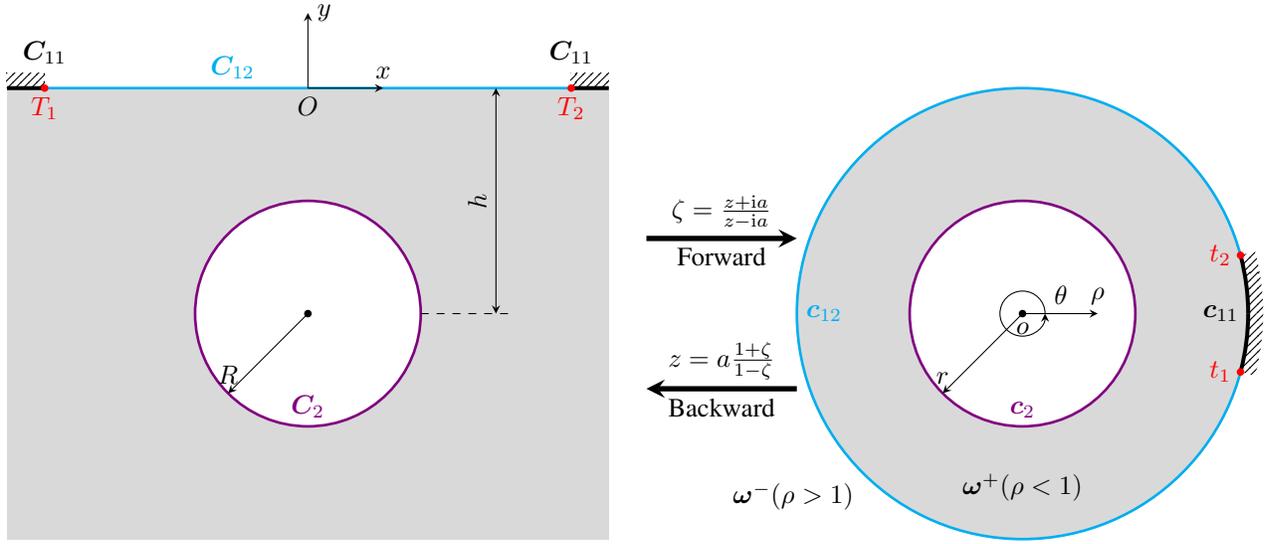

The mechanical model in Fig.~\ref{fig:2}b is no longer similar to the above mentioned ones in Refs~\cite{Strack2002phdthesis, Strack_Verruijt2002buoyancy, Lu2016, Zengguisen2019, Self2020JEM}, but is a mixed boundary value problem of elasticity. The ground surface is separated into the displacement-constrained segment $ {\bm C}_{11} $ and the partially free segment $ {\bm C}_{12} $, as shown in Fig.~\ref{fig:3}a. Both segments are axisymmetrical, and the intersection points between $ {\bm C}_{11} $ and $ {\bm C}_{12} $ are denoted by $ T_{1} $ and $ T_{2} $, respectively. The following mixed boundary conditions can be constructed according to the mechanical model in Fig.~\ref{fig:3}a: far-field ground surface ($ {\bm C}_{11} $) is constrained, and the rest part ($ {\bm C}_{12} $) is free, and the tunnel periphery ($ {\bm C}_{2} $) is subjected to tractions in Eq.~(\ref{eq:2.3}) with resultants in Eq.~(\ref{eq:2.4}). Apparently, the constrained far-field ground surface would produce a constraining force to equilibrate the unbalanced resultant along shallow tunnel periphery. In other words, the unbalanced model proposed by Strack \cite[]{Strack2002phdthesis} is modified to the static equilibrium model with mixed boundary conditions in this paper.

The boundary conditions along the ground surface can be expressed in the complex variable manner as
\begin{subequations}
  \label{eq:2.9}
  \begin{equation}
    \label{eq:2.9a}
    u(T) + {\rm i} v(T) = 0, \quad T \in {\bm C}_{11}
  \end{equation}
  \begin{equation}
    \label{eq:2.9b}
    X_{o}(T) + {\rm i} Y_{o}(T) = 0, \quad T \in {\bm C}_{12}
  \end{equation}
\end{subequations}
where $ u(T) $ and $ v(T) $ denote the horizontal and vertical displacement components along boundary $ {\bm C}_{11} $, respectively; $ X_{o}(T) $ and $ Y_{o}(T) $ denote horizontal and vertical tractions along boundary $ {\bm C}_{12} $, respectively. Apparently, Eq.~(\ref{eq:2.9}) already contains mixed boundary conditions.

Considering the mathematical equilibriums $\cos\langle\vec{n},x\rangle = \frac{{\rm d}y}{{\rm d}{S}}$ and $\cos\langle\vec{n},y\rangle = - \frac{{\rm d}x}{{\rm d}S}$ owing to outward normal vector $\vec{n}$ and clockwise positive direction of tunnel periphery, Eq.~(\ref{eq:2.2}) can be rewritten as
\begin{equation}
  \label{eq:2.2'}
  \tag{2.2'}
  \left\{
    \begin{aligned}
      X_{i}(S) = & \; - k_{0}\gamma y \frac{{\rm d}y}{{\rm d}S} \\
      Y_{i}(S) = & \; \gamma y \frac{{\rm d}x}{{\rm d}S} \\
    \end{aligned}
  \right.
  , \quad
  S: \; x^{2} + (y+h)^{2} = R^{2}
\end{equation}
Eq.~(\ref{eq:2.2'}) is the boundary condition along tunnel periphery, and is similar to those boundary conditions in the previously mentioned studies~\cite{Lu2016, Zengguisen2019, Self2020JEM}. Eqs.~(\ref{eq:2.9}) and~(\ref{eq:2.2'}) form the neccessary mathematical expressions for the mixed boundary value problem.

For better use of Eqs.~(\ref{eq:2.9}) and~(\ref{eq:2.2'}), the Verruijt's conformal mapping \cite[]{Verruijt1997traction,Verruijt1997displacement} is applied
to map the geometerial with a shallow circular tunnel in a lower half plane onto the unit annulus with inner radius of $ r $ (denoted by ${\bm \omega}$ in Fig.~\ref{fig:3}b) in a bidirectional manner via the following mapping functions:
\begin{subequations}
  \label{eq:2.10}
  \begin{equation}
    \label{eq:2.10a}
    \zeta(z) = \frac {z + {\rm i} a} {z - {\rm i} a}
  \end{equation}
  \begin{equation}
    \label{eq:2.10b}
    z(\zeta) = - {\rm i} a \frac {1 + \zeta} {1 - \zeta}
  \end{equation}
\end{subequations}
where
\begin{equation}
  \label{eq:2.11}
  \left\{
    \begin{aligned}
      a = & \; h \frac {1 - r^{2}} {1 + r^{2}} \\
      r = & \; \frac {R} {h + \sqrt{h^{2} - R^{2}}} \\
    \end{aligned}
  \right.
\end{equation}
Via the conformal mapping functions in Eq.~(\ref{eq:2.10}), the boundaries $ {\bm C}_{11} $, $ {\bm C}_{12} $, and $ {\bm C}_{2} $ in the lower half plane in Fig.~\ref{fig:3}a are respectively mapped onto $ {\bm c}_{11} $, $ {\bm c}_{12} $, and $ {\bm c}_{2} $ in the unit annulus in Fig.~\ref{fig:3}b in a bidirectional manner. Boundaries ${\bm c}_{11}$ and ${\bm c}_{12}$ together can be denoted by ${\bm c}_{1}$. Corresponding to region ${\bm \varOmega}$ in the physical plane $z$, region ${\bm \omega}$ in the mapping plane $\zeta$ is also a closure that contains ${\bm c}_{1}$, ${\bm c}_{2}$, and the unit annulus bounded by both boundaries. It should be addressed that since the infinity point $ z = \infty $ in the lower half plane is mapped onto point $ \zeta = 1 $ in the unit annulus, the infinite boundary $ {\bm C}_{11} $ is mapped onto a finite arc $ {\bm c}_{11} $ correspondingly. The intersection points $ T_{1} $ and $ T_{2} $ are also mapped onto points $ t_{1} $ and $ t_{2} $, respectively. Since $T_{1}$ and $ T_{2} $ are axisymmetrical, $t_{1}$ and $t_{2}$ would also be axisymmetrical. Assume that the horizontal coordinates of points $T_{1}$ and $T_{2}$ are $-x_{0}$ and $x_{0}$, respectively; then the poloar angles of points $t_{1}$ and $t_{2}$ would be $-\theta_{0}$ and $\theta_{0}$, respectively, where $\theta_{0} = -{\rm i}\ln \frac{x_{0}+{\rm i}a}{x_{0}-{\rm i}a}$. 

The mixed boundary conditions in the lower half plane in Eq.~(\ref{eq:2.9}) and (\ref{eq:2.2'}) can also be conformally mapped onto the ones in the unit annulus as
\begin{subequations}
  \label{eq:2.12}
  \begin{equation}
    \label{eq:2.12a}
    u(t) + {\rm i} v(t) = u(T) + {\rm i} v(T) = 0, \quad t = {\rm e}^{{\rm i}\theta} \in {\bm c}_{11}
  \end{equation}
  \begin{equation}
    \label{eq:2.12b}
    {\rm e}^{{\rm i}\theta} \frac {z^{\prime}(t)} {|z^{\prime}(t)|} \cdot [\sigma_{\rho}(t) + {\rm i} \tau_{\rho\theta}(t)] = X_{o}(T) + {\rm i} Y_{o}(T) = 0, \quad t = {\rm e}^{{\rm i}\theta} \in {\bm c}_{12}
  \end{equation}
  \begin{equation}
    \label{eq:2.12c}
    {\rm e}^{{\rm i}\theta} \frac {z^{\prime}(s)} {|z^{\prime}(s)|} \cdot [\sigma_{\rho}(s) + {\rm i} \tau_{\rho\theta}(s)] = X_{i}(S) + {\rm i} Y_{i}(S) = - k_{0}\gamma y \frac{{\rm d}y}{{\rm d}S} + {\rm i}\gamma y \frac{{\rm d}x}{{\rm d}S}, \quad s = r {\rm e}^{{\rm i}\theta} \in {\bm c}_{2}
  \end{equation}
\end{subequations}
where the points $ t $ and $ s $ in the unit annulus are respectively corresponding to the points $T$ and $S$ in the lower half plane; $ \sigma_{\rho} $ and $ \tau_{\rho\theta} $ denote radial and shear stress components, respectively; $z^{\prime}(t)$ and $z^{\prime}(s)$ denote taking derivatives of $z(t)$ and $z(s)$, respectively. Now we should solve the mixed boundary value problem in Eq.~(\ref{eq:2.12}).

\section{Analytic continuation and Riemann-Hilbert problem}
\label{sec:analyt-cont-riem}

To solve the mixed boundary value problem in Eq.~(\ref{eq:2.12}) and to obtain stress and displacement solution of shallow tunnel excavation, the analytic continuation \cite[]{Muskhelishvili1966} is applied to the proposed mechanical model. Considering the backward conformal mapping in Eq.~(\ref{eq:2.10b}), the stress and displacement in geomaterial in Eq.~(\ref{eq:2.8}) can be mapped onto the mapping plane $\zeta$ as
\begin{subequations}
  \label{eq:3.1}
  \begin{equation}
    \label{eq:3.1a}
    \sigma_{\theta}(\zeta) + \sigma_{\rho}(\zeta) = 2 \left[ \varPhi(\zeta) + \overline{\varPhi(\zeta)} \right], \quad \zeta \in {\bm \omega}
  \end{equation}
  \begin{equation}
    \label{eq:3.1b}
    \sigma_{\rho}(\zeta) + {\rm i} \tau_{\rho\theta}(\zeta) = \varPhi(\zeta) + \overline{\varPhi(\zeta)} - {\rm e}^{-2{\rm i}\theta} \left[ \frac{z(\zeta)}{z^{\prime}(\zeta)} \overline{\varPhi^{\prime}(\zeta)} + \frac{\overline{z^{\prime}(\zeta)}}{z^{\prime}(\zeta)} \overline{\varPsi(\zeta)} \right], \quad \zeta \in {\bm \omega}
  \end{equation}
  \begin{equation}
    \label{eq:3.1c}
    g(\zeta) = 2G \left[ u(\zeta) + {\rm i} v(\zeta) \right] = \kappa \varphi(\zeta) - z(\zeta) \overline{\varPhi(\zeta)} - \overline{\psi(\zeta)}, \quad \zeta \in {\bm \omega}
  \end{equation}
  where
  \begin{equation*}
    \left\{
      \begin{aligned}
        & \varPhi(\zeta) = \frac{\varphi^{\prime}(\zeta)}{z^{\prime}(\zeta)}\\
        & \varPsi(\zeta) = \frac{\psi^{\prime}(\zeta)}{z^{\prime}(\zeta)}\\
      \end{aligned}
    \right.
  \end{equation*}
\end{subequations}
Eq.~(\ref{eq:3.1a}) requires $\varPhi(\zeta)$ to be analytic within region ${\bm \omega}$, and Eq.~(\ref{eq:3.1b}) further requires the last two items on the right-hand side to be analytic within region ${\bm \omega}$, which can be expanded as
\begin{subequations}
  \label{eq:3.2}
  \begin{equation}
    \label{eq:3.2a}
    \left\{
      \begin{aligned}
        \frac{z(\zeta)}{z^{\prime}(\zeta)} \overline{\varPhi^{\prime}(\zeta)} = & \; \frac{1-\zeta^{2}}{2} \overline{\varPhi^{\prime}(\zeta)} \\
        \frac{\overline{z^{\prime}(\zeta)}}{z^{\prime}(\zeta)} \overline{\varPsi(\zeta)} = & \; - (1-\zeta)^{2} \frac{\overline{\varPsi(\zeta)}}{(1-\overline{\zeta})^{2}}
      \end{aligned}
    \right.
    , \quad \zeta \in {\bm \omega}
  \end{equation}
  The first equation in Eq.~(\ref{eq:3.2a}) is apparently analytic within region ${\bm \omega}$, and the second equation would be analytic within region ${\bm \omega}$, as long as $\varPsi(\zeta)$ can be expressed as multiplication of $(1-\zeta)$ and an analytic function without any singularity in region ${\bm \omega}$, because
  \begin{equation*}
    \lim_{\zeta\rightarrow 1} \frac{(1-\zeta)^{2}}{1-\overline{\zeta}} = \lim_{\zeta\rightarrow 0} \frac{\zeta^{2}}{\overline{\zeta}} = \lim_{\rho\rightarrow 0} \frac{\rho^{2}}{\rho} {\rm e}^{3{\rm i}\theta} = 0
  \end{equation*}
  always stands. Eq.~(\ref{eq:3.2a}) indicates that $\zeta=1$ should not be a singularity point for stress.

  The complex potentials in Eq.~(\ref{eq:3.1c}) can be expanded as
  \begin{equation}
    \label{eq:3.2b}
    \left\{
      \begin{aligned}
        \varphi(\zeta) = & \; \int \varphi^{\prime}(\zeta) {\rm d}\zeta = \int z^{\prime}(\zeta) \varPhi(\zeta) {\rm d}\zeta = -2{\rm i}a \int \frac{\varPhi(\zeta)}{(1-\zeta)^{2}} {\rm d}\zeta \\
        \psi(\zeta) = & \; \int \psi^{\prime}(\zeta) {\rm d}\zeta = \int z^{\prime}(\zeta) \varPsi(\zeta) {\rm d}\zeta  = -2{\rm i}a \int \frac{\varPsi(\zeta)}{(1-\zeta)^{2}} {\rm d}\zeta \\
      \end{aligned}
    \right.
    , \quad \zeta \in {\bm \omega}
  \end{equation}
  Eq.~(\ref{eq:3.2b}) indicates that $\zeta=1$ may be a second-order singularity point. However, according to Eq.~(\ref{eq:2.12a}), displacement along the whole boundary ${\bm c}_{11}$ should be equal to zero, including point $\zeta=1$. Eq.~(\ref{eq:3.2b}) should therfore be analytic within region ${\bm \omega}$, as long as $\varPhi(\zeta)$ and $\varPsi(\zeta)$ can be expressed as a multiplication of $(1-\zeta)^{2}$ and an analytic function in region ${\bm \omega}$ without any singularity. Eq.~(\ref{eq:3.2b}) indicates that displacement has a stronger requirement on $\varPhi(\zeta)$ and $\varPsi(\zeta)$, comparing to stress in Eq.~(\ref{eq:3.2a}). In other words, Eqs.~(\ref{eq:3.2a}) and (\ref{eq:3.2b}) both require that the possible second-order singularity $\zeta=1$ in $z^{\prime}(\zeta)$ should not be a singularity point for $z^{\prime}(\zeta)\varPhi(\zeta)$ and $z^{\prime}(\zeta)\varPsi(\zeta)$, as well as $\varPhi(\zeta)$ and $\varPsi(\zeta)$. Eq.~(\ref{eq:3.1c}) further requires the remaining item to be analytic within region ${\bm \omega}$, which can be expanded as
  \begin{equation}
    \label{eq:3.2c}
    z(\zeta) \overline{\varPhi(\zeta)} = \overline{z^{\prime}(\zeta)\varPhi(\zeta)} \cdot \frac{z(\zeta)}{\overline{z^{\prime}(\zeta)}} = - \frac{1+\zeta}{2} \overline{z^{\prime}(\zeta)\varPhi(\zeta)} \cdot \frac{(1-\overline{\zeta})^{2}}{1-\zeta}, \quad \zeta \in {\bm \omega}
  \end{equation}
  Whether or not Eq.~(\ref{eq:3.2c}) is analytic within region ${\bm \omega}$ depends on the last item on the possible singulariy point $\zeta=1$, which can be computed in a similar manner above:
  \begin{equation*}
    \lim_{\zeta\rightarrow 1}\frac{(1-\overline{\zeta})^{2}}{1-\zeta} = \lim_{\zeta \rightarrow 0} \frac{\overline{\zeta}^{2}}{\zeta} = \lim_{\rho \rightarrow 0} \frac{\rho^{2}}{\rho}{\rm e}^{-3{\rm i}\theta} = 0
  \end{equation*}
  Eqs.~(\ref{eq:3.2b}) and (\ref{eq:3.2c}) indicate that $\zeta=1$ should not be a singularity point for displacement as well.
\end{subequations}

With the multiplier $(1-\zeta)^{2}$ in $\varPhi(\zeta)$ and $\varPsi(\zeta)$, the stress components would vanish at point $\zeta=1$ in the mapping plane according to Eqs.~(\ref{eq:3.1a}) and (\ref{eq:3.1b}), correspondingly indicating that the stress components at infinity in the physical plane would vanish. Such a result is identical to our expectation. Moreover, owing to the axisymmetry of the mechanical model, no far-field rotation or moment need be considered.

Therefore, before any further discussion, we find that $\varPhi(\zeta)$ and $\varPsi(\zeta)$, as well as $z^{\prime}(\zeta)\varPhi(\zeta)$ and $z^{\prime}(\zeta)\varPsi(\zeta)$, should be analytic within region ${\bm \omega}$ without any singularity. To facilitate our discussion, we denote the regions $r \leq \rho < 1$ and $1 < \rho \leq r^{-1}$ by ${\bm \omega}^{+}$ and ${\bm \omega}^{-}$, respectively, as shown in Fig.~\ref{fig:3}b, and ${\bm \omega} = {\bm \omega}^{+}\cup{\bm c}_{1}$. Apparently, $\varPhi(\zeta)$ and $\varPsi(\zeta)$, as well as $z^{\prime}(\zeta)\varPhi(\zeta)$ and $z^{\prime}(\zeta)\varPsi(\zeta)$, are all analytic within region ${\bm \omega}^{+}$. The analytic continuation principle will be used below to find solution of the stress and displacement components in the unit annulus in Eq.~(\ref{eq:3.1}), according to the mixed boundary conditions in Eq.~(\ref{eq:2.12}).

Substituting Eq.~(\ref{eq:3.1b}) into Eq.~(\ref{eq:2.12b}) and noting $ \overline{t} = t^{-1} $ yields
\begin{equation}
  \label{eq:3.3}
  \varPhi(t) = - \overline{\varPhi}\left( t^{-1} \right) + \frac{1}{t^{2}} \left[ \frac{z\left(t\right)}{z^{\prime}(t)}  \overline{\varPhi^{\prime}}\left( t^{-1} \right) + \frac{\overline{z^{\prime}}\left( t^{-1} \right)}{z^{\prime}(t)}  \overline{\varPsi}\left( t^{-1} \right) \right], \quad t = {\rm e}^{{\rm i}\theta} \in {\bm c}_{12} 
\end{equation}
Replacing $t$ with $\zeta = \rho \cdot {\rm e}^{{\rm i}\theta}(1 < \rho \leq r^{-1})$, Eq.~(\ref{eq:3.3}) would turn to
\begin{equation}
  \label{eq:3.4}
  \varPhi(\zeta) = - \overline{\varPhi}\left( \zeta^{-1} \right) + \frac{1}{\zeta^{2}} \left[ \frac{z\left(\zeta\right)}{z^{\prime}(\zeta)} \overline{\varPhi^{\prime}}\left( \zeta^{-1} \right) + \frac{\overline{z^{\prime}}\left( \zeta^{-1} \right)}{z^{\prime}(\zeta)} \overline{\varPsi}\left( \zeta^{-1} \right) \right], \quad \zeta = \rho \cdot {\rm e}^{{\rm i}\theta} \in {\bm \omega}^{-} ( 1 < \rho \leq r^{-1})
\end{equation}
The polar radius range $1 < \rho \leq r^{-1}$ ensure all the items of the right-hand side of Eq.~(\ref{eq:3.4}) to be analytic. Eq.~(\ref{eq:3.4}) can be modified as
\begin{equation}
  \label{eq:3.5}
  z^{\prime}(\zeta) \varPhi(\zeta) = - z^{\prime}(\zeta) \overline{\varPhi}\left( \zeta^{-1} \right) + \frac{1}{\zeta^{2}} z\left(\zeta\right) \overline{\varPhi^{\prime}}\left( \zeta^{-1} \right) + \frac{1}{\zeta^{2}} \overline{z^{\prime}}\left( \zeta^{-1} \right) \overline{\varPsi}\left( \zeta^{-1} \right), \quad \zeta = \rho \cdot {\rm e}^{{\rm i}\theta} \in {\bm \omega}^{-} (1 < \rho \leq r^{-1})
\end{equation}
Eq.~(\ref{eq:3.5}) shows that $z^{\prime}(\zeta) \varPhi(\zeta)$ is analytic within region ${\bm \omega}^{-}(1 < \rho \leq r^{-1})$ without any singularity point. Combining with Eqs.~(\ref{eq:3.2}), $z^{\prime}(\zeta)\varPhi(\zeta)$ should be analytic within the region ${\bm \omega}^{+} \cup {\bm c}_{1} \cup {\bm \omega}^{-}$ without any singularity point. Replacing $ \zeta $ with $ \overline{\zeta}^{-1} $ in Eq.~(\ref{eq:3.5}) and taking conjugate yields
\begin{equation}
  \label{eq:3.6}
  z^{\prime}(\zeta) \varPsi(\zeta) = \frac{1}{\zeta^{2}} \overline{z^{\prime}} \left( \zeta^{-1} \right) \overline{\varPhi} \left( \zeta^{-1} \right) + \frac{1}{\zeta^{2}} \overline{z^{\prime}}\left( \zeta^{-1} \right) \varPhi(\zeta) - \overline{z} \left( \zeta^{-1} \right) \varPhi^{\prime}(\zeta), \quad \zeta = \rho \cdot {\rm e}^{{\rm i}\theta} \in {\bm \omega}^{+}(r \leq \rho < 1)
\end{equation}
Eq.~(\ref{eq:3.6}) shows that $ z^{\prime}(\zeta) \varPsi(\zeta) $ is defined and analytic within region ${\bm \omega}^{+}$ without any singularity point, and can be expressed by combination of $\varPhi(\zeta)$ in regions ${\bm \omega}^{+}$ and ${\bm \omega}^{-}$. Thus, we only need to focus on $ z^{\prime}(\zeta) \varPhi(\zeta) $.

Integrating Eq.~(\ref{eq:3.5}) by $\zeta$ yields
\begin{equation}
  \label{eq:3.7}
  \varphi(\zeta) = -z(\zeta) \overline{\varPhi}\left( \zeta^{-1} \right) - \overline{\psi} \left( \zeta^{-1} \right) + {\rm constant}, \quad \zeta = \rho \cdot {\rm e}^{{\rm i}\theta} \in {\bm \omega}^{-} ( 1 < \rho \leq r^{-1} )
\end{equation}
Replacing $\zeta$ with $ \overline{\zeta}^{-1} $ in Eq.~(\ref{eq:3.8}) yields
\begin{equation}
  \label{eq:3.8}
  \varphi\left( \overline{\zeta}^{-1} \right) = -z\left( \overline{\zeta}^{-1} \right) \overline{\varPhi\left( \zeta \right)} - \overline{\psi \left( \zeta \right)} + {\rm constant}, \quad \zeta = \rho \cdot {\rm e}^{{\rm i}\theta} \in {\bm \omega}^{+} ( r \leq \rho < 1)
\end{equation}
Taking derivative of $ \overline{\zeta} $ in Eq.~(\ref{eq:3.8}) yields
\begin{equation}
  \label{eq:3.9}
  \varPhi \left( \overline{\zeta}^{-1} \right) = \overline{\zeta}^{2} \frac{z\left( \overline{\zeta}^{-1} \right)} {z^{\prime}\left( \overline{\zeta}^{-1} \right)} \overline{\varPhi^{\prime}(\zeta)} + \overline{\zeta}^{2} \frac{\overline{z^{\prime}(\zeta)}} {z^{\prime}\left( \overline{\zeta}^{-1} \right)} \overline{\varPsi(\zeta)} - \overline{\varPhi(\zeta)}, \quad \zeta = \rho \cdot {\rm e}^{{\rm i}\theta} \in {\bm \omega}^{+} (r \leq \rho < 1)
\end{equation}
Eq.~(\ref{eq:3.9}) can be transformed as
\begin{equation}
  \label{eq:3.10}
  \overline{\varPhi(\zeta)} = - \varPhi \left( \overline{\zeta}^{-1} \right) + \overline{\zeta}^{2} \frac{z\left( \overline{\zeta}^{-1} \right)} {z^{\prime}\left( \overline{\zeta}^{-1} \right)} \overline{\varPhi^{\prime}(\zeta)} + \overline{\zeta}^{2} \frac{\overline{z^{\prime}(\zeta)}} {z^{\prime}\left( \overline{\zeta}^{-1} \right)} \overline{\varPsi(\zeta)} , \quad \zeta = \rho \cdot {\rm e}^{{\rm i}\theta} \in {\bm \omega}^{+} (r \leq \rho < 1)
\end{equation}
Substituting Eq.~(\ref{eq:3.10}) into Eq.~(\ref{eq:3.1b}) yields
\begin{equation}
  \label{eq:3.11}
  \begin{aligned}
    \sigma_{\rho}(\zeta) + {\rm i} \tau_{\rho\theta}(\zeta) & = \overline{\zeta}^{2} \left[ \frac{z\left( \overline{\zeta}^{-1} \right)} {z^{\prime}\left( \overline{\zeta}^{-1} \right)} - \frac{1}{\rho^{2}} \frac{z(\zeta)} {z^{\prime}(\zeta)} \right] \overline{\varPhi^{\prime}(\zeta)} + \overline{\zeta}^{2}\left[ \frac{\overline{z^{\prime}(\zeta)}} {z^{\prime}\left( \overline{\zeta}^{-1} \right) } - \frac{1}{\rho^{2}} \frac{\overline{z^{\prime}(\zeta)}} {z^{\prime}(\zeta)} \right] \overline{\varPsi(\zeta)} \\
                                                            & + \varPhi(\zeta) - \varPhi \left( \overline{\zeta}^{-1} \right), \quad \zeta = \rho \cdot {\rm e}^{{\rm i}\theta} \in {\bm \omega}^{+} (r \leq \rho < 1)
  \end{aligned}
\end{equation}

Taking deriative of $ \theta$ of Eq.~(\ref{eq:3.1c}) yields
\begin{equation}
  \label{eq:3.12}
  \frac{{\rm d}g(\zeta)}{{\rm d}\theta} = {\rm i}\zeta \cdot \kappa z^{\prime}(\zeta) \varPhi(\zeta) - {\rm i}\zeta \cdot z^{\prime}(\zeta) \overline{\varPhi(\zeta)} + {\rm i}\overline{\zeta} \cdot z(\zeta) \overline{\varPhi^{\prime}(\zeta)} + {\rm i}\overline{\zeta} \cdot \overline{z^{\prime}(\zeta)} \; \overline{\varPsi(\zeta)}, \quad \zeta = \rho \cdot {\rm e}^{{\rm i}\theta} \in {\bm \omega}^{+} (r \leq \rho < 1)
\end{equation}
Substituting Eq.~(\ref{eq:3.10}) into Eq.~(\ref{eq:3.12}) yields
\begin{equation}
  \label{eq:3.13}
  \begin{aligned}
    \frac{{\rm d}g(\zeta)}{{\rm d}\theta} & =  {\rm i}\overline{\zeta} \cdot \left[ z(\zeta) - \rho^{2} \frac{z^{\prime}(\zeta)} {z^{\prime}\left( \overline{\zeta}^{-1} \right)} z\left( \overline{\zeta}^{-1} \right) \right] \overline{\varPhi^{\prime}(\zeta)} + {\rm i}\overline{\zeta} \; \overline{z^{\prime}(\zeta)} \cdot \left[ 1 - \rho^{2} \frac{z^{\prime}(\zeta)} {z^{\prime}\left( \overline{\zeta}^{-1} \right)} \right] \overline{\varPsi(\zeta)} \\
                                          & + {\rm i}\zeta \cdot z^{\prime}(\zeta) \left[\kappa  \varPhi(\zeta) + \varPhi\left(\overline{\zeta}^{-1}\right) \right], \quad \zeta = \rho \cdot {\rm e}^{{\rm i}\theta} \in {\bm \omega}^{+} (r \leq \rho < 1)
  \end{aligned}
\end{equation}
Taking deriative of $\zeta$ of Eq.~(\ref{eq:3.1c}) with consideration of Eq.~(\ref{eq:3.13}) yields
\begin{equation}
  \label{eq:3.14}
  \begin{aligned}
    \frac{{\rm d}g(\zeta)}{{\rm d}\zeta} = \frac{{\rm d}g(\zeta)}{{\rm d}\theta} \cdot \frac{1}{{\rm i}\zeta} & =  {\rm e}^{-2{\rm i}\theta} \cdot \left[ z(\zeta) - \rho^{2} \frac{z^{\prime}(\zeta)} {z^{\prime}\left( \overline{\zeta}^{-1} \right)} z\left( \overline{\zeta}^{-1} \right) \right] \overline{\varPhi^{\prime}(\zeta)} + {\rm e}^{-2{\rm i}\theta} \overline{z^{\prime}(\zeta)} \cdot \left[ 1 - \rho^{2} \frac{z^{\prime}(\zeta)} {z^{\prime}\left( \overline{\zeta}^{-1} \right)} \right] \overline{\varPsi(\zeta)} \\
                                                                                                              & + z^{\prime}(\zeta) \left[\kappa  \varPhi(\zeta) + \varPhi\left(\overline{\zeta}^{-1}\right) \right], \quad \zeta = \rho \cdot {\rm e}^{{\rm i}\theta} \in {\bm \omega}^{+} (r \leq \rho < 1)
  \end{aligned}
\end{equation}

Respectively substituting Eqs.~(\ref{eq:3.11}) and (\ref{eq:3.14}) into Eqs.~(\ref{eq:2.12b}) and (\ref{eq:2.12a}) yields
\begin{subequations}
  \label{eq:3.15}
  \begin{equation}
    \label{eq:3.15a}
    \varPhi^{+}(t) - \varPhi^{-}(t)= 0, \quad t = {\rm e}^{{\rm i}\theta} \in {\bm c}_{12}
  \end{equation}
  \begin{equation}
    \label{eq:3.15b}
    z^{\prime}(t) \left[ \kappa \varPhi^{+}(t) + \varPhi^{-}(t) \right] = 0, \quad t = {\rm e}^{{\rm i}\theta}\in {\bm c}_{11}
  \end{equation}
\end{subequations}
where $\varPhi^{+}(t)$ and $\varPhi^{-}(t)$ denote values of $\varPhi(\zeta)$ approaching boundary ${\bm c}_{1}$ from regions ${\bm \omega}^{+}$ and ${\bm \omega}^{-}$, respectively. Considering the the indefinite integrals in Eq.~(\ref{eq:3.2b}), it would be more computationally convenient to transform Eq.~(\ref{eq:3.15}) into the following form:
\begin{subequations}
  \label{eq:3.16}
  \begin{equation}
    \label{eq:3.16a}
    z^{\prime}(\zeta)[\sigma_{\rho}(\zeta)+{\rm i}\tau_{\rho\theta}(\zeta)]|_{\zeta\rightarrow t} = \left[ z^{\prime}(t) \varPhi(t) \right]^{+} - \left[ z^{\prime}(t) \varPhi(t) \right]^{-} = 0, \quad t = {\rm e}^{{\rm i}\theta} \in {\bm c}_{12}
  \end{equation}
  \begin{equation}
    \label{eq:3.16b}
    {\frac{{\rm d}g(\zeta)}{{\rm d}\zeta}}|_{\zeta\rightarrow t} = \kappa \left[ z^{\prime}(t) \varPhi(t) \right]^{+} + \left[ z^{\prime}(t) \varPhi(t) \right]^{-} = 0, \quad t = {\rm e}^{{\rm i}\theta}\in {\bm c}_{11}
  \end{equation}
\end{subequations}
where $\left[ z^{\prime}(t) \varPhi(t) \right]^{+}$ and $\left[ z^{\prime}(t) \varPhi(t) \right]^{-}$ denote values of $ z^{\prime}(\zeta) \varPhi(\zeta) $ approaching boundary $ {\bm c}_{1} $ from regions $ {\bm \omega}^{+} $ and $ {\bm \omega}^{-} $, respectively. The transformation from Eq.~(~\ref {eq:3.15} ) to (~\ref {eq:3.16} ) is reasonable, because $ z^{ \prime } ( \zeta ) \neq 0 $ is single-valued, and would not alter the possible multi-valuedness of $ \varPhi (\zeta) $. Therefore, the mixed boundary conditions in Eq.~(\ref{eq:2.12}) turn to a homogenerous Riemann-Hilbert problem Eq.~(\ref{eq:3.16}) with extra constraints of Eq.~(\ref{eq:2.12c}).

\section{Solution of the Riemann-Hilbert problem}
\label{sec:solut-riem-hilb}

\subsection{Complex potential expansion}
\label{sec:solution-riemann-hilbert-1}

As pointed out in last section that $z^{\prime}(\zeta)\varPhi(\zeta)$ should be analytic within the annular region ${\bm \omega}^{+} \cup {\bm c}_{1} \cup {\bm \omega}^{-}$ without any singularity, the general solution for the Riemann-Hilbert problem in Eq.~(\ref{eq:3.16}) should satisfy the following form according to Plemelj formula~\cite[]{Muskhelishvili1966} as
\begin{equation}
  \label{eq:4.1}
  \varphi^{\prime}(\zeta) = z^{\prime}(\zeta) \varPhi(\zeta) = X(\zeta) \sum\limits_{n=-\infty}^{\infty} {\rm i} d_{n} \zeta^{n}, \quad \zeta = \rho \cdot {\rm e}^{{\rm i}\theta}, \quad r \leq \rho \leq r^{-1}
\end{equation}
where
\begin{equation}
  \label{eq:4.2}
  X(\zeta) = (\zeta - {\rm e}^{-{\rm i}\theta_{0}})^{-\frac{1}{2}-{\rm i}\lambda}(\zeta - {\rm e}^{{\rm i}\theta_{0}})^{-\frac{1}{2}+{\rm i}\lambda}, \quad \lambda = \frac{\ln \kappa}{2\pi}
\end{equation}
$ d_{n} $ denote coefficients to be determined, and should be real owing to the axisymmetry. The imaginary unit in front of $d_{n}$ denotes that the coefficients of the sum above are pure imaginary numbers, because the geometry of the remaining geomaterial, as well as the tractions and constraints acting upon both boundaries, are all axisymmetrical about $y$ axis. The infinite Laurent series suggest bipole $\zeta = 0$ and $\zeta = \infty$ for the annular region $r \leq \rho \leq r^{-1}$. To solve the coefficients, all items in Eq.~(\ref{eq:4.1}) should be prepared into rational series. $ X(\zeta) $ can be respectivley expanded in regions $ {\bm \omega}^{+} $ and $ {\bm \omega}^{-} $ using Taylor expansions as

\begin{subequations}
  \label{eq:4.3}
  \begin{equation}
    \label{eq:4.3a}
    X(\zeta) = \sum\limits_{k=0}^{\infty} \alpha_{k} \zeta^{k}, \quad \zeta \in {\bm \omega}^{+}
  \end{equation}
  \begin{equation}
    \label{eq:4.3b}
    X(\zeta) = \sum\limits_{k=1}^{\infty} \beta_{k} \zeta^{-k}, \quad \zeta \in {\bm \omega}^{-}
  \end{equation}
\end{subequations}
where 
\begin{equation*}
  \left\{
    \begin{aligned}
      \alpha_{0} = & \; - {\rm e}^{-2\lambda\theta_{0}} \\
      \alpha_{1} = & \; - {\rm e}^{-2\lambda\theta_{0}} \left( \cos\theta_{0} - 2\lambda \sin\theta_{0} \right) \\
      \alpha_{k} = & \; - {\rm e}^{-2\lambda\theta_{0}} \cdot (-1)^{k} \left[ \frac{\mathsf{a}_{k}}{k!} + \frac{\overline{\mathsf{a}}_{k}}{k!} + \sum\limits_{l=1}^{k-1} \frac{\mathsf{a}_{l}}{l!} \cdot \frac{\overline{\mathsf{a}}_{k-l}}{(k-l)!} \right], \quad k \geq 2 \\
      \mathsf{a}_{k} = & \; \prod_{l=1}^{k} \left( \frac{1}{2} - {\rm i}\lambda -l \right) \cdot {\rm e}^{{\rm i}k\theta_{0}}
    \end{aligned}
  \right.
\end{equation*}
\begin{equation*}
  \left\{
    \begin{aligned}
      \beta_{1} = \; & 1 \\
      \beta_{2} = \; & \cos\theta_{0} + 2\lambda \sin\theta_{0} \\
      \beta_{k} = \; & (-1)^{k-1} \left[ \frac{\mathsf{b}_{k-1}}{(k-1)!} + \frac{\overline{\mathsf{b}}_{k-1}}{(k-1)!} + \sum\limits_{l=1}^{k-2} \frac{\mathsf{b}_{l}}{l!} \cdot \frac{\overline{\mathsf{b}}_{k-1-l}}{(k-1-l)!} \right], \quad k \geq 3 \\
      \mathsf{b}_{k} = & \; \prod_{l=1}^{k} \left( \frac{1}{2} - {\rm i}\lambda - l \right) \cdot {\rm e}^{-{\rm i}k\theta_{0}}
    \end{aligned}
  \right.
\end{equation*}
Last equations show that $\alpha_{k}$ and $\beta_{k}$ are both real, which is identical to the axisymmetry of the model. The branch of $ X(\zeta) $ in Eq.~(\ref{eq:4.3b}) is correct, since it potentially guarantees $ \lim\limits_{\zeta\to\infty} \zeta X(\zeta) = 1 $.

Substituting Eq.~(\ref{eq:4.3}) into Eq.~(\ref{eq:4.1}) yields
\begin{subequations}
  \label{eq:4.4}
  \begin{equation}
    \label{eq:4.4a}
    \varphi^{\prime}(\zeta) = z^{\prime}(\zeta) \varPhi(\zeta) = \sum\limits_{k=-\infty}^{\infty} {\rm i} A_{k} \zeta^{k}, \quad A_{k} = \sum\limits_{n=-\infty}^{k} \alpha_{k-n}d_{n}, \quad \zeta \in {\bm \omega}^{+}
  \end{equation}
  \begin{equation}
    \label{eq:4.4b}
    \varphi^{\prime}(\zeta) = z^{\prime}(\zeta) \varPhi(\zeta) = \sum\limits_{k=-\infty}^{\infty} {\rm i} B_{k} \zeta^{k}, \quad B_{k} = \sum\limits_{n=k+1}^{\infty} \beta_{n-k}d_{n}, \quad \zeta \in {\bm \omega}^{-}
  \end{equation}
\end{subequations}
Since $z^{\prime}(\zeta) \neq 0$ for both regions ${\bm \omega}^{+}$ and ${\bm \omega}^{-}$, Eq.~(\ref{eq:4.4}) can be transformed as
\begin{equation}
  \label{eq:4.4a'}
  \tag{4.4a'}
  \varPhi(\zeta) = \frac{(1-\zeta)^{2}}{-2{\rm i}a} \sum\limits_{k=-\infty}^{\infty} {\rm i} A_{k} \zeta^{k}, \quad \zeta \in {\bm \omega}^{+}
\end{equation}
\begin{equation}
  \label{eq:4.4b'}
  \tag{4.4b'}
  \varPhi(\zeta) = \frac{(1-\zeta)^{2}}{-2{\rm i}a} \sum\limits_{k=-\infty}^{\infty} {\rm i} B_{k} \zeta^{k}, \quad \zeta \in {\bm \omega}^{-}
\end{equation}
Eqs.~(\ref{eq:4.4a}) and (\ref{eq:4.4a'}) should be analytic along bounday ${\bm c}_{1}$ as well. Hence, Eq.~(\ref{eq:4.4a'}) analytically validates the foresight in Eq.~(\ref{eq:3.3}) that $\varPhi(\zeta)$ contains multiplier $(1-\zeta)^{2}$ and is analytic within region ${\bm \omega}$. Substituting Eqs.~(\ref{eq:4.4b}), (\ref{eq:4.4a'}), and (\ref{eq:2.6b}) into Eq.~(\ref{eq:3.7}) yields
\begin{equation}
  \label{eq:4.5}
  \psi^{\prime}(\zeta) = z^{\prime}(\zeta) \varPsi(\zeta) = \sum\limits_{k=-\infty}^{\infty} {\rm i} \left[ \frac{k+1}{2} (A_{k-1} - A_{k+1}) - B_{-k-2} \right] \zeta^{k}
\end{equation}

\subsection{Static equilibrium and displacement single-valuedness}
\label{sec:solution-riemann-hilbert-2}

Static equilibrium between the constrained ground surface and the unbalanced resultant along tunnel periphery should be satisfied in the proposed mechanical model. In the physical plane, the static equilibrium can be expressed as:
\begin{equation}
  \label{eq:4.6}
  \int_{{\bm C}_{11}} \left[ X_{o}(T)+{\rm i}Y_{o}(T) \right]|{\rm d}T| = - \varointclockwise_{{\bm C}_{2}} \left[ X_{i}(S)+{\rm i}Y_{i}(S) \right]|{\rm d}S|
\end{equation}
where $|{\rm d}T|$ and $|{\rm d}S|$ denote length increments along boundaries ${\bm C}_{11}$ and ${\bm C}_{2}$ in the physical plane, respectively. The integral path for ${\bm C}_{11}$ is from point $T_{1}$ along ground surface towards infinity, then from infinity along ground surface towards point $T_{2}$, so that the geomaterial always lies on the left side of the integral path. Considering the zero traction along boundary ${\bm C}_{12}$ in Eq.~(\ref{eq:2.12b}), the left-hand side of Eq.~(\ref{eq:4.6}) can be modified as
\begin{equation}
  \label{eq:4.7}
  \int_{{\bm C}_{11}} \left[ X_{o}(T)+{\rm i}Y_{o}(T) \right]|{\rm d}T| = \int_{{\bm C}_{1}} \left[ X_{o}(T)+{\rm i}Y_{o}(T) \right]|{\rm d}T| = \ointctrclockwise_{{\bm c}_{1}} {\rm e}^{{\rm i}\theta}z^{\prime}({\rm e}^{{\rm i}\theta})[\sigma_{\rho}({\rm e}^{{\rm i}\theta})+{\rm i}\tau_{\rho\theta}({\rm e}^{{\rm i}\theta})]{\rm d}\theta
\end{equation}
where $ |{\rm d}T| = |z^{\prime}({\rm e}^{{\rm i}\theta})|\cdot |{\rm i}{\rm e}^{{\rm i}\theta}| \cdot |{\rm d}\theta| = |z^{\prime}({\rm e}^{{\rm i}\theta})|{\rm d}\theta $ for counter-clockwise length increment in the mapping plane. Substituting Eqs.~(\ref{eq:3.16a}) and (\ref{eq:4.4}) into Eq.~(\ref{eq:4.7}) yields
\begin{equation}
  \label{eq:4.8}
  \ointctrclockwise \left(\sum\limits_{k=-\infty}^{\infty} A_{k}{\rm e}^{{\rm i}k\theta} - \sum\limits_{k=-\infty}^{\infty} B_{k}{\rm e}^{{\rm i}k\theta}\right) \cdot {\rm i}{\rm e}^{{\rm i}\theta} {\rm d}\theta = 2\pi{\rm i}(A_{-1}-B_{-1})
\end{equation}
On the other hand, according to Eq.~(\ref{eq:2.4}), the right-hand side of Eq.~(\ref{eq:4.6}) can be written as
\begin{equation}
  \label{eq:4.9}
  -\varointclockwise_{{\bm C}_{2}} \left[ X_{i}(S)+{\rm i}Y_{i}(S) \right]|{\rm d}S| = - \varointclockwise_{{\bm C}_{2}} X_{i}(S)|{\rm d}S| - {\rm i}\varointclockwise_{{\bm C}_{2}} Y_{i}(S)|{\rm d}S| = - F_{x} - {\rm i} F_{y} = - {\rm i}\gamma\pi R^{2}
\end{equation}
The equilibrium of Eqs.~(\ref{eq:4.8}) and (\ref{eq:4.9}) gives
\begin{equation}
  \label{eq:4.10}
  A_{-1}-B_{-1} = - \frac{\gamma R^{2}}{2}
\end{equation}

In the physical plane, the displacement single-valuedness of geomaterial owing to the displacement boundary in Eq.~(\ref{eq:2.12a}) should be guaranteed and verified. Substituting Eqs.~(\ref{eq:4.4a}), (\ref{eq:4.4a'}), and (\ref{eq:4.5}) into Eq.~(\ref{eq:3.1c}) yields
\begin{equation}
  \label{eq:4.11}
  \begin{aligned}
    g(\zeta) = & \; \kappa \int z^{\prime}(\zeta)\varPhi(\zeta) {\rm d}\zeta - z(\zeta) \overline{\varPhi(\zeta)} - \int \overline{z^{\prime}(\zeta)\varPsi(\zeta)} {\rm d}\overline{\zeta} + {\rm i} C_{0} \\
    = & \; {\rm i}(\kappa A_{-1} + B_{-1}) {\rm Ln}{\rm e}^{{\rm i}\theta} + {\rm i}(\kappa A_{-1} - B_{-1})\ln\rho + {\rm single \; valued \; items}
  \end{aligned}
\end{equation}
where ${\rm Ln}$ denotes the multi-valued natural logarithm sign, $C_{0}$ denotes the integral constant, which should be real due to symmetry, to ensure $g(1)=0$ to satisfy the displacement boundary condition in Eq.~(\ref{eq:2.12a}). The displacement single-valuedness of Eq.~(\ref{eq:4.11}) requires that the multi-valued natural logarithm item should vanish. Thus, we have
\begin{equation}
  \label{eq:4.12}
  \kappa A_{-1} + B_{-1} = 0
\end{equation}
Eq.~(\ref{eq:4.12}) ensures displacement single-valuedness of the geomaterial. Eqs.~(\ref{eq:4.10}) and (\ref{eq:4.12}) give
\begin{equation}
  \label{eq:4.13}
  \left\{
    \begin{aligned}
      A_{-1} = & \; \frac{-\gamma R^{2}}{2(1+\kappa)} \\
      B_{-1} = & \; \frac{\kappa\gamma R^{2}}{2(1+\kappa)} \\
    \end{aligned}
  \right.
\end{equation}
The deduction above indicates that Eqs.~(\ref{eq:4.10}) and (\ref{eq:4.12}) are free from conformal mapping. Such a result is reasonable, since the resultant equilibrium and displacement single-valuedness should always be satisfied for arbitrary doubly-connected region of the mixed boundary conditions in Eq.~(\ref{eq:2.12}). The integral constant $C_{0}$ could be determined in a later stage.

\subsection{Traction boundary condition along tunnel periphery}
\label{sec:solution-riemann-hilbert-3}

Traction boundary condition in Eq.~(\ref{eq:2.12c}) should be also satisfied to uniquely determine the unknown coefficients $d_{n}$ in Eq.~(\ref{eq:4.1}). However, the original form in Eq.~(\ref{eq:2.12c}) would cause computational difficulties, and corresponding path integral form can be used instead. To use the path integral along the boundary ${\bm C}_{2}$ in the physical plane, as well as along the boundary ${\bm c}_{2}$ in the mapping plane, we denote the starting point by $A$ and an arbitrary integral point by $B$ in the physical plane, and the mapping points of these two points in the mapping plane are respectively denoted by $a$ and $b$. Note that point $B$ is always clockwise to point $A$ in both the physical and mapping planes, so that the integral path would always keep the geomaterial on the left hand.

The first equilibrium in Eq.~(\ref{eq:2.12c}) can be equivalently modified in the following path integral form as
\begin{equation}
  \label{eq:4.14}
  - {\rm i} \int_{A}^{B} [X_{i}(S)+{\rm i}Y_{i}(S)]\cdot|{\rm d}S| = - {\rm i} \int_{a}^{b} {\rm e}^{{\rm i}\theta} \frac{z^{\prime}(s)}{|z^{\prime}(s)|} [\sigma_{\rho}(s)+{\rm i}\tau_{\rho\theta}(s)] \cdot |z^{\prime}(s)| \cdot |{\rm d}s| = \int_{a}^{b} z^{\prime}(s)[\sigma_{\rho}(s)+{\rm i}\tau_{\rho\theta}(s)] {\rm d}{s}
\end{equation}
where $ |{\rm d}s| = |r{\rm d}{\rm e}^{{\rm i}\theta}| = r|{\rm d}\theta| = - r{\rm d}\theta = - {\rm d}s $ for clockwise length increment in mapping plane. Simultaneously, the second equilibrium in Eq.~(\ref{eq:2.12c}) can be modified in the following integral form as
\begin{equation}
  \label{eq:4.15}
  - {\rm i} \int_{A}^{B} [X_{i}(S)+{\rm i}Y_{i}(S)]\cdot|{\rm d}S| = {\rm i} \int_{A}^{B} \left( - k_{0}\gamma y\frac{{\rm d}y}{{\rm d}S} + {\rm i}\gamma y\frac{{\rm d}x}{{\rm d}S}\right){\rm d}S = - {\rm i} k_{0}\gamma \int_{A}^{B} y{\rm d}y - \gamma \int_{A}^{B} y{\rm d}x
\end{equation}
where $|{\rm d}S| = R|{\rm d}\vartheta| = - R{\rm d}\vartheta = - {\rm d}S$ for clockwise length increment in physical plane, similar to Eq.~(\ref{eq:2.4}).

As long as the difference between the indefinite integrals in Eqs.~(\ref {eq:4.14}) and (\ref {eq:4.15}) is a constant, Eq.~(\ref {eq:2.12c}) would be satisfied. Both Eqs.~(\ref{eq:4.14}) and (\ref{eq:4.15}) should be prepared into rational series in the mapping plane. Substituting Eq.~(\ref{eq:3.1b}) into Eq.~(\ref{eq:4.14}) yields
\begin{equation}
  \label{eq:4.14a}
  \tag{4.14a}
  \begin{aligned}
    \int_{a}^{b} z^{\prime}(s)[\sigma_{\rho}(s)+{\rm i}\tau_{\rho\theta}(s)] {\rm d}{s} 
    = & \; \int_{a}^{b} \varphi^{\prime}(s){\rm d}s + \frac{z(s)}{\overline{z^{\prime}(s)}}\overline{\varphi^{\prime}(s)} + \int_{a}^{b} \overline{\psi^{\prime}(s)}{\rm d}\overline{s} \\
    = & \; \frac{1}{-2{\rm i}(1-r{\rm e}^{{\rm i}\theta})}\sum\limits_{k=-\infty}^{\infty}\left[
        \begin{aligned}
          & r^{-k+2}A_{-k+1}+(1-2r^{2})r^{-k}A_{-k} \\
          + &(r^{2}-2)r^{-k}A_{-k-1}+r^{-k}A_{-k-2} \\
        \end{aligned}
        \right]{\rm e}^{{\rm i}k\theta} \\
    + & \; {\rm i}\sum\limits_{k=1}^{\infty} \left[\frac{r^{k}}{k}A_{k-1}-\frac{r^{-k}}{2}(A_{-k-2}-A_{-k})+\frac{r^{-k}}{-k}B_{k-1}\right]{\rm e}^{{\rm i}k\theta} \\
    + & \; {\rm i}\sum\limits_{k=1}^{\infty} \left[\frac{r^{-k}}{-k}A_{-k-1}-\frac{r^{k}}{2}(A_{k-2}-A_{k})+\frac{r^{k}}{k}B_{-k-1}\right]{\rm e}^{-{\rm i}k\theta} \\
    + & \; {\rm i}(A_{-1}+B_{-1})\ln{r} + {\rm i}(A_{-1}-B_{-1}){\rm Ln}{\rm e}^{{\rm i}\theta} + {\rm i} C_{a}
  \end{aligned}
\end{equation}
where $C_{a}$ denotes the integral constant to be determined and should be real due to symmetry. Considering the backward conformal mapping in Eq.~(\ref{eq:2.10b}), the real variables $x$ and $y$ along boundary ${\bm C}_{2}$ in Eq.~(\ref{eq:4.15}) can be written as
\begin{subequations}
  \label{eq:4.16}
  \begin{equation}
    \label{4.16a}
    x = - \frac{{\rm i}a}{2}\left(\frac{1+r{\rm e}^{{\rm i}\theta}}{1-r{\rm e}^{{\rm i}\theta}} - \frac{{\rm e}^{{\rm i}\theta}+r}{{\rm e}^{{\rm i}\theta}-r}\right)
  \end{equation}
  \begin{equation}
    \label{4.16b}
    y = - \frac{a}{2}\left(\frac{1+r{\rm e}^{{\rm i}\theta}}{1-r{\rm e}^{{\rm i}\theta}} + \frac{{\rm e}^{{\rm i}\theta}+r}{{\rm e}^{{\rm i}\theta}-r}\right)
  \end{equation}
  \begin{equation}
    \label{eq:4.16c}
    {\rm d}x = -{\rm i}ar\left[\frac{1}{(1-r{\rm e}^{{\rm i}\theta})^{2}}+\frac{1}{({\rm e}^{{\rm i}\theta}-r)^{2}}\right]{\rm d}{\rm e}^{{\rm i}\theta}
  \end{equation}
\end{subequations}
Substituting Eq.~(\ref{eq:4.16}) into Eq.~(\ref{eq:4.15}) with notation of $\frac{ar}{1-r^{2}}=\frac{R}{2}$ yields
\begin{equation}
  \label{eq:4.15a}
  \tag{4.15a}
  \begin{aligned}
    - {\rm i} k_{0}\gamma \int_{A}^{B} y{\rm d}y - \gamma \int_{A}^{B} y{\rm d}x 
    = & - \frac{{\rm i}k_{0}\gamma a^{2}(1-r^{2})^{2}{\rm e}^{2{\rm i}\theta}}{2(1-r{\rm e}^{{\rm i}\theta})^{2}({\rm e}^{{\rm i}\theta}-r)^{2}} - \frac{{\rm i}\gamma a^{2}}{2(1-r{\rm e}^{{\rm i}\theta})^{2}} + \frac{{\rm i}\gamma a^{2} r^{2}}{2({\rm e}^{{\rm i}\theta}-r)^{2}} \\
      & + \frac{{\rm i}\gamma a R}{2({\rm e}^{{\rm i}\theta}-r)} - \frac{{\rm i}\gamma a r R}{2(1-r{\rm e}^{{\rm i}\theta})} + \frac{{\rm i}\gamma R^{2}}{2}\ln\frac{1-r{\rm e}^{{\rm i}\theta}}{1-r{\rm e}^{-{\rm i}\theta}} - \frac{{\rm i}\gamma R^{2}}{2}{\rm Ln}{\rm e}^{{\rm i}\theta}
  \end{aligned}
\end{equation}
According to Eqs.~(\ref{eq:4.14}) and (\ref{eq:4.15}), Eqs.~(\ref{eq:4.14a}) and (\ref{eq:4.15a}) should be equal. To guarantee such a requirement, the multi-valued components of $ {\rm Ln}{\rm e}^{{\rm i}\theta} $ should be eliminated simultaneously, and we again obtain Eq.~(\ref{eq:4.10}). The deduction above further analytically emphasizes the mechanical fact that stress and traction in the geomaterial should be single-valued.

Eqs.~(\ref {eq:4.14a}) and (\ref {eq:4.15a}) should be equal, however, owing to the denominator $-2{\rm i}(1-r{\rm e}^{{\rm i}\theta})$, the equilibrium is not in rational form to facilitate coefficient comparisons, thus, we should modify the equilibrium. Eqs.~(\ref{eq:4.14a}) can be modified as
\begin{equation}
  \label{eq:4.14b}
  \tag{4.14b}
  \begin{aligned}
    & - 2{\rm i}(1-r{\rm e}^{{\rm i}\theta})\left[\int_{a}^{b} \varphi^{\prime}(s){\rm d}s + \frac{z(s)}{\overline{z^{\prime}(s)}}\overline{\varphi^{\prime}(s)} + \int_{a}^{b} \overline{\psi^{\prime}(s)}{\rm d}\overline{s}\right] \\
    = & \; \sum\limits_{k=1}^{\infty}\left[\frac{2r^{-k}}{-k}A_{-k-1}-\frac{2r^{-k}}{-k-1}A_{-k-2}+\frac{2r^{k}}{k}B_{-k-1}-\frac{2r^{k+2}}{k+1}B_{-k-2}+2(1-r^{2})r^{k}(A_{k}-A_{k-1})\right] {\rm e}^{-{\rm i}k\theta} \\
    + & \; \sum\limits_{k=2}^{\infty}\left[\frac{2r^{k}}{k}A_{k-1}-\frac{2r^{k}}{k-1}A_{k-2}+\frac{2r^{-k}}{-k}B_{k-1}-\frac{2r^{-k+2}}{-k+1}B_{k-2}+2(1-r^{2})r^{-k}(A_{-k}-A_{-k-1})\right] {\rm e}^{{\rm i}k\theta} \\
    + & \; \left[2(B_{0}-r^{2}A_{0})+2(1-r^{2})A_{-2}+2(r^{2}-1)A_{-1}+r^{2}\tilde{C}_{a}\right] \left(-r^{-1} {\rm e}^{{\rm i}\theta} \right) \\
    + & \; \left[2(1-r^{2})A_{0}+2(A_{-2}-r^{2}B_{-2})+2(r^{2}-1)A_{-1}+\tilde{C}_{a}\right] \\
  \end{aligned}
\end{equation}
where 
\begin{equation*}
  \tilde{C}_{a} = -(A_{0}-A_{-2})+2(A_{-1}+B_{-1})\ln{r}+2C_{a}
\end{equation*}
The coefficients in Eq.~(\ref{eq:4.14b}) would degenerate to the ones in Refs~\cite{Verruijt1997traction}, when we replace the unbalanced resultant along tunnel periphery by a balanced traction and cancel the far-field displacement constraint along ground surface. The degeneration and comparison details can be found in Appendix~\ref{sec:A}. Similarly, Eq.~(\ref{eq:4.15a}) can be modified as
\begin{equation}
  \label{eq:4.15b}
  \tag{4.15b}
  -2{\rm i}(1-r{\rm e}^{{\rm i}\theta})\left( - {\rm i} k_{0}\gamma \int_{A}^{B} y{\rm d}y - \gamma \int_{A}^{B} y{\rm d}x\right) = \sum\limits_{k=1}^{\infty} E_{-k} {\rm e}^{-{\rm i}k\theta} + \sum\limits_{k=2}^{\infty} E_{k} {\rm e}^{{\rm i}k\theta} + E_{1}{\rm e}^{{\rm i}\theta} + E_{0}
\end{equation}
where $E_{k}$ can be found in Appendix~\ref{sec:B}. The equilibriums of coefficients between Eqs.~(\ref{eq:4.14b}) and (\ref{eq:4.15b}) with some slight modification respectively give
\begin{subequations}
  \label{eq:4.17}
  \begin{equation}
    \label{eq:4.17a}
    A_{-k-1} - A_{-k-2} = - \frac{1}{k+1}A_{-k-2}  + r^{2k}B_{-k-1} - \frac{k}{k+1} r^{2k+2}B_{-k-2} + k(1-r^{2})r^{2k}(A_{k}-A_{k-1}) - \frac{k}{2}r^{k}E_{-k}, \quad k \geq 1
  \end{equation}
  \begin{equation}
    \label{eq:4.17b}
    B_{k-1} = \frac{k}{k-1} r^{2}B_{k-2} + r^{2k}A_{k-1} - \frac{k}{k-1} r^{2k}A_{k-2} + k(1-r^{2})(A_{-k}-A_{-k-1}) - \frac{k}{2}r^{k}E_{k}, \quad k \geq 2
  \end{equation} 
\end{subequations}
\begin{subequations}
  \label{eq:4.18}
  \begin{equation}
    \label{eq:4.18a}
    2B_{0} - 3r^{2}A_{0}+(2-r^{2})A_{-2}+2r^{2}C_{a} = -rE_{1} + 2(1-r^{2})A_{-1}-2r^{2}(A_{-1}+B_{-1})\ln{r}
  \end{equation}
  \begin{equation}
    \label{eq:4.18b}
    (1-2r^{2})A_{0} + 3A_{-2} - 2r^{2}B_{-2} + 2C_{a} = E_{0} + 2(1-r^{2})A_{-1} - 2(A_{-1}+B_{-1})\ln{r}
  \end{equation}
\end{subequations}

\subsection{Solution}
\label{sub:Solution}

With Eqs.~(\ref{eq:4.17}), (\ref {eq:4.18}), and (\ref{eq:4.13}), we can solve $d_{n}$ using an approximate method below. For better convergence, we substitute Eq.~(\ref {eq:4.17a}) into Eq.~(\ref {eq:4.17b}), and Eq.~(\ref {eq:4.17}) can be further modified as
\begin{subequations}
  \label {eq:4.19}
  \begin{equation}
    \label{eq:4.19a}
    \begin{aligned}
      A_{-k}
      = & \; \frac{k-1}{k}A_{-k-1}  + r^{2k-2}B_{-k} - \frac{k-1}{k} r^{2k}B_{-k-1} \\
        & \; + (k-1)(1-r^{2})r^{2k-2}(A_{k-1}-A_{k-2}) - \frac{k-1}{2}r^{k-1}E_{-k+1}, \quad k \geq 2
    \end{aligned}
  \end{equation} 
  \begin{equation}
    \label{eq:4.19b}
    \begin{aligned}
      B_{k} 
      = & \; \frac{k+1}{k} r^{2}B_{k-1} + r^{2k+2}A_{k} - \frac{k+1}{k} r^{2k+2}A_{k-1} - (1-r^{2})A_{-k-2} \\
        & \; + (k+1)(1-r^{2})r^{2k}B_{-k-1} - k(1-r^{2})r^{2k+2}B_{-k-2} \\
        & \; + k(k+1)(1-r^{2})^{2}r^{2k}(A_{k}-A_{k-1}) - \frac{k(k+1)}{2}(1-r^{2})r^{k}E_{-k} - \frac{k+1}{2}r^{k+1}E_{k+1}, \quad k \geq 1
    \end{aligned}
  \end{equation}
\end{subequations}
Using Eq.~(\ref {eq:4.4}), the left-hand sides of Eq.~(\ref{eq:4.19}) can be expanded as
\begin{subequations}
  \label{eq:4.20}
  \begin{equation}
    \label{eq:4.20a}
    \begin{aligned}
      \sum\limits_{n=k}^{\infty} \alpha_{n-k} d_{-n}
      = & \; \frac{k-1}{k}A_{-k-1} + r^{2k-2}B_{-k} - \frac{k-1}{k} r^{2k}B_{-k-1} \\
        & \; + (k-1)(1-r^{2})r^{2k-2}(A_{k-1}-A_{k-2}) - \frac{k-1}{2}r^{k-1}E_{-k+1}, \quad k \geq 2
    \end{aligned}
  \end{equation} 
  \begin{equation}
    \label{eq:4.20b}
    \begin{aligned}
      \sum\limits_{n=k+1}^{\infty} \beta_{n-k} d_{n}
      = & \; \frac{k+1}{k} r^{2}B_{k-1} + r^{2k+2}A_{k} - \frac{k+1}{k} r^{2k+2}A_{k-1} - (1-r^{2})A_{-k-2} \\
        & \; + (k+1)(1-r^{2})r^{2k}B_{-k-1} - k(1-r^{2})r^{2k+2}B_{-k-2} \\
        & \; + k(k+1)(1-r^{2})^{2}r^{2k}(A_{k}-A_{k-1}) - \frac{k(k+1)}{2}(1-r^{2})r^{k}E_{-k} - \frac{k+1}{2}r^{k+1}E_{k+1}, \quad k \geq 1
    \end{aligned}
  \end{equation}
\end{subequations}
Eqs.~(\ref{eq:4.20a}) and (\ref{eq:4.20b}) constrain $ d_{-n} (n \geq 2) $ and $ d_{n} (n \geq 2) $, respectively, while $ d_{-1} $, $ d_{0} $, and $ d_{1} $ are not constrained yet. Thus, Eqs.~(\ref {eq:4.18}) and (\ref {eq:4.13}) should be used and can be expanded as
\begin{subequations}
  \label {eq:4.21} 
  \begin{equation}
    \label {eq:4.21a}
    I_{-1} d_{-1} + I_{0} d_{0} + I_{1} d_{1} + 2r^{2}C_{a} = - \sum_{n=2}^{\infty} I_{-n} d_{-n} - \sum_{n=2}^{\infty} I_{n} d_{n} + I^{ \prime }
  \end{equation}
  \begin{equation}
    \label {eq:4.21b}
    J_{-1} d_{-1} + J_{0} d_{0} + J_{1} d_{1} + 2C_{a} = - \sum_{n=2}^{\infty} J_{-n} d_{-n} - \sum_{n=2}^{\infty} J_{n} d_{n} + J^{ \prime }
  \end{equation}
  \begin{equation}
    \label {eq:4.21c}
    \alpha_{0}d_{-1} = - \sum\limits_{n=2}^{\infty} \alpha_{-1+n}d_{-n} + \frac{-\gamma R^{2}}{2(1+\kappa)}
  \end{equation}
  \begin{equation}
    \label {eq:4.21d}
    \beta_{1}d_{0} + \beta_{2}d_{1} = - \sum\limits_{n=2}^{\infty} \beta_{n+1}d_{n} + \frac{\kappa \gamma R^{2}}{2(1+\kappa)}
  \end{equation}
\end{subequations}
where
\begin{equation*}
  \left\{
    \begin{aligned}
      I_{n} = & \; 2\beta_{n}, \quad n \geq 1 \\
      I_{0} = & \; - 3r^{2}\alpha_{0} \\
      I_{-1} = & \; - 3r^{2}\alpha_{1} \\
      I_{-n} = & \; \left( 2-r^{2} \right) \alpha_{-2+n} - 3r^{2}\alpha_{n}, \quad n \geq 2 \\
      I^{ \prime } = & \; \frac{-\gamma R^{2}}{1+\kappa} \left( 1-r^{2} \right) + \frac{ 1 - \kappa }{ 1 + \kappa } \gamma R^{2} r^{2}\ln r - rE_{1}
    \end{aligned}
  \right.
\end{equation*}
\begin{equation*}
  \left\{
    \begin{aligned}
      J_{n} = & \; -2r^{2}\beta_{n+2}, \quad n \geq 1 \\
      J_{0} = & \; \left( 1-2r^{2} \right) \alpha_{0} - 2r^{2}\beta_{2} \\
      J_{-1} = & \; \left( 1-2r^{2} \right) \alpha_{1} - 2r^{2}\beta_{1} \\
      J_{-n} = & \; \left( 1-2r^{2} \right) \alpha_{n} + 3\alpha_{-2+n}, \quad n \geq 2 \\
      J^{ \prime } = & \; \frac{-\gamma R^{2}}{1+\kappa} \left( 1-r^{2} \right) + \frac{1-\kappa}{1+\kappa}\gamma R^{2}\ln{r} + E_{0}
    \end{aligned}
  \right.
\end{equation*}
Eq.~(\ref {eq:4.21}) contains four linear equations on four variables ($ d_{-1} $, $ d_{0} $, $ d_{1} $, and $ C_{a} $). Thus, the linear system on $ d_{k} $ and $ C_{a} $ composed of Eqs.~(\ref {eq:4.20}) and (\ref {eq:4.21}) are definite, as long as the values of $ A_{k} $ and $ B_{k} $ are determined.

Assume that $d_{n}$ and $C_{a}$ can be written as
\begin{subequations}
  \label{eq:4.22}
  \begin{equation}
    \label{eq:4.22a}
    d_{n} = \sum\limits_{q=0}^{\infty} d_{n}^{(q)}, \quad n \in {\bm Z}
  \end{equation}
  \begin{equation}
    \label{eq:4.22b}
    C_{a} = \sum\limits_{q=0}^{\infty} C_{a}^{(q)}
  \end{equation}
\end{subequations}
Then the initial values of $ d_{-n}^{(0)} (n \geq 2) $ and $ d_{n}^{(0)} (n \geq 2) $ can be determined by the following linear system based on Eq.~(\ref{eq:4.20}) as
\begin{subequations}
  \label{eq:4.23}
  \begin{equation}
    \label{eq:4.23a}
    \sum\limits_{n=k}^{\infty} \alpha_{-k+n} d_{-n}^{(0)} = - \frac{k-1}{2}r^{k-1}E_{-k+1}, \quad k \geq 2
  \end{equation}
  \begin{equation}
    \label{eq:4.23b}
    \sum\limits_{n=k+1}^{\infty} \beta_{n-k} d_{n}^{(0)} = - \frac{k(k+1)}{2}(1-r^{2})r^{k}E_{-k} - \frac{k+1}{2}r^{k+1}E_{k+1}, \quad k \geq 1
  \end{equation}
\end{subequations}
Then the initial value of $ d_{-1}^{(0)} $, $ d_{0}^{(0)} $, $ d_{1}^{(0)} $, and $ C_{a}^{(0)} $ can be determined by the following linear system based on Eq.~(\ref{eq:4.21}) as
\begin{subequations}
  \label {eq:4.24} 
  \begin{equation}
    \label {eq:4.24a}
    I_{-1} d_{-1}^{(0)} + I_{0} d_{0}^{(0)} + I_{1} d_{1}^{(0)} + 2r^{2}C_{a}^{(0)} = - \sum_{n=2}^{\infty} I_{-n} d_{-n}^{(0)} - \sum_{n=2}^{\infty} I_{n} d_{n}^{(0)} + I^{ \prime }
  \end{equation}
  \begin{equation}
    \label {eq:4.24b}
    J_{-1} d_{-1}^{(0)} + J_{0} d_{0}^{(0)} + J_{1} d_{1}^{(0)} + 2C_{a}^{(0)} = - \sum_{n=2}^{\infty} J_{-n} d_{-n}^{(0)} - \sum_{n=2}^{\infty} J_{n} d_{n}^{(0)} + J^{ \prime }
  \end{equation}
  \begin{equation}
    \label {eq:4.24c}
    \alpha_{0}d_{-1}^{(0)} = - \sum\limits_{n=2}^{\infty} \alpha_{-1+n}d_{-n}^{(0)} + \frac{-\gamma R^{2}}{2(1+\kappa)}
  \end{equation}
  \begin{equation}
    \label {eq:4.24d}
    \beta_{1}d_{0}^{(0)} + \beta_{2}d_{1}^{(0)} = - \sum\limits_{n=2}^{\infty} \beta_{n+1}d_{n}^{(0)} + \frac{\kappa \gamma R^{2}}{2(1+\kappa)}
  \end{equation}
\end{subequations}
Eq.~(\ref{eq:4.24}) indicates that $ d_{-1}^{(0)} $, $ d_{0}^{(0)} $, $ d_{1}^{(0)} $, and $ C_{a}^{(0)} $ are dependent on $ d_{-n}^{(0)} (n \geq 2) $ and $ d_{n}^{(0)} (n \geq 2) $ solved via Eq.~(\ref{eq:4.23}).

With the initial values in Eqs.~(\ref{eq:4.23}) and (\ref{eq:4.24}), for iteration rep $ q \geq 0 $, the coefficients of $ \varphi(\zeta) $ in Eq.~(\ref{eq:4.4}) can be computed as
\begin{subequations}
  \label{eq:4.25}
  \begin{equation}
    \label{eq:4.25a}
    A_{k}^{(q)} = \sum\limits_{n=-\infty}^{k} \alpha_{k-n}d_{n}^{(q)}
  \end{equation}
  \begin{equation}
    \label{eq:4.25b}
    B_{k}^{(q)} = \sum\limits_{n=k+1}^{\infty} \beta_{n-k}d_{n}^{(q)}
  \end{equation}
\end{subequations}
Then for the next iteration rep $q+1$, $ d_{-n}^{(q+1)} (q \geq 0, n \geq 2) $ and $ d_{n}^{(q+1)} (q \geq 0, n \geq 2) $ can be determined via the following linear system based on Eq.~(\ref{eq:4.20}) as
\begin{subequations}
  \label{eq:4.26}
  \begin{equation}
    \label{eq:4.26a}
    \begin{aligned}
      \sum\limits_{n=k}^{\infty} \alpha_{n-k} d_{-n}^{(q+1)}
      = & \; \frac{k-1}{k}A_{-k-1}^{(q)} + r^{2k-2}B_{-k}^{(q)} - \frac{k-1}{k} r^{2k}B_{-k-1}^{(q)} \\
        & \; + (k-1)(1-r^{2})r^{2k-2}(A_{k-1}^{(q)}-A_{k-2}^{(q)}), \quad k \geq 2
    \end{aligned}
  \end{equation}
  \begin{equation}
    \label{eq:4.26b}
    \begin{aligned}
      \sum\limits_{n=k+1}^{\infty} \beta_{n-k} d_{n}^{(q+1)}
      = & \; \frac{k+1}{k}r^{2}B_{k-1}^{(q)} + r^{2k+2}A_{k}^{(q)} - \frac{k+1}{k} r^{2k+2}A_{k-1}^{(q)} - (1-r^{2})A_{-k-2}^{(q)} \\
        & \; + (k+1)(1-r^{2})r^{2k}B_{-k-1}^{(q)} - k(1-r^{2})r^{2k+2}B_{-k-2}^{(q)} \\
      & \; + k(k+1)(1-r^{2})^{2}r^{2k}(A_{k}^{(q)} - A_{k-1}^{(q)}), \quad k \geq 1
    \end{aligned}
  \end{equation}
\end{subequations}
Subsequently, $ d_{-1}^{(q+1)}$, $ d_{0}^{(q+1)}$, $ d_{1}^{(q+1)}$, and $ C_{a}^{(q+1)} (q \geq 0) $ can be determined via the following linear system based on Eq.~(\ref{eq:4.21}) as
\begin{subequations}
  \label {eq:4.27} 
  \begin{equation}
    \label {eq:4.27a}
    I_{-1} d_{-1}^{(q+1)} + I_{0} d_{0}^{(q+1)} + I_{1} d_{1}^{(q+1)} + 2r^{2}C_{a}^{(q+1)} = - \sum_{n=2}^{\infty} I_{-n} d_{-n}^{(q+1)} - \sum_{n=2}^{\infty} I_{n} d_{n}^{(q+1)}
  \end{equation}
  \begin{equation}
    \label {eq:4.27b}
    J_{-1} d_{-1}^{(q+1)} + J_{0} d_{0}^{(q+1)} + J_{1} d_{1}^{(q+1)} + 2C_{a}^{(q+1)} = - \sum_{n=2}^{\infty} J_{-n} d_{-n}^{(q+1)} - \sum_{n=2}^{\infty} J_{n} d_{n}^{(q+1)}
  \end{equation}
  \begin{equation}
    \label {eq:4.27c}
    \alpha_{0}d_{-1}^{(q+1)} = - \sum\limits_{n=2}^{\infty} \alpha_{-1+n}d_{-n}^{(q+1)}
  \end{equation}
  \begin{equation}
    \label {eq:4.27d}
    \beta_{1}d_{0}^{(q+1)} + \beta_{2}d_{1}^{(q+1)} = - \sum\limits_{n=2}^{\infty} \beta_{n+1}d_{n}^{(q+1)}
  \end{equation}
\end{subequations}
Eq.~(\ref{eq:4.27}) indicates that $ d_{-1}^{(q+1)} $, $ d_{0}^{(q+1)} $, $ d_{1}^{(q+1)} $, and $ C_{a}^{(q+1)} (q \geq 0) $ are dependent on $ d_{-n}^{(q+1)} (n \geq 2, q \geq 0) $ and $ d_{n}^{(q+1)} (n \geq 2, q \geq 0) $ solved in Eq.~(\ref{eq:4.26}). Then we can set $ q:=q+1 $ to continue the iteration, until 
\begin{equation*}
  \max|d_{n}^{(q+1)}| \leq \varepsilon
\end{equation*}
is reached, where $ \varepsilon $ is a small numeric, $ \varepsilon = 10^{-16} $ for instance, which is the default machine accuracy of double precision of the programming code {\tt{fortran}}.

\subsection{Solution convergence and numerical truncation}
\label{sub:Numerical-convergence}

The approximate solution consists of the initial value determination phase in Eqs.~(\ref {eq:4.23}) and (\ref {eq:4.24}) and the iteration phase in Eqs.~(\ref {eq:4.26}) and (\ref {eq:4.27}). Obviously, the latter determines the solution convergence. As long as the absolute values of the coefficients in front of $ A_{k} $ and $ B_{k} $ in Eq.~(\ref {eq:4.26}) are less than {\tt{1}}, the magnitudes of $ d_{n}^{(q)} $ would approach zero as iteration procedes ($q$ increases). Henceforth, we should examine the value ranges of these coefficients in front of $ A_{k} $ and $ B_{k} $ in Eq.~(\ref {eq:4.26}).

It is common in shallow tunnel engineering that tunnel buried depth $h$ is 2 times larger than tunnel radius $R$ ($ h \geq 2R$), which leads to $ r \leq (2+\sqrt{3})^{-1} $ according to Eq.~(\ref{eq:2.11}). It can be easily verified that such a range of $ r $ would guarantee that the absolute values of all coefficients in front of $ A_{k} $ and $ B_{k} $ in Eq.~(\ref{eq:4.26}) are less than {\tt{1}}, and the computation is so plain that it is not necessary to be expanded here. Therefore, such a range of $ r $ would further guarantee the convergence of the iteration procedure in Eqs.~(\ref{eq:4.23})-(\ref{eq:4.27}). In other words, in the iteration procedure in Eqs.~(\ref{eq:4.23})-(\ref{eq:4.27}), as $q$ increases, the magnitudes of $ d_{n}^{(q)} $ would approach zero.

To obtain numerical solution in practical computation, we have to truncate the infinite bilaterial series of $ d_{n} $ in Eq.~(\ref {eq:4.1}) into $ 2N+1 $ items ($ -N \leq n \leq N $), and the coefficients in Eq.~(\ref {eq:4.4}) would correspondingly turn to the following form:
\begin{equation*}
  \left\{ 
    \begin{aligned}
      A_{k} = & \; \sum_{n=-N}^{k} \alpha_{k-n} d_{n} \\
      B_{k} = & \; \sum_{n=k+1}^{N} \beta_{n-k} d_{n} \\
    \end{aligned}
  \right. 
\end{equation*}
Subsequently, all the equations in the Section~\ref {sub:Solution} would turn to finite, as well as the $ E_{k} $ series in Eq.~(\ref {eq:4.15b}). No matter in the initial value determination phase in Eqs.~(\ref {eq:4.23}) and  (\ref {eq:4.24}), or in the iteration phase in Eqs.~(\ref {eq:4.26}) and (\ref {eq:4.27}), three linear equatition system sets can be established: (1) Set 1: Eqs.~(\ref {eq:4.23a}) and (\ref {eq:4.26a}), which independently determine $ d_{-n}^{(q)} (2 \leq n \leq N, q \geq 0)$; (2) Set 2: Eqs.~(\ref {eq:4.23b}) and (\ref {eq:4.26b}), which independently determine $ d_{n}^{(q)} (2 \leq n \leq N, q \geq 0)$; (3) Set 3: Eqs.~(\ref {eq:4.24}) and (\ref {eq:4.27}), which further determine $ d_{1}^{(q)} $, $ d_{0}^{(q)} $, $ d_{-1}^{(q)} $, and $ C_{a}^{(q)} (q \geq 0) $ for completion. The three sets can be written into the following form:
\begin{equation}
  \label{eq:4.28}
  {\bf A}_{l} {\bf \cdot} {\bf x}_{l} = {\bf b}_{l}, \quad l = 1,2,3
\end{equation}
where $ {\bf A}_{l} $, $ {\bf x}_{l} $, and $ {\bf b}_{l} $ respectively denote the coefficient matrix, the variable vector, and the constant vector for sets of linear systems $ l = 1,2,3 $. The coefficient matrices of linear system sets 1 and 2 can be expanded below:
\begin{subequations}
  \label{eq:4.29}
  \begin{equation}
    \label{eq:4.29a}
    {\bf A}_{1} (i, j) = \alpha_{j-i}, \quad 1 \leq i \leq j \leq N-1
  \end{equation}
  \begin{equation}
    \label{eq:4.29b}
    {\bf A}_{2} (i, j) = \beta_{j-i+1}, \quad 1 \leq i \leq j \leq N-1
  \end{equation}
\end{subequations}
The other components of these three linear system sets are plain in Eqs.~(\ref {eq:4.23}), (\ref {eq:4.24}), (\ref {eq:4.26}), and (\ref {eq:4.27}).

With the coefficient matrices in Eq.~(\ref {eq:4.29}), we can explain the reason of applying the approximate solution. Intuitively, it seems that applying the approximate solution above is very indirect and complicated, since the problem is definitely a linear elastic one and could be solved via a more direct solution containing only one single linear system, instead of such a great deal of linear systems piled in iteration. Indeed, such a problem can be solved via one single linear system, but the numerical stability can not be theoretically guaranteed. When dealing with solution of linear system, we should consider numerical stability as well, which is greatly related to condition number of the coefficient matrix. For the direct solution, the condition number would be very large, since the elements would contain both positive and negative powers of $r$, as can be seen in the coefficients in Eq.~(\ref {eq:4.14b}). Then the condition number of the coefficient matrix measured by 2-norm would be very large, and would make the solution of the direct method potentially unstable, and might cause great error. In contrast, Eqs.~(\ref {eq:4.24}), (\ref {eq:4.27}), and (\ref {eq:4.29}) show that the condition numbers of the coefficient matrices $ {\bf A}_{l} (l=1,2,3) $ would be small, which is further verified in the numerical verification and case discussions.

\section{Stress and displacement in geomaterial}
\label{sec:stress-displ-solut}

The solution of $ d_{n} $ gives $ A_{k} $ and $ B_{k} $ in Eq.~(\ref{eq:4.4}) to reach $ \varphi(\zeta) $ in Eq.~(\ref{eq:4.4}) and $ \varPhi(\zeta) $ in Eq.~(\ref{eq:4.4a'}), as well as $ \psi(\zeta) $ in Eq.~(\ref{eq:4.5}). The stress and displacement components within the annulus $ {\bm \omega} $ can be obtained via Eq.~(\ref{eq:3.1}) as
\begin{subequations}
  \label{eq:5.1}
  \begin{equation}
    \label{eq:5.1a}
    \sigma_{\theta}(\zeta) + \sigma_{\rho}(\zeta) = 4 \Re \left[ \frac{(1-\zeta)^{2}}{-2a} \sum\limits_{k=-\infty}^{\infty} \zeta^{k} A_{k} \right], \quad \zeta = \rho \cdot {\rm e}^{{\rm i}\theta} \in {\bm \omega}
  \end{equation}
  \begin{equation}
    \label{eq:5.1b}
    \begin{aligned}
      \sigma_{\rho}(\zeta) + {\rm i} \tau_{\rho\theta}(\zeta) = 
      & \sum\limits_{k=-\infty}^{\infty} \left[ 
        \begin{aligned}
          & \frac{(1-\overline{\zeta})^{2}}{-2a} \overline{\zeta}^{k} + \frac{k+2}{\sigma^{2}}\frac{(1-\zeta)^{2}}{-4a}\overline{\zeta}^{k+1} + \frac{k}{\sigma^{2}} \frac{(1-\zeta)^{2}}{4a} \overline{\zeta}^{k-1} \\
          & + \frac{1-\zeta^{2}}{-2a} \frac{1-\overline{\zeta}}{\sigma^{2}}\overline{\zeta}^{k} + \frac{1-\zeta^{2}}{4a} \frac{(1-\overline{\zeta})^{2}}{\sigma^{2}} k \overline{\zeta}^{k-1} +  \frac{(1-\zeta)^{2}}{-2a} \zeta^{k}
        \end{aligned}
        \right] A_{k} \\
      & +  \sum\limits_{k=-\infty}^{\infty} \frac{(1-\zeta)^{2}}{2a \sigma^{2}} \overline{\zeta}^{-k-2} B_{k}
    \end{aligned}
    , \quad \zeta = \rho \cdot {\rm e}^{{\rm i}\theta} \in {\bm \omega}
  \end{equation}
  \begin{equation}
    \label{eq:5.1c}
    \begin{aligned}
      g(\zeta) = 2G \left[ u(\zeta) + {\rm i} v(\zeta) \right] = 
      & {\rm i} \sum\limits_{k=1}^{\infty} \left\{ \kappa A_{k-1} \frac{\zeta^{k}}{k} + \left[ \frac{1}{2}(A_{k-2}-A_{k})-\frac{B_{-k-1}}{k}\right] \overline{\zeta}^{k}\right\} \\
      & + {\rm i} \sum\limits_{k=1}^{\infty} \left\{ \kappa A_{-k-1} \frac{\zeta^{-k}}{-k} + \left[ \frac{1}{2}(A_{-k-2}-A_{-k})-\frac{B_{k-1}}{-k} \right]\overline{\zeta}^{-k}\right\} \\
      & - {\rm i} \frac{1+\zeta}{2} \frac{(1-\overline{\zeta})^{2}}{1-\zeta} \sum\limits_{k=-\infty}^{\infty} A_{k} \overline{\zeta}^{k} + {\rm i} (\kappa A_{-1} - B_{-1}) \ln \rho + {\rm i} C_{0} \\
    \end{aligned}
    , \quad \zeta = \rho \cdot {\rm e}^{{\rm i}\theta} \in {\bm \omega}
  \end{equation}
\end{subequations}
where $ \sigma_{\rho}(\zeta) $, $ \sigma_{\theta}(\zeta) $, and $ \tau_{\rho\theta}(\zeta) $ denote radial, hoop, and tangential stress components in the unit annulus, respectively; $ u(\zeta) $ and $ v(\zeta) $ denote horizontal and vertical displacement components in the unit annulus, respectively.

When $\rho = 1$, the displacement in Eq.~(\ref{eq:5.1c}) would turn to
\begin{equation}
  \label{eq:5.2}
  \begin{aligned}
    g({\rm e}^{{\rm i}\theta}) = 
    & {\rm i} \sum\limits_{k=1}^{\infty} \left\{ \kappa A_{k-1} \frac{{\rm e}^{{\rm i}k\theta}}{k} + \left[ \frac{1}{2}(A_{k-2}-A_{k})-\frac{B_{-k-1}}{k}\right] {\rm e}^{-{\rm i}k\theta}\right\} \\
    & + {\rm i} \sum\limits_{k=1}^{\infty} \left\{ \kappa A_{-k-1} \frac{{\rm e}^{-{\rm i}k\theta}}{-k} + \left[ \frac{1}{2}(A_{-k-2}-A_{-k})-\frac{B_{k-1}}{-k} \right]{\rm e}^{{\rm i}k\theta}\right\} \\
    & - {\rm i} \frac{1-{\rm e}^{2{\rm i}\theta}}{2{\rm e}^{2{\rm i}\theta}} \sum\limits_{k=-\infty}^{\infty} A_{k} {\rm e}^{-{\rm i}k\theta} + {\rm i} C_{0} \\
  \end{aligned}
  , \quad \theta \in [0, 2\pi)
\end{equation}
Then the undetermined coefficient $ C_{0} $ in Eq.~(\ref{eq:4.11}) can be determined when $\theta = 0$ in Eq.~(\ref{eq:5.2}) as
\begin{equation}
  \label{eq:5.3}
  C_{0} = - \sum\limits_{k=1}^{\infty} \left[ \kappa \frac{A_{k-1}}{k} + \frac{1}{2}(A_{k-2}-A_{k})-\frac{B_{-k-1}}{k}\right] - \sum\limits_{k=1}^{\infty} \left[ \kappa \frac{A_{-k-1}}{-k} + \frac{1}{2}(A_{-k-2}-A_{-k})-\frac{B_{k-1}}{-k}\right]
\end{equation}
Eq.~(\ref{eq:5.2}) indicates that the displacement components are zero, when $\zeta = 1$, instead of infinity in Refs~\cite{Strack2002phdthesis, Strack_Verruijt2002buoyancy, Lu2016, Zengguisen2019, Self2020JEM} mentioned in Setion~\ref{sec:intro}. Till now, we can say that the displacement singularity has been cancelled.

Due to the abrupt change of boundary conditions near polar points $(1, -\theta_{0})$ and $(1, \theta_{0})$ in Eqs.~(\ref{eq:2.12a}) and (\ref{eq:2.12b}), the Gibbs phenomena would occur in Eqs.~(\ref{eq:5.1}) and (\ref{eq:5.2}), and cause oscillations of the stress and displacement components. To weaken the influence of the Gibbs phenomena, the Lanczos filtering technique is applied~\cite[]{Lanczos1956,singh2019simplified,chawde2021mixed}. To be specific, all $ A_{k} $ and $ B_{k} $ in Eqs.~(\ref{eq:5.1}) and (\ref{eq:5.3}) are replaced by $ L_{k} \cdot A_{k} $ and $ L_{k} \cdot B_{k} $, where $ L_{k} $ denote the Lanczos filtering parameters, and can be expressed as
\begin{equation}
  \label{eq:5.4}
  L_{k} =
  \left\{
    \begin{aligned}
      & 1, \quad k=0 \\
      & \sin\left(\frac{k}{N} \pi \right)/\left(\frac{k}{N} \pi\right), \quad {\rm otherwise} \\
    \end{aligned}
  \right.
\end{equation}

With the bidirectional conformal mappings, the stress and displacement components in the unit annulus can be backwardly mapped onto the ones in the geomaterial in the lower half plane:
\begin{subequations}
  \label{eq:5.5} 
  \begin{equation}
    \label{eq:5.5a}
    \left\{
      \begin{aligned}
        & \sigma_{y}(z) + \sigma_{x}(z) = \sigma_{\theta}(\zeta) + \sigma_{\rho}(\zeta) \\
        & \sigma_{y}(z) - \sigma_{x}(z) + 2{\rm i}\tau_{xy}(z) = \left[ \sigma_{\theta}(\zeta) - \sigma_{\rho}(\zeta) + 2{\rm i} \tau_{\rho\theta}(\zeta) \right] \cdot \frac {\overline{\zeta}} {\zeta} \frac {\overline{z^{\prime}(\zeta)}} {z^{\prime}(\zeta)} \\
      \end{aligned}
    \right.
  \end{equation}
  \begin{equation}
    \label{eq:5.5b}
    u(z) + {\rm i}v(z) = u(\zeta) + {\rm i}v(\zeta)
  \end{equation}
\end{subequations}
where $ \sigma_{x}(z) $, $ \sigma_{y}(z) $, and $ \tau_{xy}(z) $ denote the horizontal, vertical, and shear stress components due to excavation in the lower half plane, respectively; $ u(z) $ and $ v(z) $ denote horizontal and vertical displacement components in the lower half plane, respectively. The final stress field within geomaterial is the sum of Eqs.~(\ref{eq:2.1}) and~(\ref{eq:5.5a}). Till now, the stress and displacement distributions in the geomaterial region $ {\bm \varOmega} $ are obtained, and our problem is solved.

\section{Case verification}
\label{sec:numer-verifi-discu}

In this section, we use numerical cases to illustrate the neccessity of Lanczos filtering in Eq.~(\ref{eq:5.4}), the convergence of the proposed solution in this study, and the comparisons between the solution in this study and an existing analytic solution for verification. The parameters in Table~\ref{tab:1} are shared in all the following numerical cases in this study, while the free ground surface length $ x_{0} $ is not included, since it is changable and would greatly affect the solution convergence. The geomaterial is set to be plane-strain ($ \kappa=3-4\nu $). All numerical cases are realized in the programming code {\tt{fortran}} of {\tt{gcc-13.1.1}}. The condition numbers of the coefficient matrices in all numerical cases are computed via {\tt{dgesvd}} package of {\tt{lapack-3.11.0}}, and the linear systems in all numerical cases are solved via {\tt{dgesv}} package of {\tt{lapack-3.11.0}}. All the figures are constructed via {\tt{gnuplot-5.4}}. According to {\tt{GNU General Public License}}, all the source codes in the numerical cases are released in author {\sf{Luobin Lin}}'s github repository {\emph{github.com/luobinlin987/eliminate-shallow-tunnel-displacement-singularity}}.

The condition numbers of all the iterative linear systems in the numerical cases for the propose solution below are all less than {\tt{$10^{2}$}} after computation, indicating that the solutions of the linear systems for the proposed solution in this study are numerically stable. Moreover, the iterative reps of all the numerical cases below are less than {\tt{30}}, indicating that the convergence of the iteration procedure in Eqs.~(\ref{eq:4.22})-(\ref{eq:4.27}) is very fast. The time cost of each numerical case is less than {\emph{0.5 sec}}.

\begin{table}[htb]
  \centering
  \begin{tabular}{ccccccc}
    \toprule
    $ R {\rm (m)} $ & $h/R$ & $ \gamma {\rm (kN/m^{3})}$ & $ k_{0} $ & $ E {\rm (MPa)} $ & $ \nu $ & $N$ \\
    \midrule
    5 & 2 & 20 & 0.8 & 20 & 0.3 & 50 \\
    \bottomrule
  \end{tabular}
  \caption{Parameters in numerical cases}
  \label{tab:1}
\end{table}

\subsection{Lanczos filtering}
\label{sec:lanczos-filtering}

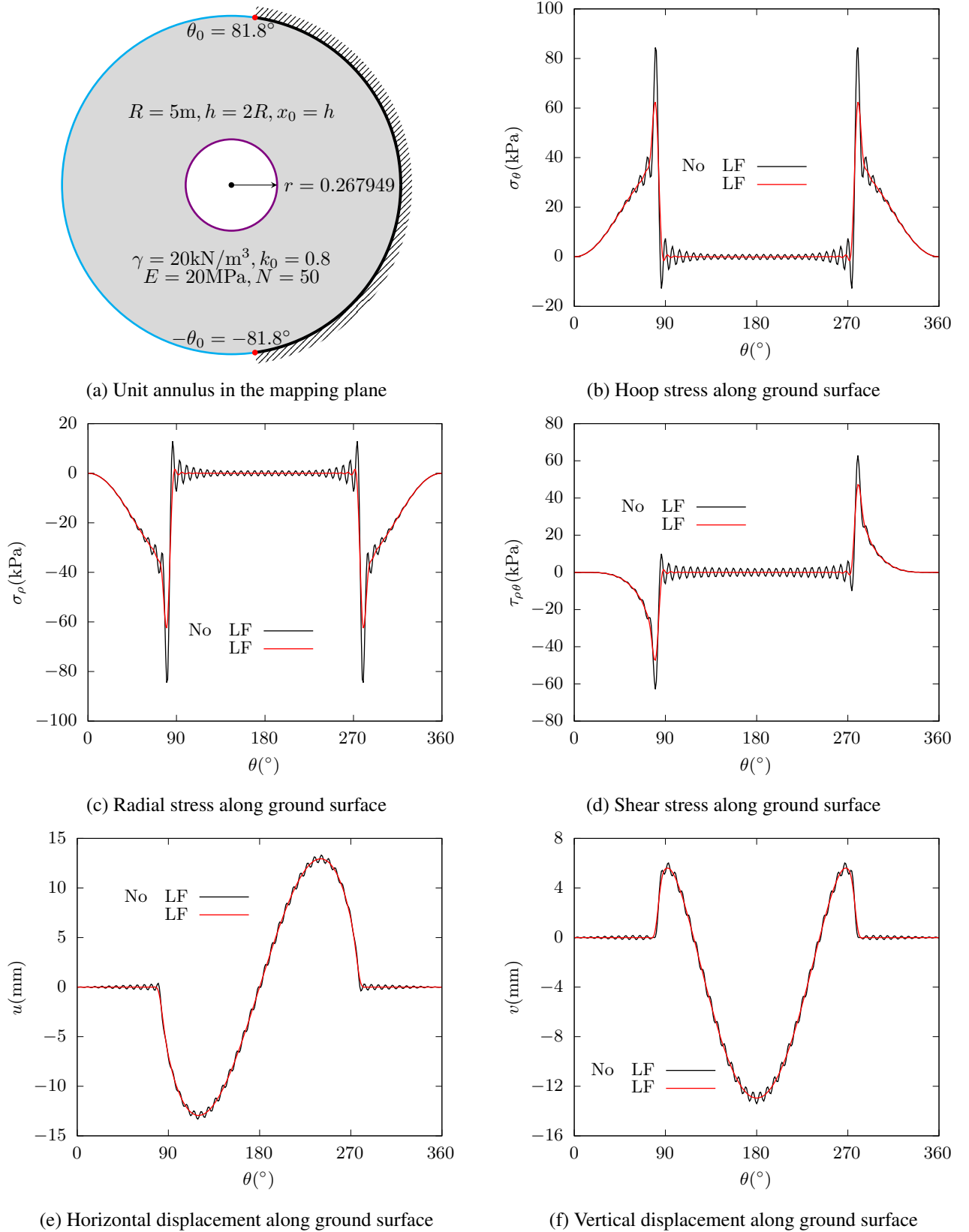
\begin{figure}[htb]
  \centering
  \begin{tabular}{cc}
    \begin{tikzpicture}
      \fill [gray!30] (0,0) circle [radius = 3];
      \fill [white] (0,0) circle [radius = 3*0.27];
      \draw [cyan, line width = 1pt] circle [radius = 3];
      \draw [violet, line width = 1pt] circle [radius = 3*0.27];
      \fill [pattern = north east lines] (3,0) -- (3.2,0) arc [start angle = 0, end angle = 82, radius = 3.2] -- ({3*cos(82)},{3*sin(82)}) arc [start angle = 82, end angle = 0, radius = 3];
      \fill [pattern = north east lines] (3,0) -- (3.2,0) arc [start angle = 0, end angle = -82, radius = 3.2] -- ({3*cos(-82)},{3*sin(-82)}) arc [start angle = -82, end angle = 0, radius = 3];
      \draw [line width = 1.5pt] (3,0) arc [start angle = 0, end angle = 82, radius = 3];
      \draw [line width = 1.5pt] (3,0) arc [start angle = 0, end angle = -82, radius = 3];
      \fill [red] ({3*cos(82)},{3*sin(82)}) circle [radius = 0.05];
      \fill [red] ({3*cos(-82)},{3*sin(-82)}) circle [radius = 0.05];
      \draw [->] (0,0) -- ({3*0.27},0) node [right] {$ r=0.267949 $};
      \fill (0,0) circle [radius = 0.05];
      \node at (0,-3) [above] {$ -\theta_{0} = -81.8^{\circ} $};
      \node at (0,3) [below] {$ \theta_{0} = 81.8^{\circ} $};
      \node at (0,1) [above] {$ R=5{\rm m}, h=2R, x_{0}=h $};
      \node at (0,-1) [below, align=center] {
        $ \gamma=20{\rm kN/m^{3}}, k_{0}=0.8 $ \\
        $ E = 20{\rm MPa}, N = 50 $
      };
    \end{tikzpicture}
    &
      \includegraphics{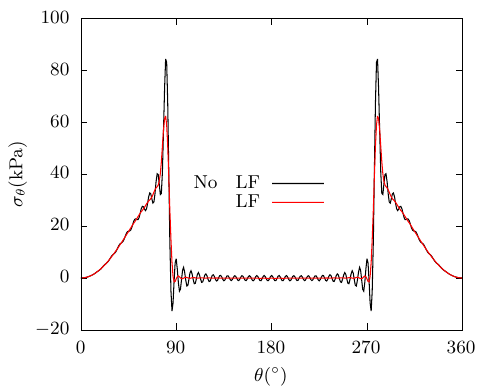}
    \\
    (a) Unit annulus in the mapping plane
    &
      (b) Hoop stress along ground surface
    \\
    \includegraphics{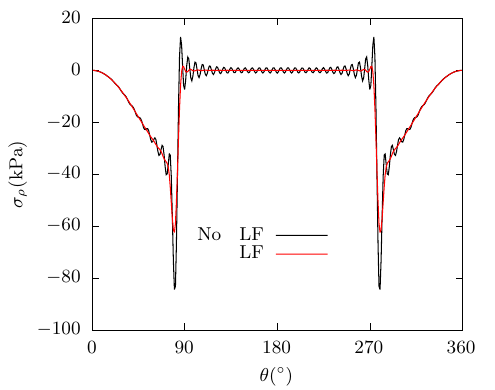}
    &
      \includegraphics{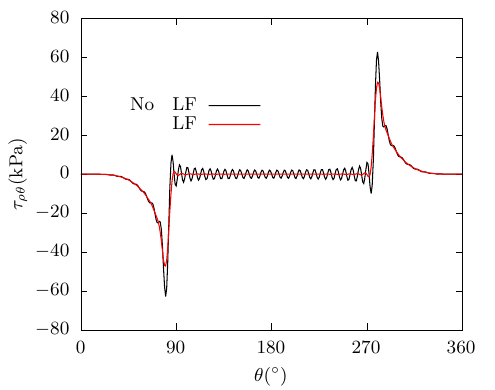}
    \\
    (c) Radial stress along ground surface
    &
      (d) Shear stress along ground surface
    \\
    \includegraphics{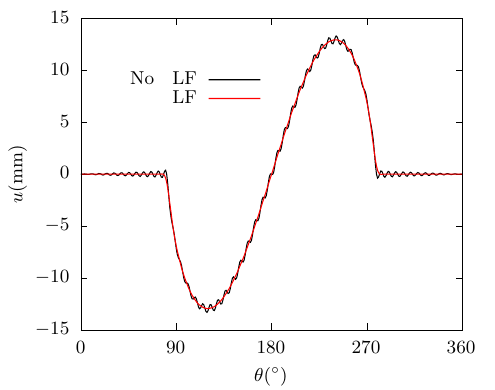}
    &
      \includegraphics{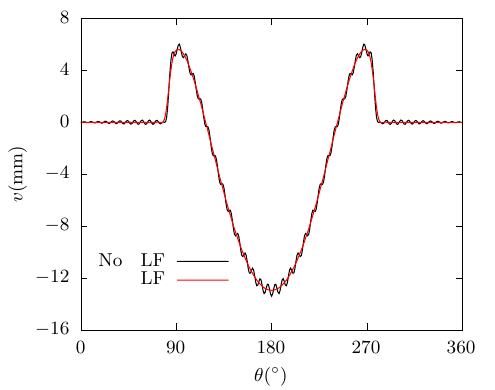}
    \\
    (e) Horizontal displacement along ground surface
    &
      (f) Vertical displacement along ground surface
    \\
  \end{tabular}
  \caption{The Lanczos filtering (LF) for stress and displacement components along ground surface in the mapping plane}
  \label{fig:4}
\end{figure}

We first illustrate the neccessity of the Lanczos filtering in Eq.~(\ref{eq:5.4}). The partially free ground surface length is selected as $ x_{0}/h = 1 $ for better demonstration, since a large $ x_{0} $ would make the displacement constraint arc in the mapping plane too small, referring to the equation $\theta_{0} = -{\rm i}\ln \frac{x_{0}+{\rm i}a}{x_{0}-{\rm i}a}$ mentioned in Section~\ref{sec:problem-3} and arc $ \wideparen{t_{1}t_{2}} $ Fig.~\ref{fig:3}b. Substituting the parameters in Table~\ref{tab:1} and $ x_{0}/h = 1 $ into the solution in this study, we obtain the stress and displacement components along ground surface in the mapping plane when using the Lanczos filtering in Eq.~(\ref{eq:5.4}) or not, as shown in Fig.~\ref{fig:4}, where the abbreviation {\emph{LF}} in the figure denotes Lanczos filtering.

As mentioned in Section~\ref{sec:stress-displ-solut}, stress and displacement would show obvious oscillations of Gibbs phenomena due to abrupt change of boundary condition along ground surface. The results in Fig.~\ref{fig:4} may serve as a strong evidence of the Gibbs phenomena along the ground surface, and also indicates that the Lanczos filtering is neccessary to reduce the Gibbs phenomena. It should be addressed that the Gibbs phenomena {\emph{can not be fully eliminated}}, but can only be greatly reduced by the Lanczos filtering in Eq.~(\ref{eq:5.4}).

\subsection{Solution convergence}
\label{sec:solution-convergence}

Now we illustrate the convergence of the solution proposed in this study. In tunnel engineering, the ground surface is generally assumed free from any traction~\cite[]{Self2020JEM,Lu2019new_solution,Zengguisen2019,Lu2016,Verruijt_Strack2008buoyancy,Verruijt1997traction,Verruijt1997displacement}. In the solution in this study, the ground surface between points $ T_{1} $ and $ T_{2} $ is free from traction, which is identical to the assumptions in the studies mentioned above. What is different from these studies is that the ground surface outside of points $ T_{1} $ and $ T_{2} $ is displacement-constrained, which is different from the assumptions~\cite[]{Self2020JEM,Lu2019new_solution,Zengguisen2019,Lu2016,Verruijt_Strack2008buoyancy,Verruijt1997traction,Verruijt1997displacement}.

\begin{figure}[htb]
  \centering
  \begin{tabular}{cc}
    \begin{tikzpicture}
      \fill [gray!30] (0,0) rectangle (7,-6);
      \draw [cyan, line width = 1pt] (0,0) -- (7,0);
      \draw [line width = 1.5pt] (6,0) -- (7,0);
      \draw [<-] (3.5,0) -- (3.5,-0.5) node [below, align=left] {
        Stress and displacement \\
        along ground surface
      };
      \fill [pattern = north east lines] (6,0) rectangle (7,0.2);
      \fill [white] (0,-1) arc [start angle = 90, end angle = -90, radius = 1] -- (0,-1);
      \draw [dash dot] (0,0.5) -- (0,-6);
      \draw [violet, line width = 1pt] (0,-1) arc [start angle = 90, end angle = -90, radius = 1];
      \fill (0,-2) circle [radius = 0.05];
      \draw [->] (0,-2) -- (1,-2) node [right] {$ R=5{\rm m} $};
      \draw (0,0) -- (-0.5,0);
      \draw (0,-2) -- (-0.5,-2);
      \draw [<->] (-0.2,0) -- (-0.2,-2);
      \node at (-0.2,-1) [rotate=90, above] {$ h = 2R $};
      \draw (6,0) -- (6,0.5);
      \draw [<->] (0,0.2) -- (6,0.2);
      \node at (3,0.2) [above] {$ x_{0} $};
      \node at (0,-4) [below right] {$ \gamma = 20{\rm kN/m^{3}}, k_{0} = 0.8, E = 20{\rm MPa} $};
      \node at (0,-5) [below right] {$ N =50 $};
    \end{tikzpicture}
    &
      \includegraphics{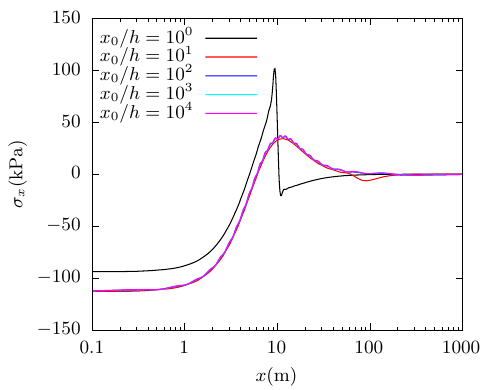}
    \\
    (a) Tunnel and geomaterial model in the physical plane
    &
      (b) Horizontal stress along ground surface
    \\
    \includegraphics{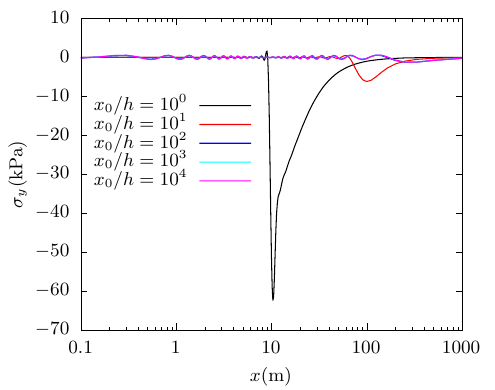}
    &
      \includegraphics{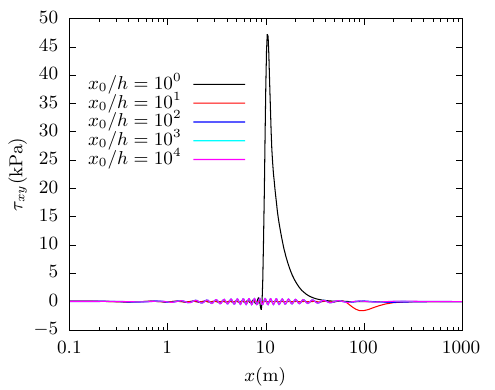}
    \\
    (c) Vertical stress along ground surface
    &
      (d) Shear stress along ground surface
    \\
    \includegraphics{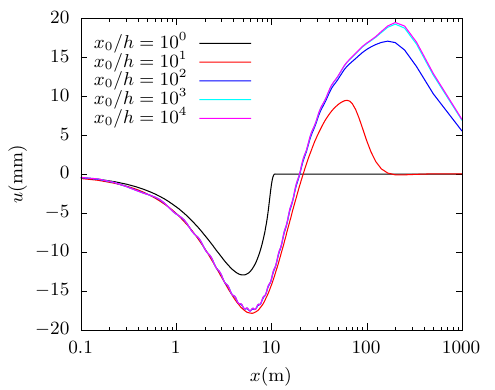}
    &
      \includegraphics{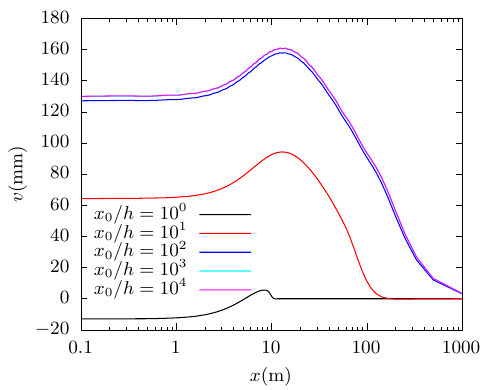}
    \\
    (e) Horizontal displacement along ground surface
    &
      (f) Vertical displacement along ground surface
    \\
  \end{tabular}
  \caption{Convergence of stress and displacement components along ground surface}
  \label{fig:5}
\end{figure}

\begin{figure}[htb]
  \centering
  \begin{tabular}{cc}
    \begin{tikzpicture}
      \fill [gray!30] (0,0) rectangle (7,-6);
      \draw [cyan, line width = 1pt] (0,0) -- (7,0);
      \draw [line width = 1.5pt] (6,0) -- (7,0);
      \fill [pattern = north east lines] (6,0) rectangle (7,0.2);
      \fill [white] (0,-1) arc [start angle = 90, end angle = -90, radius = 1] -- (0,-1);
      \draw [dash dot] (0,0.5) -- (0,-6);
      \draw [violet, line width = 1pt] (0,-1) arc [start angle = 90, end angle = -90, radius = 1];
      \draw [<-] ({cos(45)},{-2+sin(45)}) -- ({1.5*cos(45)},{-2+1.5*sin(45)}) -- ({0.5+1.5*cos(45)},{-2+1.5*sin(45)}) node [right, align = left] {
        Stress and displacement\\
        along tunnel periphery
      };
      \fill (0,-2) circle [radius = 0.05];
      \draw [->] (0,-2) -- (1,-2) node [right] {$ R=5{\rm m} $};
      \draw (0,0) -- (-0.5,0);
      \draw (0,-2) -- (-0.5,-2);
      \draw [<->] (-0.2,0) -- (-0.2,-2);
      \node at (-0.2,-1) [rotate=90, above] {$ h = 2R $};
      \draw (6,0) -- (6,0.5);
      \draw [<->] (0,0.2) -- (6,0.2);
      \node at (3,0.2) [above] {$ x_{0} $};
      \node at (0,-4) [below right] {$ \gamma = 20{\rm kN/m^{3}}, k_{0} = 0.8, E = 20{\rm MPa} $};
      \node at (0,-5) [below right] {$ N =50 $};
    \end{tikzpicture}
    &
      \includegraphics{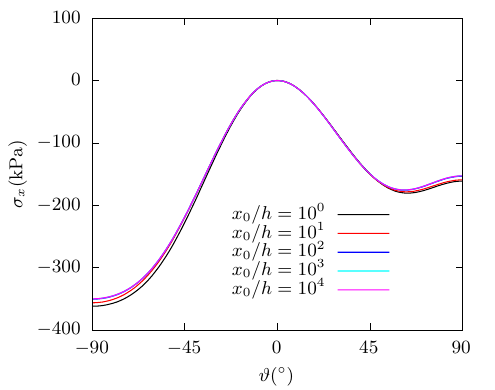}
    \\
    (a) Tunnel and geomaterial model in the physical plane
    &
      (b) Horizontal stress along tunnel periphery
    \\
    \includegraphics{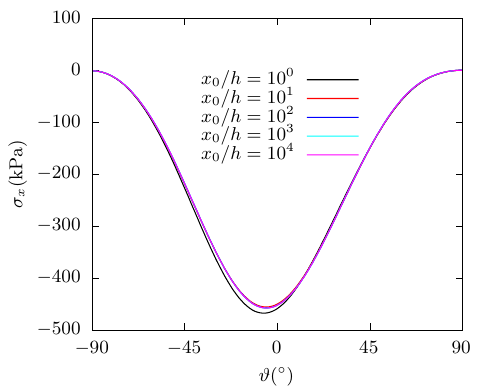}
    &
      \includegraphics{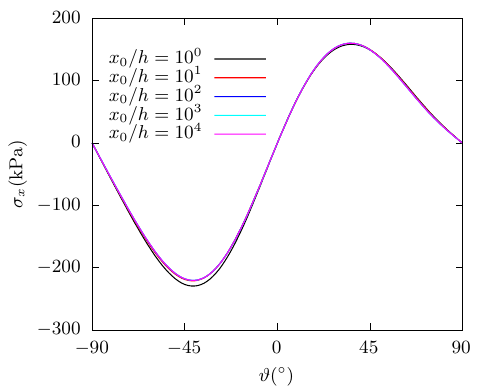}
    \\
    (c) Vertical stress along tunnel periphery
    &
      (d) Shear stress along tunnel periphery
    \\
    \includegraphics{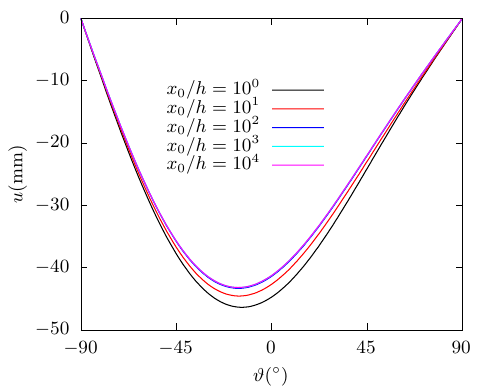}
    &
      \includegraphics{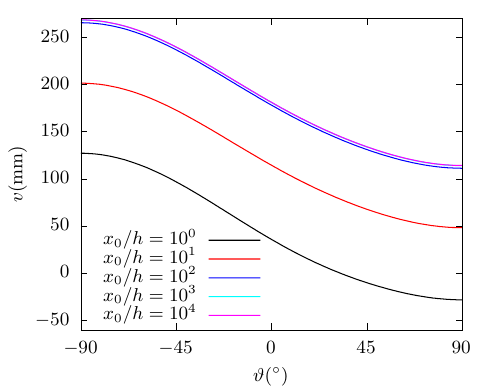}
    \\
    (e) Horizontal displacement along tunnel periphery
    &
      (f) Vertical displacement along tunnel periphery
    \\
  \end{tabular}
  \caption{Convergence of stress and displacement along tunnel periphery}
  \label{fig:6}
\end{figure}

Due to boundary effect, we should expect that as the horizontal coordinates of points $ T_{1} $ and $ T_{2} $ get larger, in other words, the absolute value of $ x_{0} $ gets larger, we should observe convergence of all stress and displacement components. To validate such an expectation, we should substitute the parameters in Table~\ref{tab:1} and different values of $ x_{0}/h $ into the proposed solution to observe the trends of stress and displacement components. To be representative, we select the values of $ x_{0}/h = 10^{p_{1}} $, where $ p_{1} = 0, 1, 2, 3, 4 $. Owing to maximum modulus theorem of complex variable, only the stress and displacement components along both boundaries (the ground surface and tunnel periphery) need to be considered. After substitution and computation, the results are shown in Figs.~\ref{fig:5} and~\ref{fig:6}, respectively. Due to axisymmetry, only the right half geomaterial and corresponding stress and displacement components are illustrated.

Figs.~\ref{fig:5} and~\ref{fig:6} indicate that when the value of $ x_{0}/h $ gets larger, all the stress and displacement components trend to certain curves. The results strongly suggest convergence of the solution. We should also notice that as the value of $ x_{0}/h $ gets larger, the Gibbs phenomena are first reduced, and are then magnified, indicating that there is an optimum value of $ x_{0}/h $. The reason is that when $ x_{0}/h $ is very large, the angle of the displacement constraint arc $ \theta_{0} $ in the mapping plane would be very small correspondingly, and the numerical computation would be erratic, as can be seen in Eq.~(\ref{eq:4.3}). Though conceptually a larger value of $ x_{0}/h $ is more mechanically reasonable, the numerical computation accuracy should also be considered for tradeoff. When $ x_{0}/h = 10^{2} $, all the stress and displacement along the ground surface or tunnel periphery are convergent enough, and the numerical results are also satisfactory. Therefore, the value of $ x_{0}/h = 10^{2} $ is recommanded for actual computation.

\subsection{Verification with existing solution}
\label{sec:comp-with-exist}

The solution convergence does not neccessarily indicates the correctness of the proposed solution. Therefore, a solution comparison should be conducted for verification. We choose the existing modified solution by Lin et al.~\cite[]{Self2020JEM} for comparisons. The reasons of choosing this solution are listed below: (1) Both solutions aim at cancelling the displacement singularity at infinity, but take different measures, which makes the comparisons more valuable than simply verifying the correctness of the proposed solution; to be specific, the proposed solution in this paper takes a mechanical measure of a mixed boundary value problem, while the other solution takes a mathematical measure to simply cancelling the effect of logarithmic item. (2) The assumptions of the proposed solution in this study and the existing modified solution by Lin et al.~\cite[]{Self2020JEM} are similar, and the former would degenerate to the latter, when cancelling the displacement constraint along the ground surface. For consistency of accuracy, the progamming codes of the modified solution by Lin et al.~\cite[]{Self2020JEM} are also rewritten in {\tt{fortran}}.

\begin{figure}[htb]
  \centering
  \begin{tabular}{cc}
    \includegraphics{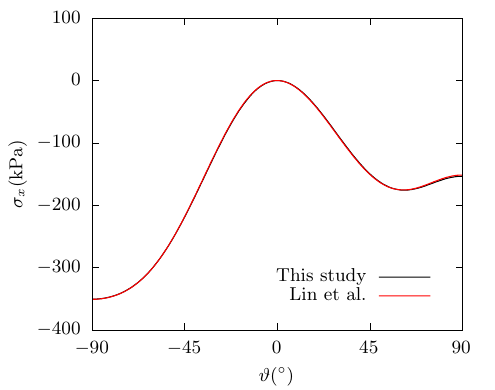}
    &
      \includegraphics{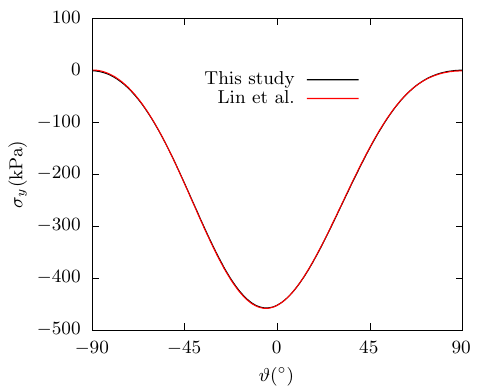}
    \\
    (a) Horizontal stress along tunnel periphery
    &
      (b) Vertical stress along tunnel periphery
    \\
    \includegraphics{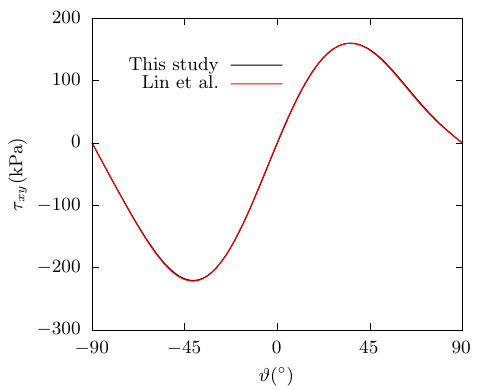}
    &
      \includegraphics{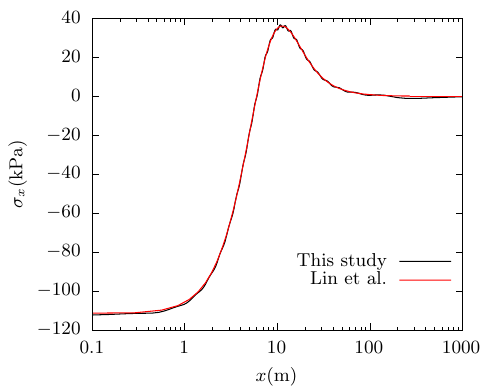}
    \\
    (c) Shear stress along tunnel periphery
    &
      (d) Horizontal stress along ground surface
    \\
    \includegraphics{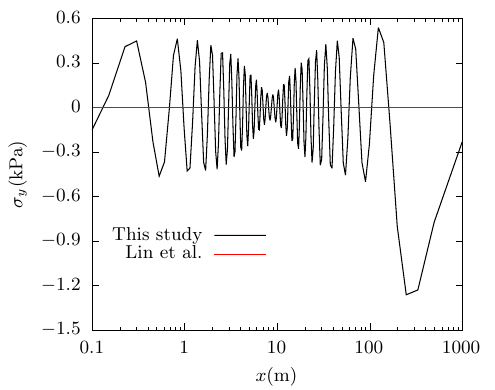}
    &
      \includegraphics{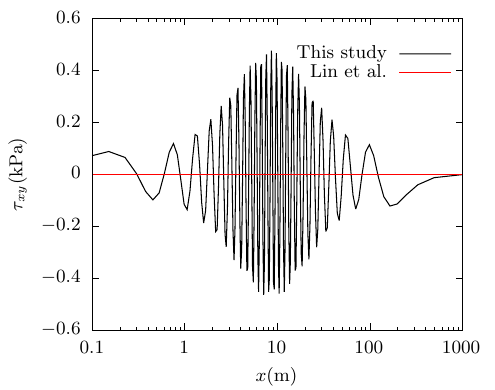}
    \\
    (e) Vertical stress along ground surface
    &
      (f) Shear stress along ground surface
  \end{tabular}
  \caption{Stress comparisons between present study and Lin et al.~\cite[]{Self2020JEM}}
  \label{fig:7}
\end{figure}

Substituting the parameters in Table~\ref{tab:1} into both solutions, and $ x_{0}/h = 10^{2} $ into the proposes solution as well, we obtain the comparisons of stress components along ground surface and tunnel periphery in Fig.~\ref{fig:7}. Fig.~\ref{fig:7} indicates that the stress components for both solutions along tunnel periphery are in agreements. For the ground surface, the stress components between these two solutions slighty deviate. Although the deviation of Figs.~\ref{fig:7}e and~\ref{fig:7}f seems obvious and acute, the absolute value of deviation is small, especially comparing to the horizontal stress in Fig.~\ref{fig:7}d. The reason of the deviation is the remaining effect of the Gibbs phenomena, which {\emph{can not be fully eliminated}} after reduction, as mentioned at the end of Section~\ref{sec:lanczos-filtering}. Fig.~\ref{fig:7} indicates that the stress components computed via the proposed solution are identical to those computed via the existing modified solution by Lin et al.~\cite[]{Self2020JEM}. The identical stress results further indicate that the solution of $ d_{n} $ in Eq.~(\ref{eq:4.1}) is correct.

\begin{figure}[htb]
  \centering
  \begin{tabular}{cc}
    \includegraphics{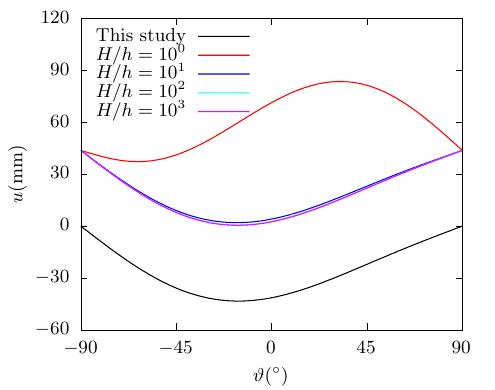}
    &
      \includegraphics{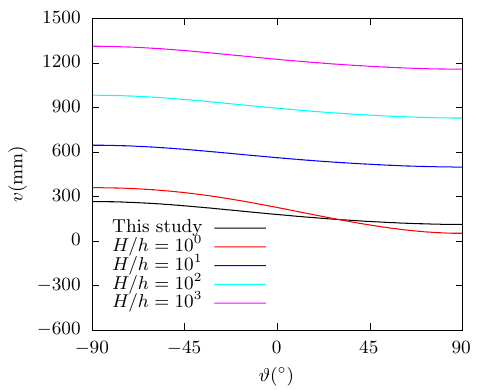}
    \\
    (a) Horizontal displacement along tunnel periphery
    &
      (b) Vertical displacement along tunnel periphery
    \\
    \includegraphics{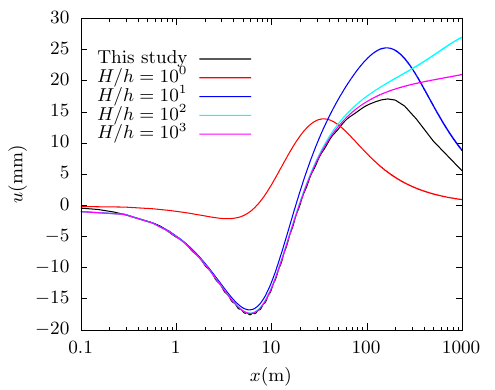}
    &
      \includegraphics{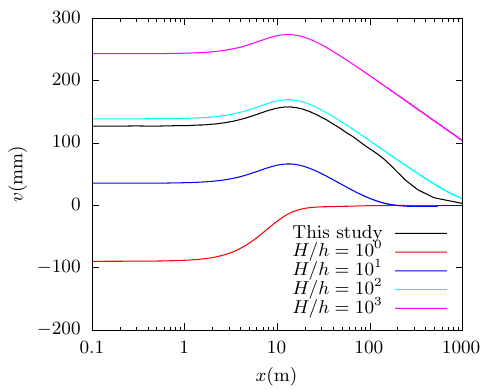}
    \\
    (c) Horizontal displacement along ground surface
    &
      (d) Vertical displacement along ground surface
    \\
  \end{tabular}
  \caption{Displacement comparisons between present study and Lin et al.~\cite[]{Self2020JEM}}
  \label{fig:8}
\end{figure}

The same parameter set obtaining stress is again substituted into the proposed solution to obtain the displacement when $ x_{0}/h = 10^{2} $. Meanwhile, it should be noted that the displacement of the modified solution by Lin et al.~\cite[]{Self2020JEM} further depends on the modified depth $H$, which is the key parameter to cancelling the displacement singularity at infinity. Such modification is only a purely mathematical strategy and lacks a strong mechanical foundation. We choose the value of $ H/h = 10^{p_{2}} $, where $ p_{2} = 0,1,2,3 $ to compare with the displacement of the proposed solution in this study. The comparing results are shown in Fig.~\ref{fig:8}.

Figs.~\ref{fig:8}a and~\ref{fig:8}b indicate that when $ H/h $ is larger than $ 10^{1} $, the curve trends of displacement components along tunnel periphery of both solutions are almost the same, and the displacement difference between these two solutions along tunnel periphery is almost a changable constant for different value of $ H/h $. Owing to axisymmetry, the horizontal displacement of angles $ \vartheta = \pm 90^{\circ} $ in Fig.~\ref{fig:8}a is expected be zero, which is identical to the result of the proposed solution in this study, while the horizontal displacement of angles $ \vartheta = \pm 90^{\circ} $ for all four curves of the modified solution by Lin et al.~\cite[]{Self2020JEM} are nonzero. Moreover, the parameter $ x_{0}/h = 10^{2} $ is chosen for the proposed solution in this study, because the displacement components are close to the convergent ones, no matter for the horizontal or the vertical one, as shown in last section. Fig.~\ref{fig:8}b suggests that the absolute value of the vertical displacement along tunnel periphery of the modified solution by Lin et al.~\cite[]{Self2020JEM} would always get larger for a larger value of $ H/h $.

Figs.~\ref{fig:8}c and~\ref{fig:8}d provide more information about the displacement difference between the proposed solution in this study and the one by Lin et al.~\cite[]{Self2020JEM}. Fig.~\ref{fig:8}c indicates that as the value of $ H/h $ gets larger, the horizontal displacement of the modified solution along the ground surface near the tunnel ($ x \leq 100{\rm m} $) by Lin et al.~\cite[]{Self2020JEM} is more approximate to that by the proposed solution, while the far-field horizontal displacements ($ x > 100{\rm m} $) between these two solution are more deviated. Fig.~\ref{fig:8}d indicates that as the value of $ H/h $ gets larger, the vertical displacement of the modified solution along the ground surface near the tunnel ($ x \leq 100{\rm m} $) by Lin et al.~\cite[]{Self2020JEM} is more deviated from the one of the proposed solution in this study.

The results in Section~\ref{sec:comp-with-exist} indicate that the modified solution by Lin et al.~\cite[]{Self2020JEM} is not fully convergent, though this solution eliminates the displacement singularity at infinity in both theoretical and computational aspects. In contrast, the proposed solution has a better convergence than the one by Lin et al.~\cite[]{Self2020JEM}. The reason of convergence comes from the mechanical assumption that the far-field ground surface is displacement-constrained ($ |x| \geq x_{0} $), and the constraint provides a constraining force to equilibrate the unbalanced resultant applied along tunnel periphery due to shallow tunnel excavation. In other words, the displacement constraint assumption turn the originally unbalanced resultant problem in Refs~\cite[]{Self2020JEM,Lu2019new_solution,Zengguisen2019,Lu2016,Verruijt_Strack2008buoyancy,Verruijt1997traction,Verruijt1997displacement} into a balanced one with mixed boundaries in this study.

For simpler and deeper understanding, the problem and solution in this study is mechanically and conceptually similar to a one-dimensional cantilever beam subjected to a fixed constraint at one end and an axisymmetrically and longitudinally concentrated force at the other end. The difference is that the dimension of the problem and corresponding solution measure in this paper is elevated from one dimension to two dimensions, and the region in the problem also changes from a simply connected one to a doubly connected one, but the conceptual solution measure remains the same. To be more specific, the static equilibrium in Eq.~(\ref{eq:4.6}) and the single-valuedness of displacement in Eq.~(\ref{eq:4.12}) consist of the mechanical foundation to construct the whole solution and determine the unbalanced resultant along tunnel periphery, while the traction boundary condition in Section~\ref{sec:solution-riemann-hilbert-3} establishes the remaining neccessary mathematical equations to reach the unique solution of $ d_{n} $ in Eq.~(\ref{eq:4.1}).

\section{Parametric investigation and discussion}
\label{sec:param-invest}

In this section, a parametric investigation is performed regarding the effects of three distinct parameters in the proposed mechanical model on stress and displacement distribution in geomaterial. The parameters of benchmark model are the same to those in Table~\ref{tab:1} (including plain strain condition $ \kappa = 3-4\nu $), except that $ k_{0} = 1 $ for hydrostatic condition. The free surface range is selected as $ x_{0}/h = 10^{2} $ in the following Sections~\ref{sec:tunnel-depth} and~\ref{sec:lateral-coefficient} due to convergence analysis in Section~\ref{sec:solution-convergence}, and will be altered in Section~\ref{sec:free-surface-range} for an extending discussion of the convergence analysis in Section~\ref{sec:solution-convergence}. To conduct non-dimensional parametric investigation, the stress and displacement in geomaterial are normalized by $ \gamma h $ (far-field hydrostatic stress at tunnel center depth) and $ u_{0} = \frac{\gamma h R}{2G} $ (radial displacement along tunnel periphery subjected to far-field hydrostatic stress $ \gamma h $), respectively. To better evaluate the stress distribution, the pair of non-dimensional maximun and minimum principle stresses ($ \sigma_{\max}/\gamma h $ and $ \sigma_{\min}/\gamma h $) is used, instead of the three stress components in Eq.~(\ref{eq:5.4}). 

\subsection{Tunnel depth}
\label{sec:tunnel-depth}

\begin{figure}[htb]
  \centering
  \begin{tabular}{cc}
    \includegraphics{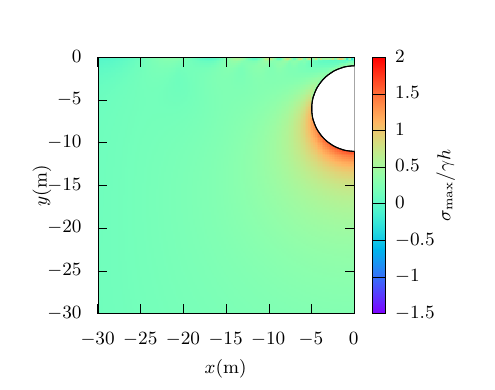}
    &
      \includegraphics{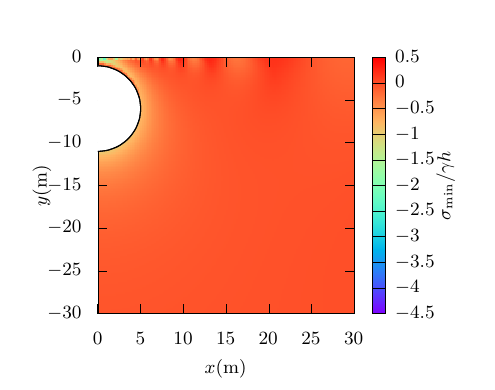}
    \\
    (a) Maximum principle stress, $ h/R=1.1 $
    &
      (b) Minimum principle stress, $h/R=1.1$
    \\
    \includegraphics{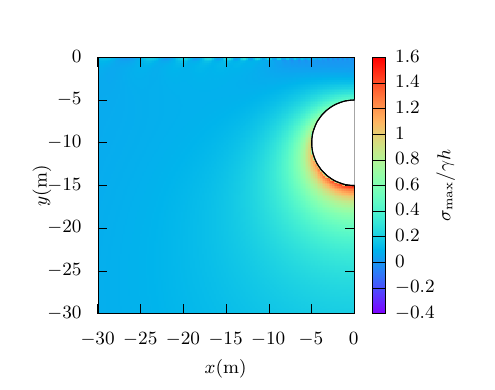}
    &
      \includegraphics{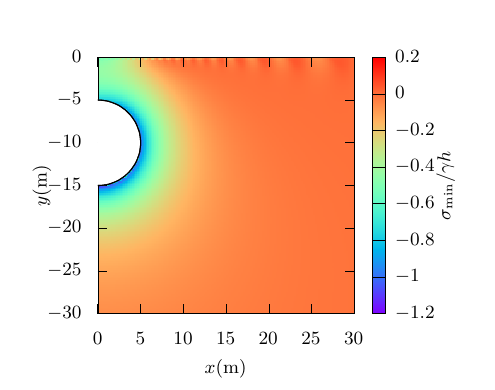}
    \\
    (c) Maximum principle stress, $ h/R=2 $
    &
      (d) Minimum principle stress, $h/R=2$
    \\
    \includegraphics{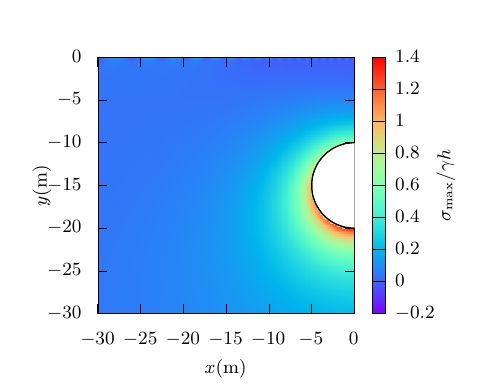}
    &
      \includegraphics{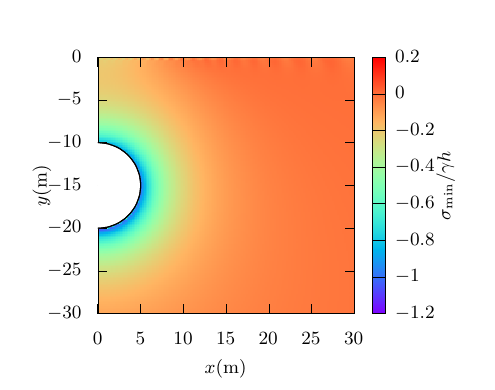}
    \\
    (e) Maximum principle stress, $ h/R=3 $
    &
      (f) Minimum principle stress, $h/R=3$
    \\
  \end{tabular}
  \caption{Normalized principle stress distribution against tunnel depth}
  \label{fig:9}
\end{figure}

\begin{figure}[htb]
  \centering
  \begin{tabular}{cc}
    \includegraphics{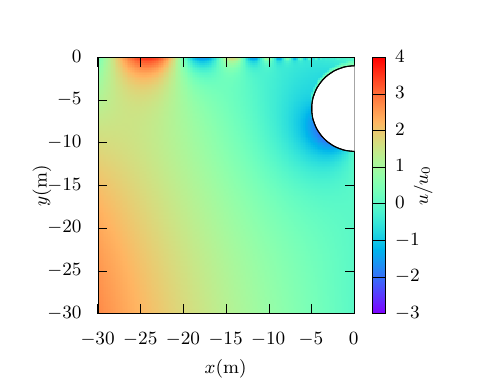}
    &
      \includegraphics{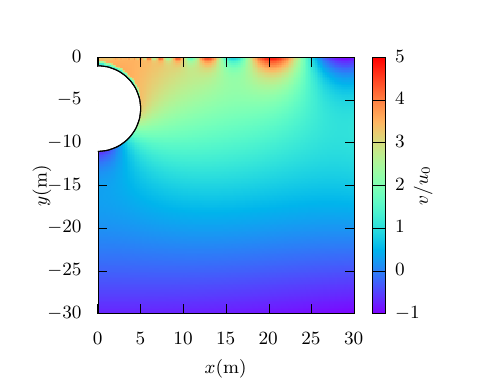}
    \\
    (a) Horizontal displacement, $ h/R=1.1 $
    &
      (b) Vertical displacement, $h/R=1.1$
    \\
    \includegraphics{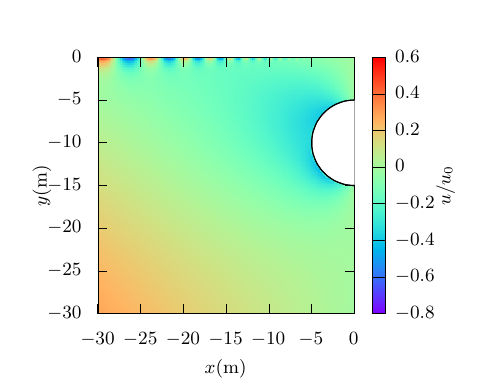}
    &
      \includegraphics{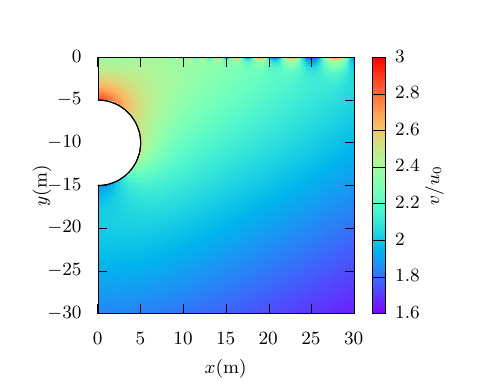}
    \\
    (c) Horizontal displacement, $ h/R=2 $
    &
      (d) Vertical displacement, $h/R=2$
    \\
    \includegraphics{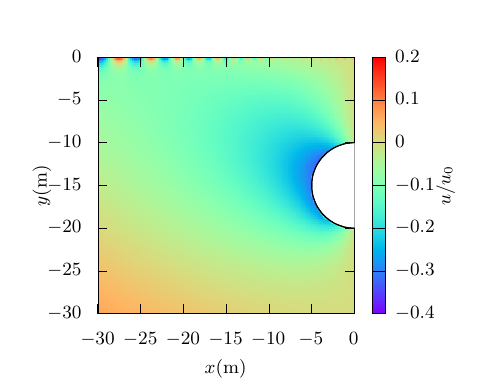}
    &
      \includegraphics{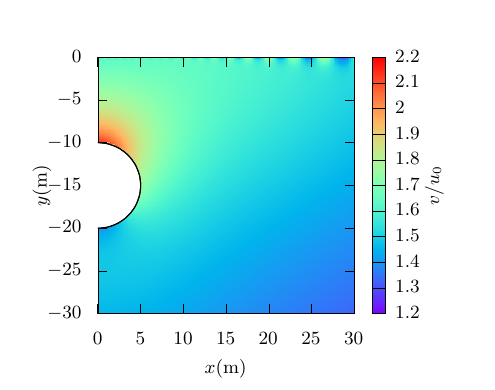}
    \\
    (e) Horizontal displacement, $ h/R=3 $
    &
      (f) Vertical displacement, $h/R=3$
    \\
  \end{tabular}
  \caption{Normalized displacement distribution against tunnel depth}
  \label{fig:10}
\end{figure}

Tunnel depth would greatly affect stress and displacement distribution in geomaterial in the proposed mechanical model. The non-dimensional tunnel depths are chosen as $ h/R = 1.1, 2, 3 $. Substituting all the parameters into the proposed mechanical model, the normalized principle stress and displacement distributions are obtained in Figs.~\ref{fig:9} and~\ref{fig:10}, respectively. Owing to axisymmetry of the proposed mechanical model, only the left or right half geomaterial in physical plane $xOy$ is presented. Such a presentation feature will be used in the following analyses.

The principle stress distributions in Fig.~\ref{fig:9} suggest that the maximum and minimum principle stresses generally concentrate at the bottom geomaterial ($ \vartheta \in [180^{\circ}, 360^{\circ}]$) along tunnel periphery and along the whole of tunnel periphery, respectively, where the local radial coordinate system $\varrho o \vartheta$ in Fig.~\ref{fig:1}b is used. Therefore, the stress at bottom of tunnel periphery may be monitored during tunnel excavation, and maybe neccessary engineering measures may be conducted to prevent possible hazards. Specifically, Fig.~\ref{fig:9}b suggests a severe tensile stress concentration above tunnel roof, when tunnel depth is very small. Thus, a extremely shallow tunnel should be avoided in design.

Figs.~\ref{fig:10}a, c, and e suggest that geomaterial near tunnel side walls ($ \vartheta \in [-45^{\circ}, 45^{\circ}] \cup [135^{\circ}, 225^{\circ}] $) horizontally deforms away tunnel cavity. When tunnel depth is small the relative horizontal displacement is larger. Figs.~\ref{fig:10}b, d, and f suggest that geomaterial generally deforms in a upward manner, which is identical to the upward resultant of shallow tunnel excavation, and the upward deformation of geomaterial near tunnel roof ($ \vartheta \in [60^{\circ}, 120^{\circ}] $) is the most severe. Displacement oscillations are observed along ground surface due to the reduction of Gibbs phenomena. The deformation of tunnel roof should be monitored during excavation.

\subsection{Lateral coefficient}
\label{sec:lateral-coefficient}

\begin{figure}[htb]
  \centering
  \begin{tabular}{cc}
    \includegraphics{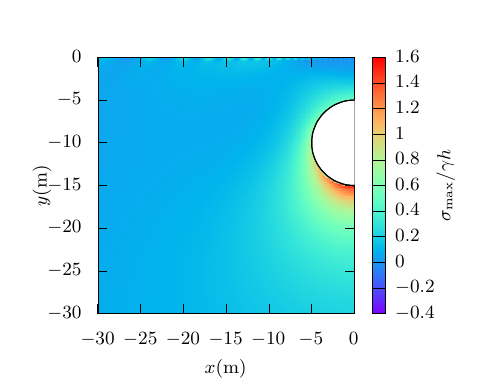}
    &
      \includegraphics{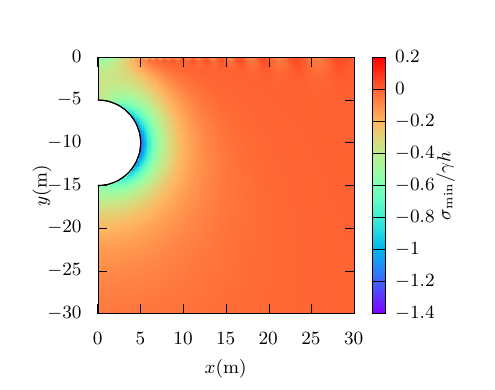}
    \\
    (a) Maximum principle stress, $ k_{0}=0.8 $
    &
      (b) Minimum principle stress, $ k_{0}=0.8 $
    \\
    \includegraphics{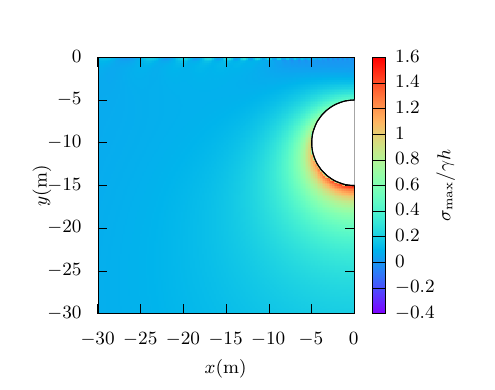}
    &
      \includegraphics{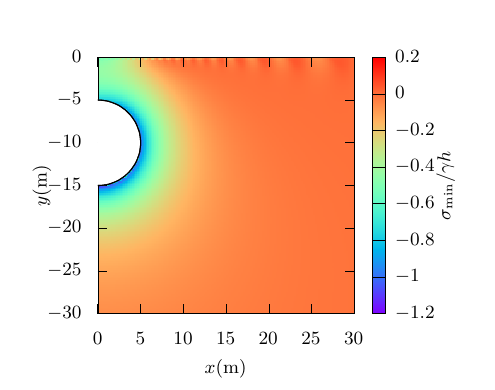}
    \\
    (c) Maximum principle stress, $ k_{0}=1 $
    &
      (d) Minimum principle stress, $ k_{0}=1 $
    \\
    \includegraphics{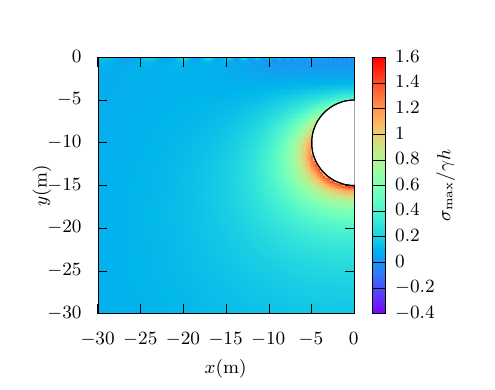}
    &
      \includegraphics{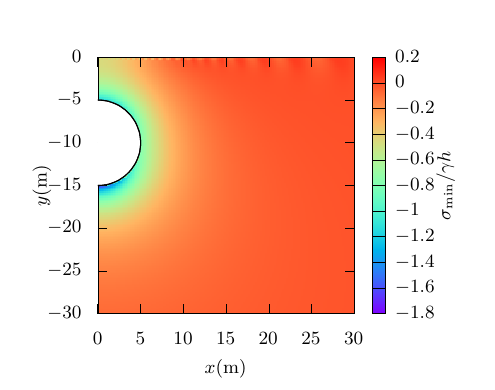}
    \\
    (e) Maximum principle stress, $ k_{0}=1.2 $
    &
      (f) Minimum principle stress, $ k_{0}=1.2 $
    \\
  \end{tabular}
  \caption{Normalized principle stress distribution against lateral coefficient}
  \label{fig:11}
\end{figure}

\begin{figure}[htb]
  \centering
  \begin{tabular}{cc}
    \includegraphics{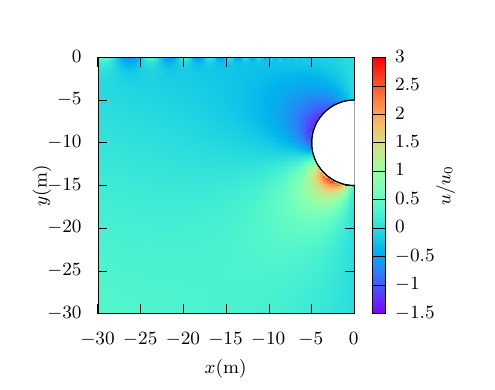}
    &
      \includegraphics{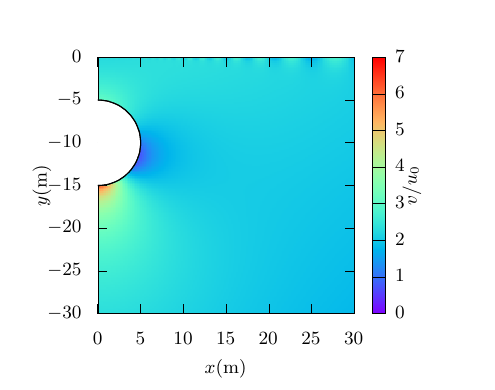}
    \\
    (a) Horizontal displacement, $ k_{0}=0.8 $
    &
      (b) Vertical displacement, $ k_{0}=0.8 $
    \\
    \includegraphics{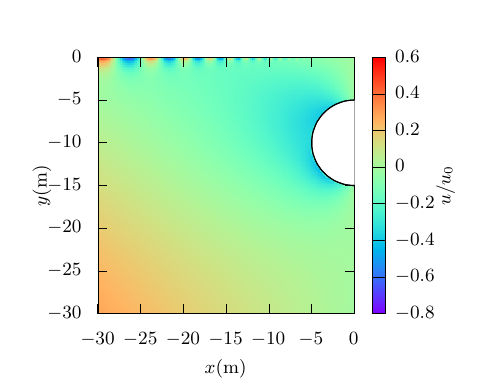}
    &
      \includegraphics{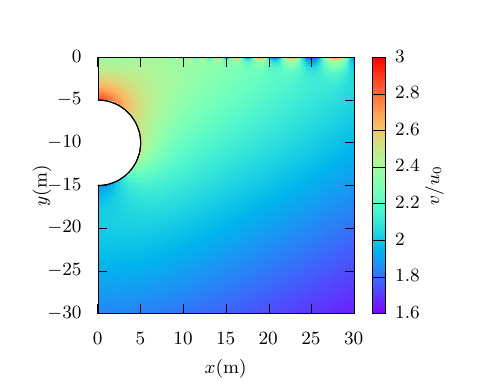}
    \\
    (c) Horizontal displacement, $ k_{0}=1 $
    &
      (d) Vertical displacement, $ k_{0}=1 $
    \\
    \includegraphics{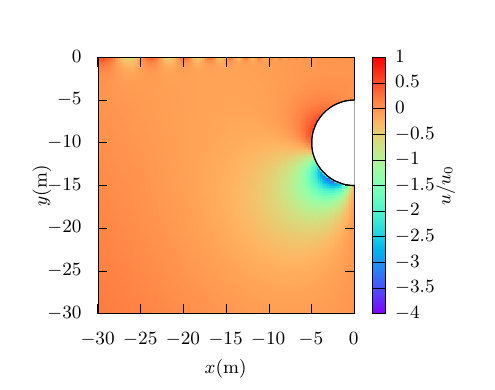}
    &
      \includegraphics{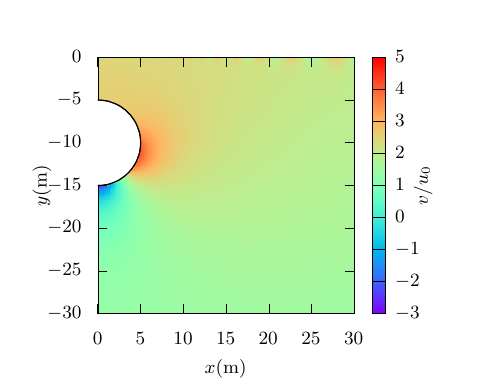}
    \\
    (e) Horizontal displacement, $ k_{0}=1.2 $
    &
      (f) Vertical displacement, $ k_{0}=1.2 $
    \\
  \end{tabular}
  \caption{Normalized displacement distribution against lateral coefficient}
  \label{fig:12}
\end{figure}

The in-situ lateral coefficient $ k_{0} $ varies in a certain range, and we select $ k_{0} = 0.8, 1, 1.2 $, which typically present vertical-stress-dominating, hydrostatic, and horizontal-stress-dominating in-situ conditions, respectively. With different $ k_{0} $ and other parameters, the normalized principle stress and displacement distributions are obtained in Figs.~\ref{fig:11} and~\ref{fig:12}, respectively.

Figs.~\ref{fig:11}a, c, and e suggest an enlarging area of maximum principle concentration at the bottom of tunnel periphery for an increasing lateral coefficient, and the magnitude of normalized maximum principle stress remains almost the same. Figs.~\ref{fig:11}b, d, and f, however, suggest a reverse distribution feature of minimum principle stress, when lateral coefficient increases from 0.8 to 1.2. The reason is that $ k_{0} = 0.8 $ and $ k_{0} = 1.2 $ denote the typical in-situ stress conditions of vertical domination and horizontal domination, and tensile areas should be near both tunnel walls and roof-bottom regions, respectively, which is identical to the results in Figs.~\ref{fig:11}b and f. Stress monitoring at side walls and roof-bottom regions is recommended for excavation in geomaterials of low and high lateral in-situ coefficients, respectively.

Fig.~\ref{fig:12}a suggests that when lateral coefficient is small ($ k_{0} = 0.8 $), geomaterial near tunnel walls ($ \vartheta \in [0^{\circ}, 60^{\circ}] \cup [120^{\circ}, 180^{\circ}] $) deforms horizontally away from tunnel cavity, while geomaterial near tunnel bottom ($ \vartheta \in [210^{\circ}, 255^{\circ}] \cup [285^{\circ}, 330^{\circ}] $) deforms horizontally towards tunnel cavity. In contrast, Fig.~\ref{fig:12}e suggests a completely oppisite deformation pattern when lateral coefficient is large ($ k_{0} = 1.2 $). A similar contrast can also be found in Figs.~\ref{fig:12}b and f. The whole tunnel deforms upwards when lateral coefficient is small ($ k_{0} = 0.8 $), indicating a vertical-stress-dominating in-situ conditions. Meanwhile, the obvious downward displacement near tunnel bottom in Fig.~\ref{fig:12}f suggest that the tunnel is horizontally squeezed and deforms to a pseudo-ellipse with vertical major axis, indicating a horizontal-stress-dominating in-situ condition. In other words, the tunnel is vertically tensile for large lateral coefficient, which is identical to the minimum principle stress distribution in Fig.~\ref{fig:11}f. 

\subsection{Free surface range}
\label{sec:free-surface-range}

\begin{figure}[htb]
  \centering
  \begin{tabular}{cc}
    \includegraphics{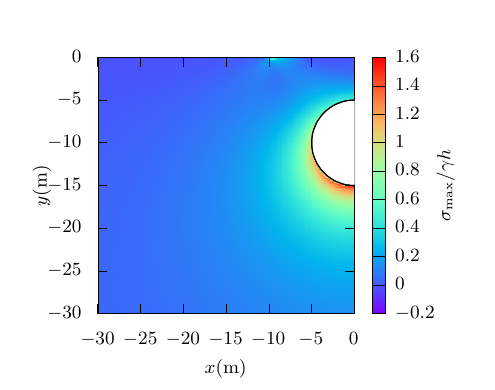}
    &
      \includegraphics{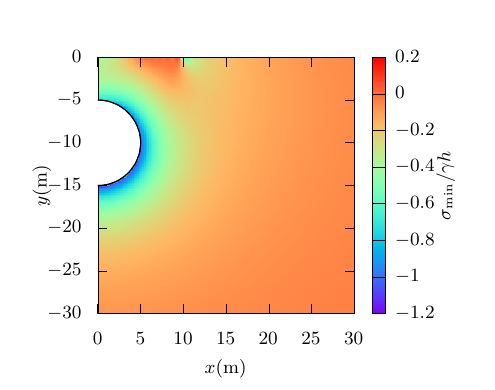}
    \\
    (a) Maximum principle stress, $ x_{0}/h = 10^{0} $
    &
      (b) Minimum principle stress, $x_{0}/h = 10^{0} $
    \\
    \includegraphics{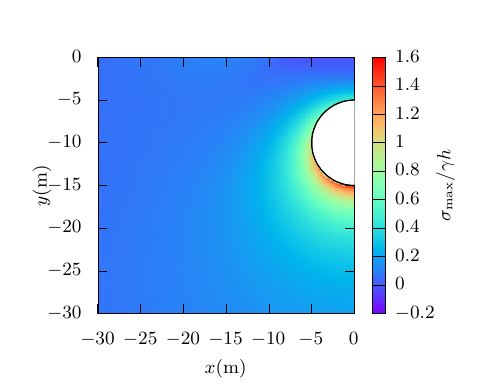}
    &
      \includegraphics{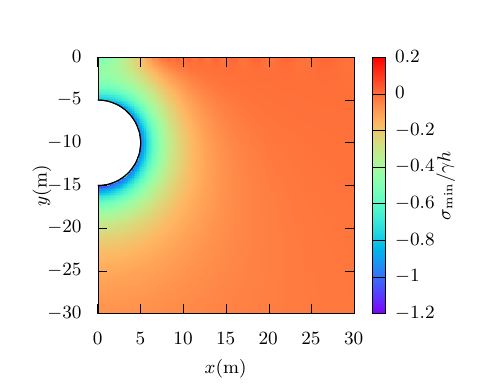}
    \\
    (c) Maximum principle stress, $ x_{0}/h = 10^{1} $
    &
      (d) Minimum principle stress, $x_{0}/h = 10^{1} $
    \\
    \includegraphics{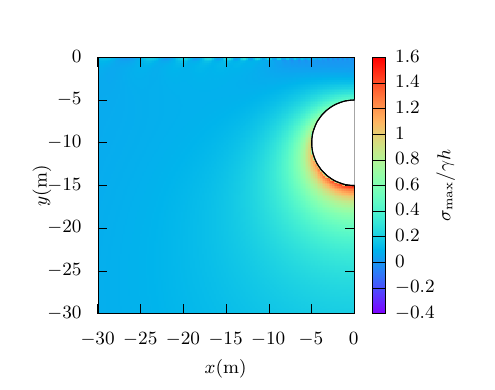}
    &
      \includegraphics{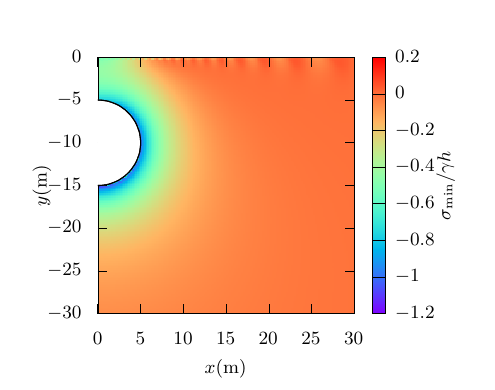}
    \\
    (e) Maximum principle stress, $ x_{0}/h = 10^{2} $
    &
      (f) Minimum principle stress, $x_{0}/h = 10^{2} $
    \\
  \end{tabular}
  \caption{Normalized principle stress distribution against free surface range}
  \label{fig:13}
\end{figure}

\begin{figure}[htb]
  \centering
  \begin{tabular}{cc}
    \includegraphics{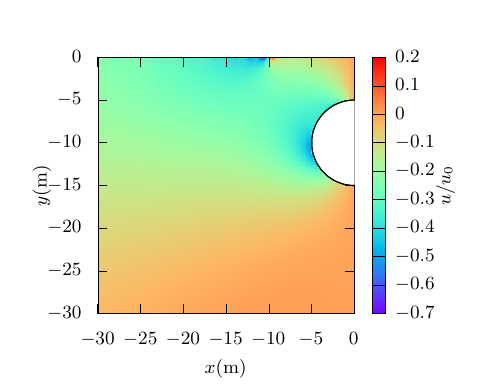}
    &
      \includegraphics{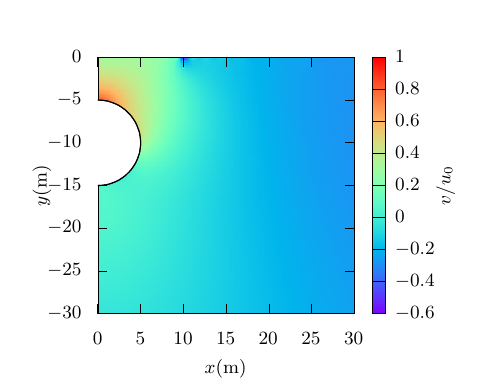}
    \\
    (a) Horizontal displacement, $ x_{0}/h = 10^{0} $
    &
      (b) Vertical displacement, $x_{0}/h = 10^{0} $
    \\
    \includegraphics{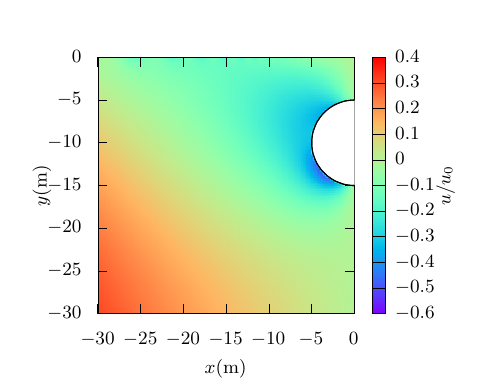}
    &
      \includegraphics{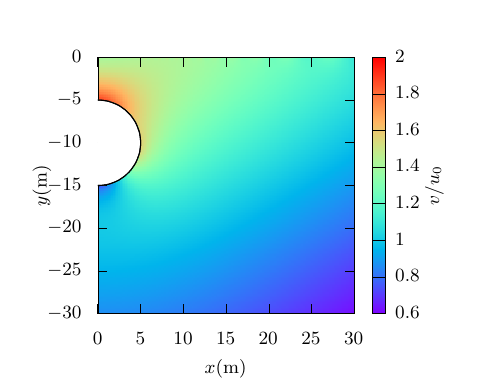}
    \\
    (c) Horizontal displacement, $ x_{0}/h = 10^{1} $
    &
      (d) Vertical displacement, $x_{0}/h = 10^{1} $
    \\
    \includegraphics{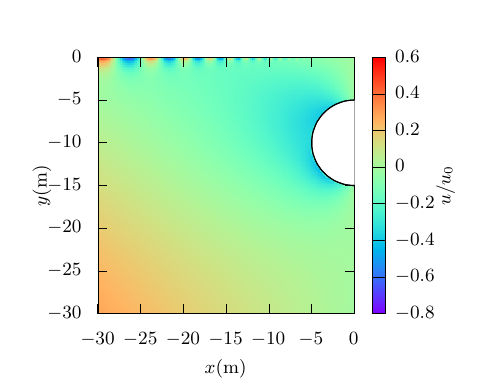}
    &
      \includegraphics{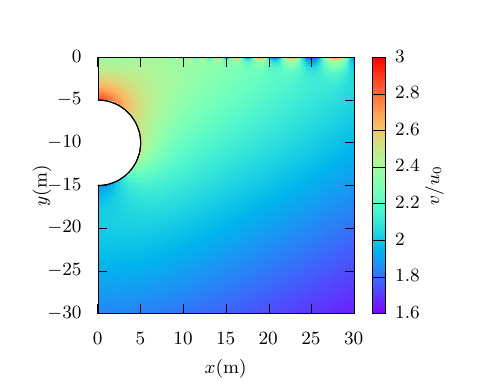}
    \\
    (e) Horizontal displacement, $ x_{0}/h = 10^{2} $
    &
      (f) Vertical displacement, $x_{0}/h = 10^{2} $
    \\
  \end{tabular}
  \caption{Normalized displacement distribution against free surface range}
  \label{fig:14}
\end{figure}

The solution convergence regarding the free surface range $ x_{0}/h $ has been validated in Section~\ref{sec:solution-convergence}, but only along ground surface and tunnel periphery. In this section, validation enhancement is conducted by choosing $ x_{0}/h = 10^{0}, 10^{1}, 10^{2} $, and further discussion is made. Substituting three values of $ x_{0}/h $ and other neccessary parameters into the proposed mechanical model, the normalized principle stress and displacement distributions are obtained in Figs.~\ref{fig:13} and~\ref{fig:14}.

Fig.~\ref{fig:13} shows that the stress distribution in geomterial remains almost the same for different free surface range $ x_{0}/h $, which indicates that the free surface range would not greatly affect stress distribution in geomterial (except for the ground surface). The reason is that the displacement constraint along ground surface is a relatively far-field constraint for the geometry combination of $h/R=2$ in this numerical case, and would spontaneously not greatly alter the stress distribution near tunnel according to Saint Venant Principle.

Fig.~\ref{fig:14} suggests an overall increasing trend of upward displacement in geomaterial for a larger free surface range $ x_{0}/h $, and a similar overall increasing trend of upward displacement in geomaterial is also reported in the analytical solution conducted by Lu et al.~\cite[]{lu2021reasonable} that a Gaussian-distributed load acts along ground surface. The corresponding finite element verification~\cite[]{lu2021reasonable} suggests identical results that a larger geomaterial size would greatly increase the overall vertical displacement, even if the Gaussian distributed load remains the same form. As for this study, the increasing free surface range is similar to the larger geomaterial size of geomaterial in the finite element model in Ref~\cite[]{lu2021reasonable}. Yet it is impossible to build a finite element model of infinitely large lower-half plane, thus, the verification is performed by comparisons with an existing analytical solution in Section~\ref{sec:numer-verifi-discu}, instead of a finite element model.

\section{Conclusion}
\label{sec:conclusion}

This paper proposes a new mechanical model in complex variable method to confront the displacement singularity caused by the unbalanced resultant of shallow tunnel excavation for reasonable stress and displacement in geomaterial. The far-field ground surface displacement is constrained to produce a corresponding constraining force to equilibrate the unbalanced resultant and the original unbalanced problem turns to a static equilibrium one. The constrained far-field ground surface together with the partially free ground surface above the shallow tunnel and the traction distribution along shallow periphery forms mixed boundary value problem, which is converted into a homogenerous Riemann-Hilbert problem by applying the analytic continuation. With extra and implicit boundary conditions of the static equilibrium and the displacement single-valuedness in geomaterial, the homogenerous Riemann-Hilbert problem is solved using an approximate and iterative method, and the reasonable stress and displacement in geomaterial are obtained. The Lanczos filtering is used to reduce the Gibbs phenomena caused by abrupt change of boundary condition along the ground surface. Several numerical cases are conducted to show that the newly proposed mechanical model sucessfully and simultaneously guarantee the the convergence and correctness of stress and displacement in geomaterial and eliminate the displacement singularity at infinity. A parametric investigation is made to further discuss the influence of tunnel depth, lateral coefficient, and free surface range on stress and displacement in geomaterial. The proposed mechanical model might also imply that any unbalanced resultant problem of statics in tunnel engineering is better to be converted into a mixed boundary value problem of static equilibrium to obtain reasonable results for both stress and displacement. 

\clearpage
\section*{Acknowlegement}
\label{sec:acknowlegement}

This study is financially supported by the Natural Science Foundation of Fujian Province, China (Grant No. 2022J05190), the Scientific Research Foundation of Fujian University of Technology (Grant No. GY-Z20094), the National Natural Science Foundation of China (Grant No. 52178318), and the Education Foundation of Fujian Province (Grant No. JAT210287). The authors would like to thank Professor Changjie Zheng, Ph.D. Yiqun Huang, and Associate Professor Xiaoyi Zhang for their suggestions on this study.

\clearpage
\appendix
\section{Coefficient degeneration}
\label{sec:A}

The complex potentials within region ${\bm \omega}$ can be obtained by repective integrations of Eqs.~(\ref{eq:4.4a}) and (\ref{eq:4.5}).
\begin{subequations}
  \label{eqa:1}
  \begin{equation}
    \label{eqa:1a}
    \varphi(\zeta) =  {\rm i}A_{-1}\ln\zeta + a_{0} + \sum\limits_{k=1}^{\infty} a_{k} \zeta^{k} + \sum\limits_{k=1}^{\infty} b_{k} \zeta^{-k}
  \end{equation}
  \begin{equation}
    \label{eqa:1b}
    \psi(\zeta) = - {\rm i}B_{-1}\ln\zeta + c_{0} + \sum\limits_{k=1}^{\infty} c_{k} \zeta^{k} + \sum\limits_{k=1}^{\infty} d_{k} \zeta^{-k}
  \end{equation}
\end{subequations}
where $a_{k}$, $b_{k}$, $c_{k}$, and $d_{k}$ corresponding to the same symbols in Refs~\cite{Verruijt1997traction, Strack_Verruijt2002buoyancy, Lu2016}, and can be expanded as: 
\begin{subequations}
  \label{eqa:2}
  \begin{equation}
    \label{eqa:2a}
    \left\{
      \begin{aligned}
        a_{k} = & \; {\rm i}\frac{A_{k-1}}{k} \\
        b_{k} = & \; {\rm i}\frac{A_{-k-1}}{-k} \\
      \end{aligned}
    \right.
    , \quad k \geq 1
  \end{equation}
  \begin{equation}
    \label{eqa:2b}
    \left\{
      \begin{aligned}
        c_{k} = & \; \frac{1}{2}({\rm i}A_{k-2}-{\rm i}A_{k})-{\rm i}\frac{B_{-k-1}}{k} \\
        d_{k} = & \; \frac{1}{2}({\rm i}A_{-k-2}-{\rm i}A_{-k})-{\rm i}\frac{B_{k-1}}{-k} \\
      \end{aligned}
    \right.
    , \quad k \geq 1
  \end{equation}
\end{subequations}
and $a_{0}$ and $c_{0}$ denotes undetermined complex constants due to integration. If we replace the unbalanced resultant along tunnel periphery with a zero resultant traction (similar to the case in Refs~\cite{Verruijt1997traction}), the unbalanced resultant vanishes, and $A_{-1}=0$ and $B_{-1}=0$ should be correspondingly satisfied to eliminate the possible multi-valuedness of the complex potentials in Eq.~(\ref{eqa:1}). If we further cancel the displacement constraint along the surface near infinity, $A_{k}=B_{k} (k \neq -1)$ in Eq.~(\ref{eq:4.4}) is guaranteed and henceforth.

It can be easily verified that with the notations in Eq.~(\ref{eqa:2}), the equilibriums in Eqs.~(21) and (22) in Ref~\cite{Verruijt1997traction}) would be guaranteed. As for $k=0$, we have
\begin{equation}
  \label{eqa:3}
  c_{0} + \overline{a}_{0} = \frac{{\rm i}}{2}(A_{-2}-A_{0})
\end{equation}
Since both $a_{0}$ and $c_{0}$ are integration constants, Eq.~(\ref{eqa:3}) indicates that $(A_{-2}-A_{0})$ should be a constant as well. Since $A_{0}$ and $A_{-2}$ are both real due to symmetry, Eq.~(\ref {eqa:3}) also indicates that $c_{0} + \overline{a}_{0}$ should a be pure imaginary constant. Such a property will be used later on. Subsequently, the coefficients in Eq.~(\ref{eq:4.14b}) can be correspondingly rewritten as
\begin{subequations}
  \label{eqa:4}
  \begin{equation}
    \label{eqa:4a}
    {\rm e}^{-{\rm i}k\theta}: \quad (r^{-k}-r^{k}){\rm i}\frac{A_{-k-1}}{-k} + (r^{k+2}-r^{-k}){\rm i}\frac{A_{-k-2}}{-k-1} + (1-r^{2})r^{k}({\rm i}A_{k}-A_{k-1}), \quad k \geq 1
  \end{equation}
  \begin{equation}
    \label{eqa:4b}
    {\rm e}^{{\rm i}k\theta}: \quad (r^{k}-r^{-k}){\rm i}\frac{A_{k-1}}{k} + (r^{-k+2}-r^{k}){\rm i}\frac{A_{k-2}}{k-1} + (1-r^{2})r^{-k}({\rm i}A_{-k}-A_{-k-1}), \quad k \geq 2
  \end{equation}
  \begin{equation}
    \label{eqa:4c}
    -r^{-1}{\rm e}^{{\rm i}\theta}: \quad {\rm i} (1-r^{2})(A_{0}+A_{-2}) + r^{2}C
  \end{equation}
  \begin{equation}
    \label{eqa:4d}
    {\rm Constant}: \quad {\rm i} (1-r^{2})(A_{0}+A_{-2}) + C 
  \end{equation}
\end{subequations}
where
\begin{equation*}
  C = \frac{{\rm i}}{2}(A_{-2}-A_{0}) + {\rm i} C_{a} = (c_{0} + \overline{a}_{0}) + {\rm i} C_{a}
\end{equation*}
should be a pure imaginary constant due to Eq.~(\ref{eqa:3}). All the coefficients in Eq.~(\ref{eqa:4}) are multiplied by a constant $\frac{{\rm i}}{2}$ for simplicity. It can be verified that the coefficients in Eq.~(\ref{eqa:4}) are respectively identical to the those in Eqs.~(29), (28), (31), and (32) in Ref~\cite{Verruijt1997traction}.

\clearpage
\section{Coefficient expansion in Eq.~(\ref{eq:4.15b})}
\label{sec:B}

The following expansions are neccessary:
\begin{subequations}
  \label{eqb:1}
  \begin{equation}
    \label{eqb:1a}
    \frac{1}{1-r{\rm e}^{{\rm i}\theta}} = \sum\limits_{k=0}^{\infty} r^{k}{\rm e}^{{\rm i}k\theta}
  \end{equation}
  \begin{equation}
    \label{eqb:1b}
    \ln(1-r{\rm e}^{{\rm i}\theta}) = -\sum\limits_{k=1}^{\infty} \frac{r^{k}}{k}{\rm e}^{{\rm i}k\theta}
  \end{equation}
\end{subequations}

\begin{subequations}
  \label{eqb:2}
  \begin{equation}
    \label{eqb:2a}
    \frac{1}{{\rm e}^{{\rm i}\theta}-r} = \sum\limits_{k=1}^{\infty} r^{k-1}{\rm e}^{{-{\rm i}k\theta}}
  \end{equation}
  \begin{equation}
    \label{eqb:2b}
    \frac{1}{({\rm e}^{{\rm i}\theta}-r)^{2}} = \sum\limits_{k=1}^{\infty} (k-1) r^{k-2}{\rm e}^{-{\rm i}k\theta}
  \end{equation}
  \begin{equation}
    \label{eqb:2c}
    \ln(1-r{\rm e}^{-{\rm i}\theta}) = -\sum\limits_{k=1}^{\infty} \frac{r^{k}}{k}{\rm e}^{-{\rm i}k\theta}
  \end{equation}
\end{subequations}

\begin{equation}
  \label{eqb:3}
  \frac{1}{(1-r{\rm e}^{{\rm i}\theta})({\rm e}^{{\rm i}\theta}-r)^{2}} = \frac{K_{1}}{1-r{\rm e}^{{\rm i}\theta}} + \frac{K_{2}{\rm e}^{{\rm i}\theta}+K_{3}}{({\rm e}^{{\rm i}\theta}-r)^{2}}
\end{equation}
where
\begin{equation*}
  \left\{
    \begin{aligned}
      K_{1} = & \; \frac{r^{2}}{(1-r^{2})^{2}} \\
      K_{2} = & \; \frac{r}{(1-r^{2})^{2}} \\
      K_{3} = & \; \frac{1-2r^{2}}{(1-r^{2})^{2}} \\
    \end{aligned}
  \right.
\end{equation*}
Then with substitution of Eqs.~(\ref{eqb:1a}) and (\ref{eqb:2b}), Eq.~(\ref{eqb:3}) can be expanded as
\begin{equation}
  \label{eqb:3'}
  \tag{B.3'}
  \frac{1}{(1-r{\rm e}^{{\rm i}\theta})({\rm e}^{{\rm i}\theta}-r)^{2}} = \sum\limits_{k=0}^{\infty} f_{k} {\rm e}^{{\rm i}k\theta} + \sum\limits_{k=1}^{\infty} g_{k} {\rm e}^{-{\rm i}k\theta}
\end{equation}
where
\begin{equation*}
  \left\{
    \begin{aligned}
      f_{k} = & \; K_{1}r^{k} \\
      g_{k} = & \; K_{2}\cdot kr^{k-1}+K_{3}\cdot (k-1)r^{k-2}
    \end{aligned}
  \right.
\end{equation*}

The following notation is used:
\begin{equation}
  \label{eqb:4}
  \left\{
    \begin{aligned}
      L_{1} = & \; - k_{0}\gamma a^{2}(1-r^{2})^{2} \\
      L_{2} = & \; - \gamma a^{2} \\
      L_{3} = & \; \gamma a^{2} r^{2} \\
      L_{4} = & \; \gamma a R \\
      L_{5} = & \; -\gamma a r R \\
      L_{6} = & \; \gamma R^{2} \\
    \end{aligned}
  \right.
\end{equation}

With Eqs.~(\ref{eqb:1}), (\ref{eqb:2}), and (\ref{eqb:3'}), $E_{k}$ in Eq.~(4.15b) can be expressed as
\begin{equation}
  \label{eqb:5}
  \left\{
    \begin{aligned}
      E_{k} = & \; L_{1}f_{k-2} + L_{2}r^{k} - L_{6}\frac{r^{k}}{k} + L_{6}\frac{r^{k}}{k-1}, \quad k \geq 2 \\
      E_{1} = & \; L_{1}g_{1} + L_{2}r - L_{6}r \\
      E_{0} = & \; L_{1}g_{2} + L_{2} - L_{4}r + L_{5} - L_{6}r^{2} \\
      E_{-k} = & \; L_{1}g_{k+2} + L_{3}(k-1)r^{k-2} - L_{3}kr^{k} + L_{4}r^{k-1} - L_{4}r^{k+1} + L_{6}\frac{r^{k}}{k} - L_{6}\frac{r^{k+2}}{k+1}, \quad k \geq 1 \\
    \end{aligned}
  \right.
\end{equation}

\clearpage

\clearpage

\clearpage

\clearpage

\clearpage

\clearpage

\clearpage

\clearpage

\clearpage

\clearpage

\clearpage

\clearpage

\clearpage

\clearpage

\clearpage
\bibliographystyle{plainnat}
\bibliography{biblio}  

\end{document}